\documentclass[english,12pt,twoside,reqno,letter]{amsart}
\usepackage[square,comma]{natbib}
\usepackage[english,german]{babel} 
\usepackage{exscale}   
\usepackage{amsmath}
\usepackage{amsfonts}
\usepackage{stmaryrd}
\usepackage{xspace}
\usepackage{mathrsfs}

\newcommand{\Cset}{\mathbb{C}} 
\newcommand{\Nset}{\mathbb{N}}
 
\newcommand{\Rset}{\mathbb{R}}
\newcommand{\Zset}{\mathbb{Z}} 
\renewcommand{\d}{\hskip 1pt\mathrm{d}}
\newcommand{\e}{\mathrm{e}}
\newcommand{\qedopt}{}
\newcommand{\qedopteqn}{\tag*{\qedhere}} 
\renewcommand{\theequation}{\thesection.\arabic{equation}}
\numberwithin{equation}{subsection}
\newtheorem{Thm}{Theorem}[subsection]
\newtheorem{Lem}[Thm]{Lemma}
\newtheorem{Cor}[Thm]{Corollary}
\newtheorem{Prop}[Thm]{Proposition}
\newtheorem{Defn}[Thm]{Definition}
\newtheorem{Rem}[Thm]{Remark}
\newtheorem{Exmp}[Thm]{Example}
\newtheorem{Notation}[Thm]{Notation}
\newcommand{\spc}{~}
\newtheorem{ThmApp}{Theorem}[section]
\newtheorem{LemApp}[ThmApp]{Lemma}
\newtheorem{CorApp}[ThmApp]{Corollary}

\newtheorem{DefnApp}[ThmApp]{Definition}
\renewcommand{\theenumi}{\roman{enumi}}     
\newcommand{\cA}{\ensuremath{{\mathcal A}}\xspace}         
\newcommand{\cB}{\ensuremath{{\mathcal B}}\xspace}         
\newcommand{\cC}{\ensuremath{{\mathcal C}}\xspace}         
\newcommand{\cE}{\ensuremath{{\mathcal E}}\xspace}         
\newcommand{\cF}{\ensuremath{{\mathcal F}}\xspace}         
\newcommand{\cH}{\ensuremath{{\mathcal H}}\xspace}         
\newcommand{\cK}{\ensuremath{{\mathcal K}}\xspace}         
\newcommand{\cI}{\ensuremath{{\mathcal I}}\xspace}         
\newcommand{\cL}{\ensuremath{{\mathcal L}}\xspace}         
\newcommand{\cN}{\ensuremath{{\mathcal N}}\xspace}         
\newcommand{\cP}{\ensuremath{{\mathcal P}}\xspace}         
\newcommand{\cT}{\ensuremath{{\mathcal T}}\xspace}         
\newcommand{\cU}{\ensuremath{{\mathcal U}}\xspace}         
\newcommand{\cV}{\ensuremath{{\mathcal V}}\xspace}         
\newcommand{\cZ}{\ensuremath{{\mathcal Z}}\xspace}         
\newcommand{\al}{\ensuremath{\alpha}\xspace}               
\newcommand{\be}{\ensuremath{\beta}\xspace}                
\newcommand{\Ga}{\ensuremath{\Gamma}\xspace}               
\newcommand{\de}{\ensuremath{\delta}\xspace}               
\newcommand{\De}{\ensuremath{\Delta}\xspace}               
\newcommand{\la}{\ensuremath{\lambda}\xspace}              
\newcommand{\La}{\ensuremath{\Lambda}\xspace}              
\newcommand{\si}{\ensuremath{\sigma}\xspace}               
\newcommand{\Si}{\ensuremath{\Sigma}\xspace}               
\newcommand{\ve}{\varepsilon}                              
\newcommand{\om}{\ensuremath{\omega}\xspace}               
\newcommand{\Om}{\ensuremath{\Omega}\xspace}               
\renewcommand{\phi}{\ensuremath{\varphi}\xspace}           
\renewcommand{\rho}{\ensuremath{\varrho}\xspace}           
\newcommand{\1}{\ensuremath{{\rm 1\kern-.25em l}}\xspace}  
\newcommand{\Exp}{\ensuremath{\operatorname{Exp}}\xspace}  
\renewcommand{\i}{\mathrm{i}}                              
\newcommand{\id}{\operatorname{id}}                        
\newcommand{\IM}{\operatorname{Im}}                        
\newcommand{\Ln}{\ensuremath{\operatorname{Ln}}\xspace}    
\newcommand{\Lno}{\ensuremath{\operatorname{Ln}_0}\xspace} 
\newcommand{\RE}{\operatorname{Re}}                        
\newcommand{\Aut}{\ensuremath{\operatorname{Aut}}\xspace}  
\newcommand{\Int}[1]{\operatorname{Int}#1}         
\newcommand{\Mor}{\ensuremath{\operatorname{Mor}}\xspace}  
\newcommand{\set}[2]{\mathopen{\{}#1\mathop{|}#2\mathclose{\}}}
\newcommand{\Ao}{\ensuremath{\cA_{0}}\xspace}           
\newcommand{\Ho}{\ensuremath{\cH_{0}}\xspace}           
\newcommand{\srS}{\mathscr{S}\xspace}                      
\newcommand{\Apsi}{\ensuremath{(\cA,\psi)}\xspace}      
\newcommand{\Lapneu}{\ensuremath{\cE}\xspace}           
\newcommand{\LapIneu}{\ensuremath{\cE_I}\xspace}        
\newcommand{\LapJneu}{\ensuremath{\cE_J}\xspace}        
\newcommand{\LapXneu}[1]{\ensuremath{\cE_{#1}}\xspace}  
\newcommand{\CBS}{\ensuremath{%
            (\cA,\psi,S,(\cA_I)_{I\in \cI})}\xspace}    
\newcommand{\CBSE}{\ensuremath{%
            (\cE,\1,S,(\cE_I)_{I\in \cI})}\xspace}      
\newcommand{\bimod}[1]{{\phantom{_{#1}}\llap{$_{#1}$}\cE_{\!\Ao}}}
\newcommand{\Cm}{\mathscr{C}_{0}^{}}        
\newcommand{\Cmo}{\mathscr{C}_{0}^{0}}      
\newcommand{\Ca}{\mathscr{C}_{0}^{}}        
\newcommand{\Cao}{\mathscr{C}_{0}^{0}}      
\newcommand{\cac}{{c}}                   
\newcommand{\Eo}{E_{0}\xspace}           
\newcommand{\ncac}{{b}}                  
\newcommand{\tii}{t_i}                              
\newcommand{\tiemti}{t_{i+1}-t_i}                   
\newcommand{\sij}{s_{i,j}}                          
\newcommand{\sije}{s_{i,j+1}}                       
\newcommand{\sijemij}{s_{i,j+1}-s_{i,j}}            
\newcommand{\nmolim}[1]{\nmot\text{\rm -}\lim_{#1}}          
\newcommand{\stp}{\ensuremath{\operatorname{stop}}\xspace}   
\newcommand{\skp}[2]{\mathopen{\langle}%
                 #1\mathop{|}#2\mathclose{\rangle}} 
\newcommand{\bskp}[2]{\mathopen{\bigl\langle}%
                 #1\mathop{\big|}#2\mathclose{\bigr\rangle}} 
\newcommand{\Bskp}[2]{\mathopen{\Bigl\langle}%
                 #1\mathop{\Big|}#2\mathclose{\Bigr\rangle}} 
\newcommand{\skpo}[2]{\skp{#1}{#2}_{0}^{}}          
\newcommand{\skpO}[2]{\skp{#1}{#2}_{0}}             
\newcommand{\bskpo}[2]{\bskp{#1}{#2}_{0}^{}}        
\newcommand{\Bskpo}[2]{\Bskp{#1}{#2}_{0}^{}}        
\newcommand{\nm}[1]{\mathord{\mathopen{\|}#1\mathclose{\|}}} 
\newcommand{\bnm}[1]{\mathord{\mathopen{\bigl\|}#1%
                              \mathclose{\bigr\|}}} 
\newcommand{\Bnm}[1]{\mathord{\mathopen{\Bigl\|}#1%
                              \mathclose{\Bigr\|}}} 
\newcommand{\nmo}[1]{\nm{#1}_{0}^{}}                
\newcommand{\nmO}[1]{\nm{#1}_{0}}                   
\newcommand{\Bnmo}[1]{\Bnm{#1}_{0}^{}}              
\newcommand{\BnmO}[1]{\Bnm{#1}_{0}}                 
\newcommand{\ab}[1]{\mathord{\lvert#1\rvert}}       
\newcommand{\Bab}[1]{\mathord{\Bigl\lvert#1\Bigr\rvert}} 
\newcommand{\abo}[1]{\ab{#1}_{0}^{}}                
\newcommand{\abO}[1]{\ab{#1}_{0}}                   
\newcommand{\Babo}[1]{\Bab{#1}_{0}^{}}              
\newcommand{\BabO}[1]{\Bab{#1}_{0}}                 
\newcommand{\nmt}{\raisebox{1.5pt}{${\scriptstyle\|\,\|}$}} 
\newcommand{\nmot}{\raisebox{1.5pt}{${\scriptstyle\|\,\|%
                   \raisebox{-1pt}{$_{0}$}}$}}              

\begin{document}  

\selectlanguage{english}

\begin{center}
{\large \textbf{NON-COMMUTATIVE\\ 
         CONTINUOUS BERNOULLI SHIFTS}}\\[24pt]
\textsc{J{\"u}rgen Hellmich}\\[3pt]
{\small \textit{Mathematisches Institut, 
              Universit\"{a}t T\"{u}bingen,\\[0pt] 
    Auf der Morgenstelle 10, D-72076 T\"{u}bingen,Germany}\\[0pt]
    (e-mail:\texttt{juergen.hellmich@uni-tuebingen.de})}\\[12pt] 

\textsc{Claus K{\"o}stler}\\[3pt]
{\small  \textit{Department of Mathematics and Statistics, 
                  Queen's University, \\[0pt] 
             Kingston, Ontario K7L 3N6, Canada}\\[0pt]
             (e-mail:\texttt{koestler@mast.queensu.ca})}\\[12pt]

\textsc{Burkhard K\"{u}mmerer}\\[3pt] 
{\small \textit{Fachbereich Mathematik, 
        Technische Universit\"{a}t Darmstadt,\\[0pt]
        Schlo{\ss}gartenstr.~7, D-64289 Darmstadt, Germany}\\[0pt]
     (e-mail:\texttt{kuemmerer@mathematik.tu-darmstadt.de})}\\[12pt]
\end{center}
\vspace{12pt}
\begin{center}
November 15th, 2004 
\end{center}
\vspace{1cm}
\noindent
\textbf{Abstract:} We introduce a non-commutative extension of
Tsirelson-Vershik's noises \cite{TsVe98a,Tsir04a}, called
\emph{(non-commu\-ta\-tive) continuous Bernoulli shifts}. These shifts
encode stochastic independence in terms of commuting squares, as they
are familiar in subfactor theory \cite{Popa83a,GHJ89a}. Such shifts
are, in particular, capable of producing Arveson's product system of
type~I and type~II \cite{Arve03a}. We investigate the structure of
these shifts and prove that the von Neumann algebra of a
(scalar-expected) continuous Bernoulli shift is either finite or of
type~III.
  
The role of (`classical') $G$-stationary flows for Tsirelson-Vershik's 
noises is now played by \emph{cocycles of continuous Bernoulli
shifts}. We show that these cocycles provide an operator algebraic 
notion for L\'evy processes. They lead, in particular, to units and 
`logarithms' of units in Arveson's product systems \cite{Koes04b}. 
Furthermore, we introduce \emph{(non-commutative)
white noises}, which are operator algebraic versions of Tsirelson's
`classical' noises. We give examples coming from probability, quantum
probability and from Voiculescu's theory of free probability
\cite{VDN92a}.
   
Our main result is a bijective correspondence between additive and
unital shift cocycles. For the proof of the correspondence we develop
tools which are of interest on their own: non-commutative extensions of
stochastic It\^o integration, stochastic logarithms and exponentials.
 
\pagebreak 

{%
\def\widedotfill{\leaders\hbox to 10pt{\hfil.\hfil}\hfill}
\def\pg#1{\widedotfill\rlap{\hbox to 15pt{\hfill{\small#1}}}\par}
\rightskip=15pt\leftskip=10pt%
\newcommand{\sct}[2]{\noindent\llap{\hbox to%
    10pt{\textbf{#1}\hfill}}~\mbox{#2~}}
\newcommand{\separ}{\vspace{.2em}}

\section*{Contents}
\sct{}{\textbf{Introduction}}%
\pg{\textbf{\pageref{section:introduction}}}%
\separ
\sct{1}{\textbf{Preliminaries}}%
\pg{\textbf{\pageref{section:preliminaries}}}
\sct{}{1.1}{\small General terminology}%
\pg{\pageref{subsection:prelim-general}}
\sct{}{1.2}{\small Non-commutative probability 
            spaces and their morphisms}\hfill%
\pg{\pageref{subsection:spaces-morphisms}}
\sct{}{1.3}{\small Filtrations}%
\pg{\pageref{subsection:filtrations}}
\sct{2}{\textbf{Non-commutative independence}}%
\pg{\textbf{\pageref{section:independence}}}
\sct{}{2.1}{\small $\Ao$-independence and Popa's commuting squares}%
\pg{\pageref{subsection:independence-Popa}}
\sct{}{2.2}{\small Commuting subalgebras and $\Cset$-independence}%
\pg{\pageref{subsection:independence-examples}}
\sct{}{2.3}{\small From $\Cset$-independence to $\Ao$-independence}%
\pg{\pageref{subsection:C2A0-independence}}
\sct{}{2.4}{\small Non-commutative examples of $\Ao$-independence}%
\pg{\pageref{subsection:nc-independence-examples}}
\separ
\sct{3}{\textbf{Continuous Bernoulli shifts I}}%
\pg{\textbf{\pageref{section:CBSI}}}
\sct{}{3.1}{\small Continuous Bernoulli shifts and their basic 
properties}\pg{\pageref{subsection:cbs-basic}}
\sct{}{3.2}{\small The type of a $\Cset$-expected continuous 
Bernoulli shift}\pg{\pageref{subsection:type}}
\sct{}{3.3}{\small Composition of continuous Bernoulli shifts}%
\pg{\pageref{subsection:composition}}
\sct{}{3.4}{\small Decomposition of continuous Bernoulli shifts}%
\pg{\pageref{subsection:decomposition}}
\separ
\sct{4}{\textbf{Continuous Bernoulli shifts II}}%
\pg{\textbf{\pageref{section:CBSII}}}
\sct{}{4.1}{\small Local minimality and local maximality}%
\pg{\pageref{subsection:local-min-max}}
\sct{}{4.2}{\small Enriched independence}%
\pg{\pageref{subsection:enriched-independence}}
\sct{}{4.3}{\small Commuting past and future}%
\pg{\pageref{subsection:commuting-past/future}}
\sct{}{4.4}{\small Commutative von Neumann algebras}%
\pg{\pageref{subsection:commutative-BS}}
\sct{}{4.5}{\small Local minimality and compressions}%
\pg{\pageref{subsection:local-min-comp}}
\sct{}{4.6}{\small Examples from probability theory}%
\pg{\pageref{subsection:CBS-examples}}
\sct{}{4.7}{\small Examples from quantum probability theory}%
\pg{\pageref{subsection:nCBS-examples}}
\separ
\sct{5}{\textbf{Continuous GNS Bernoulli shifts}}%
\pg{\textbf{\pageref{section:CGNSBS}}}
\sct{}{5.1}{\small Hilbert bimodules of 
            $\Ao$-expected probability spaces}%
\pg{\pageref{subsection:hm-spaces}}
\sct{}{5.2}{\small GNS representation of morphisms}%
\pg{\pageref{subsection:hm-morphisms}}
\sct{}{5.3}{\small The product of $\Ao$-independent elements}%
\pg{\pageref{subsection:hm-mult}}
\sct{}{5.4}{\small Continuous GNS Bernoulli shifts}%
\pg{\pageref{subsection:CBS-GNS}}
\separ
\sct{6}{\textbf{Cocycles of continuous (GNS) Bernoulli shifts}}%
\pg{\textbf{\pageref{section:cocycles}}}
\sct{}{6.1}{\small Multiplicative cocycles of 
            continuous Bernoulli shifts}%
\pg{\pageref{subsection:uc-cbs}}
\sct{}{6.2}{\small Multiplicative cocycles of continuous 
            GNS Bernoulli shifts}%
\pg{\pageref{subsection:uc-cgnsbs}}%
\sct{}{6.3}{\small Additive cocycles of continuous 
            GNS Bernoulli shifts}%
\pg{\pageref{subsection:ac-cgnsbs}}
\sct{}{6.4}{\small The correspondence}%
\pg{\pageref{subsection:correspondence}}
\sct{}{6.5}{\small Non-commutative white noises}%
\pg{\pageref{subsection:CBS-WN}}
\sct{}{6.6}{\small Examples for the correspondence}%
\pg{\pageref{subsection:examples-correspondence}}
\separ
\sct{7}{\textbf{Non-commutative It{\^o} integration}}%
\pg{\textbf{\pageref{section:nc-ito-integration}}}
\sct{}{7.1}{\small Non-commutative It{\^o} integrals for 
            simple adapted processes}%
\pg{\pageref{subsection:ito-simple}}
\sct{}{7.2}{\small An extension of the non-commutative 
            It{\^o} integral}%
\pg{\pageref{subsection:ito-extension}}
\sct{}{7.3}{\small Non-commutative It\^o differential equations}%
\pg{\pageref{subsection:IDE}}
\separ
\pagebreak
\sct{8}{\textbf{Non-commutative exponentials and logarithms}}%
\pg{\textbf{\pageref{section:nc-ln-exp}}}
\sct{}{8.1}{\small Non-commutative exponentials of additive cocycles}%
\pg{\pageref{subsection:nc-exp}}
\sct{}{8.2}{\small Non-commutative logarithms of unital cocycles}%
\pg{\pageref{subsection:nc-ln}}
\sct{}{8.3}{\small Proof of the correspondence}%
\pg{\pageref{subsection:proof-correspondence}}
\separ
\sct{}{\textbf{Appendix A: Hilbert W*-modules}}%
\pg{\textbf{\pageref{sec-appendix-hm}}}
\sct{}{\textbf{Appendix B: The $\boldsymbol{\psi}$-adjoint 
            of morphisms}}%
\pg{\textbf{\pageref{sec-appendix-psiadjoint}}}
\sct{}{\textbf{References}}%
\pg{\textbf{\pageref{section:bibliography}}}
}

\allowdisplaybreaks 
\section*{Introduction}                    \label{section:introduction}

Recently, Tsirelson and Vershik established the existence of
intrinsically non-linear random fields \cite{TsVe98a}. Their
surprising result was stimulated by the existence of Arveson-Power's
product systems of non-type~I \cite{Powe87a,Arve89a,Arve03a}. These
random fields or `noises' go beyond the realm of the L\'evy-Khintchine
formula and provide, in particular, a rich probabilistic source of
non-type~I product systems \cite{Tsir03a,Tsir04a,Liebscher1}. Here we
are interested in Tsirelson-Vershik's noises, as they are defined in
\cite{Tsir98a,Tsir04a}. We will introduce a non-commutative extension
of these noises, called \emph{(non-commutative) continuous Bernoulli
shifts}. These shifts incorporate so-called \emph{$\Ao$-independence}
which extends amalgamated stochastic independence to an operator
algebraic frame, known as commuting squares in subfactor theory
\cite{Popa83a,GHJ89a,JoSu97a}. These shifts may be regarded as
two-sided `time-continuous' analogues of shifts on towers of von
Neumann algebras \cite{Rupp95a,GoKo04a} (as they are also implicitly
present in \cite{JoSu97a}, for example). But the notion of
$\Ao$-independence comprises also Voiculescu's
amalgamated free independence \cite{VDN92a}. Consequently, our
approach is, in particular, in close contact with free probability
theory. Further examples of continuous Bernoulli shifts arise, aside
of fermionic and bosonic noises, on deformed Fock spaces 
\cite{BoSp91a,BoSp94a,BoGu02a,BKS97a,GuMa02a}.   

Non-commutative continuous Bernoulli shifts comprise, in an algebraic
form, all 'classical' examples of Tsirelson-Vershik's
noises. 'Classical noises' are generated by additive (square
integrable adapted) stationary flows, called L\'evy processes,
and are classified via the L\'evy-Khinchin formula (see e.g.\ 
\cite[Corollary 6a7]{Tsir04a}). 

What is the operator algebraic 
counterpart of a 'classical noise'? Here, the situation is much more complex
and we are only at the beginning of understanding this complexity (see also
\cite{KoSp04a}). Let us illustrate this in the case $\cA_0= \Cset$: 
Brownian motion is unique (up to stochastic equivalence), but there 
exist many different non-commutative Brownian motions: 
$q$-Brownian motions $(-1 < q <1)$ (including free Brownian motion), 
bosonic and fermionic Brownian motions (parametrized by 
'temperature' or, equivalently, the period of the modular 
automorphism group), and this is just the 
beginning of a long list. But all these diverse examples have in
common that they appear as additive (adapted) cocycles in the GNS 
Hilbert space of a ($\Cset$-expected) continuous Bernoulli shift. 
We will show in this paper that such additive cocycles are in 
correspondence to multiplicative (adapted) cocycles in the GNS 
Hilbert space, called unital cocycles. This will imply, in particular,
that unitary (adapted) cocycles, i.e.\ multiplicative cocycles in 
the unitary operators of the von Neumann algebra of a continuous 
Bernoulli shift, are also in correspondence with additive cocycles. 
This means in the terminology of Tsirelson and Vershik that unitary
cocycles are 'linearizable'. Thus we have available an 
operator algebraic notion of 'classical noise': a continuous 
Bernoulli shift is called a (non-commutative) white noise if it is
generated by the set of all unitary cocycles.  

Why do we avoid to say `A continuous Bernoulli shift is called a
white noise if it is generated by the set of all additive cocycles'? This 
would become conceptually cumbersome already for bosonic white 
noises in Araki-Woods representations: the vector space of 
additive cocycles does not capture the type of the von Neumann 
algebra (see also Example \ref{CCR-correspondence}). In all  
examples this kind of information is normally encoded into the choice 
of a functor or into mixed higher moments which we do not have
available in our general setting.  

Continuous Bernoulli shifts comprise also all 'non-classical noises'.
But we do yet not know a single example of a 'non-classical quantum 
noise' which is a non-commutative continuous Bernoulli shift, even 
though the latter object is a straightforward extension of 
Tsirelson-Vershik's noises. Such an example would produce, similarly 
as 'non-classical noises' do \cite{Tsir04a}, an Arveson product 
system of Hilbert spaces of type $II$ \cite{Koes04a}.  

Some clarifying remarks on the terminology are appropriate at this
point. The attribute `non-commutative' will always be used the sense of
`not necessarily commutative'. Frequently, we will drop it at all and
will just write, for example, `continuous Bernoulli shift' instead of
`non-commutative continuous Bernoulli shift'. The attribute `quantum'
will be reserved for situations beyond the realm of probability theory.
We will avoid the attribute `classical' and use instead `commutative'.
This convention is motivated from a conflict with the terminology in
\cite{Tsir04a}. A `classical noise' there will be a `commutative white
noise' here, a `white noise' therein will be a `Gaussian white noise'
herein.

We have put emphasis on a self-contained, comprehensive presentation
of the subject, since our approach is not available easily in the
published literature. Moreover, we refrained here to go for the most
general cases. This is mainly a technical matter and can be
accomplished later on. From doing so we hope that our work is
better accessible for the reader coming from probability, quantum
probability, quantum dynamics, quantum symmetries or operator
algebras, being interested in connecting these fields.

A starting point of this work has been the operator algebraic approach
to stationary quantum Markov processes in \cite{Kuem85a} and results in
\cite{Kuem84aUP,Prin89a,Kuem93UP,Rupp95a,Koes00a,Hell01a}. 
The present conceptual form is also stimulated in parts by
\cite{Arve03a,Tsir04a}. In the context of product systems and their
classification are, in particular, of relevance 
\cite{Bhat99a,Bhat01a,MuSo02a,Liebscher1,BBLS04a,BhSr04a}.  
Finally, we want to bring to the
reader's attention also the work of Gohm \cite{Gohm04a,Gohm04bTA}, 
which is in close contact to our approach. Further information on
related approaches is provided within each Section. 

Unital cocycles of continuous Bernoulli shifts are in closest contact
to units of product systems of Hilbert spaces or modules. More details
about this relation will be provided in \cite{Koes04b}.

Unitary cocycles of continuous Bernoulli shifts or white noises give
immediately rise to Markovian cocycles or stationary quantum 
Markov processes, as they are relevant in quantum dynamics and 
quantum probability. There are various approaches to the construction of
quantum Markov processes in the literature. An operator algebraic 
setting is used in \cite{AFL,Kuem02a,Arve03a,Koes03a}. For accounts 
on bosonic Fock space and Hudson-Parthasarathy's quantum stochastic 
calculus we refer the reader to \cite{Part92a,Meye93a}. Meanwhile, 
this approach is further developed and a modern account can be found 
in \cite{Lind04a}.    

Next we will provide an outline of the major results and contents of all
sections. More detailed introductions are contained at the beginning of
each section.

\subsection*{Section \ref{section:preliminaries}:}
Throughout a non-commutative probability space $(\cA,\psi)$ is modeled
by a von Neumann algebra $\cA$ together with a faithful normal state
$\psi$ on $\cA$. The predual $\cA_*$ is always assumed to be separable.
If $\Ao$ is a subalgebra of $\cA$, such that the conditional
expectation from $\cA$ onto $\Ao$ exists, then we will say that $(\cA,
\psi$) is an $\Ao$-expected non-commutative probability space (see
Subsection \ref{subsection:spaces-morphisms} for a motivation).

\subsection*{Section  \ref{section:independence}:}
Our non-commutative extension of (amalgamated) stochastic independence
is connected intimately to \emph{commuting squares}, as they have been
introduced by Popa in subfactor theory \cite{Popa83a,GHJ89a,Popa90a}.
Let $(\cA,\psi)$ be an $\Ao$-expected probability space and denote by
$E_0$ the conditional expectation from $\cA$ onto $\Ao$.

\begin{Defn}
  Let $\cB$ and $\cC$ be two von Neumann subalgebras of $\cA$ such
  that, respectively, the conditional expectations $E_\cB$, $E_\cC$
  from $\cA$ onto $\cB$, $\cC$ exist. Then $\cB$ and $\cC$ are called
  \emph{$\Ao$-independent} if $\cB \cap \cC =\Ao$ and $E_\cB E_\cA
  = E_0$.
\end{Defn}

In other words, the four von Neumann algebras $\cA, \Ao, \cB$ and
$\cC$ form a commuting square:
\begin{align*}
  \begin{matrix}
    \cC           & \subset & \cA                               \\
    \cup          &         & \cup                              \\
    \Ao           & \subset & \cB
  \end{matrix}
\end{align*}
If the von Neumann algebra $\cA$ is commutative and $\Ao \simeq
\Cset$, then one recovers the usual notion of stochastic independence
in probability theory. In the general setting, $\Ao$-independence
encloses amalgamated stochastic independence, tensor product
independence and Voiculescu's amalgamated free independence
\cite{VDN92a}. But most importantly, it is not restricted to
non-commutative notions of stochastic independence with  universal
product rules \cite{Spei97a,BeSc02a}. Further examples of
$\Ao$-independence are accessible by `white noise functors'
\cite{Kuem96a,BKS97a,GuMa02a}, applied to von Neumann algebras
generated by `generalized Brownian motions'
\cite{BoSp91a,BoSp94a}. Moreover, we expect that Anshelevich's
$q$-L\'evy processes will provide further examples of
$\Ao$-independence \cite{Ansh04aTA}.
 
\subsection*{Section \ref{section:CBSI}:}
We will introduce \emph{non-commutative continuous  Bernoulli shifts}
and study their structure. These shifts provide a non-commutative
extension of Tsirelson-Vershik's noises, or speaking more technically,
of  homogeneous continuous products of probability spaces
\cite{TsVe98a,Tsir04a}. Their infrastructure, as stated in Definition
\ref{def:cbs}, is integral for this paper,  so that we will 
introduce them informally next.

An $\Ao$-expected {\em non-commutative continuous Bernoulli shift}
encodes from an algebraic point of view the following structure.
Consider an $\Ao$-expected probability space $(\cA,\psi)$ together with
an pointwise weakly*-continuous automorphism group $(S_t)_{t\in \Rset}$
(called \emph{shift}) and a family of von Neumann subalgebras
$(\cA_{[r,s]})_{-\infty \le r \le s \le \infty}$ (called
\emph{filtration}), indexed by (`time'-)intervals, such that
$S_{t}\cA_{[r,s]}=\cA_{[r+t,s+t]}$. Assume that all conditional
expectations from $\cA$ onto $\cA_{[r,s]}$ ($r\le s$) exist and that
$\psi$ is $S$-invariant. Now we encode $\Ao$-independence by the
requirement that the filtration forms a family of \emph{commuting
squares} which is moreover shifted covariantly in `time' by the
action of $S$:
\begin{align*}
  \begin{pmatrix}
    \cA_{[r,s]}     & \subset & \cA_{[r,v]}                       \\
    \cup            & \quad   & \cup                              \\
    \Ao             & \subset & \cA_{[u,v]}
  \end{pmatrix}_{r < s <u < v}
  \hskip-8pt
  \xrightarrow{\phantom{x}\overline{S}_{t}\phantom{x}}
  \begin{pmatrix}
    \cA_{[r+t,s+t]} & \subset & \cA_{[r+t,v+t]}                   \\
    \cup            & \quad   & \cup                              \\
    \Ao             & \subset & \cA_{[u+t,v+t]}
  \end{pmatrix}_{r < s < u < v}
\end{align*}
The von Neumann algebra $\Ao$ is not shifted. But we will stipulate a
stronger condition: $\Ao$ is required to be the fixed point algebra of
the shift $S$. Moreover, we will impose the system to minimality: the
family $(\cA_{[r,s]})_{-\infty < r \le s <\infty}$ generates already
$\cA$ (`minimal filtration'). A priori we do not require `local
minimality', $\cA_{[r,s]}\vee \cA_{[s,t]} \neq \cA_{[r,t]}$ may occur.
Denoting by the set of all (closed) intervals by $\cI$, an object
\[
      \big(\cA,\psi, S, (\cA_{I})_{I \in \cI}\big),
\]  
enjoying all these properties, will be called an \emph{$\Ao$-expected
(non-commutative) continuous Bernoulli shift} (see Definition
\ref{def:cbs}, where we will also allow non-closed intervals in the
index set). Notice that we do not assume continuity properties for the 
filtration itself. Assuming for a moment the `local minimality'
condition $\cA_{[r,s]}\vee \cA_{[s,t]} =\cA_{[r,t]}$, we observe the
following.

\begin{enumerate}
\item[(i)] Putting $\Ao\simeq \Cset$ and requiring that $\cA$ is
  commutative, one obtains an algebraic version of Tsirelson-Vershik's
  noises (see Subsection \ref{subsection:commutative-BS}).
\item[(ii)] Putting as `time' $\Zset$ (or $\Nset_0$) instead of
  $\Rset$, dropping all notions of continuity, above scheme reduces to
  a (one-sided) Bernoulli shift which finds its examples in subfactor
  theory \cite{Rupp95a,KuMa98a,GoKo04a}, in particular.
\end{enumerate}
Here we will focus onto the case of continuous `time', which includes
Tsirelson-Vershik's noises. The `discrete time' case and its connection
to subfactor theory goes beyond the limits of this work and is
postponed to sequel publications.

In Section \ref{section:CBSI} we will concentrate on properties of
$\Ao$-expected continuous Bernoulli shifts which follow
out of Definition \ref{def:cbs} without any further assumptions. Among
these properties are strongly mixing properties of the shift,
stability with respect to compositions (like tensor products and
direct sums) and decompositions (like compressions with conditional
expectations). Most importantly, the structure of a $\Cset$-expected
non-commutative continuous Bernoulli shift $\CBS$ is already
sufficient to give results on the type of its von Neumann algebra (see
Subsection \ref{subsection:type}):

\begin{Thm} 
  Let $\CBS$ be a $\Cset$-expected continuous Bernoulli shift. Then
  $\cA$ is either finite or of type~III. The state $\psi$ is
  non-tracial if and only if $\cA$ is of type~III.
\end{Thm}

This result extends to $\Ao$-expected continuous Bernoulli shifts if
one passes to `derived' continuous Bernoulli shifts, similarly as it is
understood for towers of von Neumann algebras in subfactor theory. We
will present this line of research in more depth in \cite{HeKo04}.

\subsection*{Section \ref{section:CBSII}:}
Up to the present, all examples of continuous Bernoulli shifts enjoy
much more algebraic structure and continuity properties, as they are
stipulated in Definition \ref{def:cbs}. Among such additional
(algebraic) properties, which we will study in this section, are
\emph{local minimality}, \emph{local maximality}, \emph{enriched
$\Ao$-independence}, \emph{commuting past and future}. These are
properties which are quite familiar in Arveson's approach to quantum
dynamics \cite{Arve03a}. At this place we will also investigate the
relationship of continuous Bernoulli shifts and Tsirelson-Vershik's
noises (see Subsection \ref{subsection:commutative-BS}). Finally, we
will give examples of continuous Bernoulli shifts, coming both from
probability theory and quantum probability theory. Among the first
ones are Gaussian, Poisson white noise and Tsirelson-Vershik's black
noise (see Subsection \ref{subsection:CBS-examples}). Among the second
ones are fermionic white noises and bosonic white noises in `finite
temperature' representations of Araki-Woods type \cite{ArWo67a}, as
well as $q$-white noises, in particular free white noise (see
Subsection \ref{subsection:nCBS-examples}). As already stated for
$\Ao$-independence, more examples of $\Ao$-expected white noises can
be constructed easily from generalized Brownian motions, using again
the properties of white noise functors. Moreover, Anshelevich's
$q$-L'evy processes lead to further examples, as soon as it can be
proven that the vacuum vector is separating for the von Neumann
algebras generated by these processes \cite{Ansh04aTA}.

\subsection*{Section \ref{section:CGNSBS}:}
The example of Gaussian white noise makes it already evident that most
interesting processes, here Brownian motion, are not contained in the
$L^\infty$-space over the underlying measure space, but they are
contained in the corresponding $L^2$-space. In consequence, we will
extend the infrastructure of an $\Ao$-expected continuous Bernoulli shift
$\CBS$ to the GNS Hilbert bimodule $\bimod{\cA}$ (see Subsection
\ref{subsection:hm-spaces}). Thus, we will obtain a `homogeneous
continuous commuting square system of pointed Hilbert bimodules', in
allegory to \cite{Tsir04a} or alternatively, an
\emph{$\Ao$-expected (non-commutative) continuous GNS Bernoulli shift}.
We will develop the representation theory of such shifts only as far
as it is necessary for this paper. A key result is Proposition
\ref{new-product} which establishes the product of two Hilbert bimodule
elements, as long as they are $\Ao$-independent. This provides a
non-commutative extension of the well-known result in probability
theory that the product of two stochastically independent
$L^2$-functions is again an $L^2$-function. Furthermore, we will show
that the shift and all the conditional expectations extend to
adjointable operators on the Hilbert bimodules (see Theorem
\ref{lap-morph}). This will put us finally into the position to
introduce in Definition \ref{HMBS} an $\Ao$-expected non-commutative
continuous Bernoulli shift 
\[ 
     \big(\bimod{\cA},\1, \overline{S},
     (\bimod{\cA_{I}})_{I \in \cI}\big) 
\] 
that encodes, similarly as before, the covariant shift of a commuting
 square system of pointed Hilbert bimodules:
\begin{align*}
  \begin{pmatrix}
    \bimod{\cA_{[r,s]}}     & \subset & \bimod{\cA_{[r,v]}}        \\
    \cup                    & \quad   & \cup                       \\
    \bimod{\Ao}             & \subset & \bimod{\cA_{[u,v]}}
  \end{pmatrix}_{r < s <u < v}
  \hskip-14pt
  \xrightarrow{\phantom{x}\overline{S}_{t}\phantom{x}}
  \begin{pmatrix}
    \bimod{\cA_{[r+t,s+t]}} & \subset & \bimod{\cA_{[r+t,v+t]}}    \\
    \cup                    & \quad   & \cup                       \\
    \bimod{\Ao}             & \subset & \bimod{\cA_{[u+t,v+t]}}
  \end{pmatrix}_{r < s < u < v}
\end{align*}
Here denotes $\overline{S}$ the extension of the shift $S$ to the
bounded $\Ao$-linear operators on $\bimod{\cA}$. The Hilbert modules
are `pointed', because the cyclic separating vector $\1$ is contained
in each Hilbert bimodule and provides a (trivial) multiplicative shift
cocycle (see Subsection \ref{subsection:uc-cgnsbs}). Now it is
elementary to see that
\begin{enumerate}
\item[(iii)] if $\Ao \simeq \Cset$ and $\cA_{[r,t]}\simeq
  \cA_{[r,s]}\otimes \cA_{[s,t]}$ for $r < s < t$,
  then the continuous GNS Bernoulli shift is a homogeneous continuous
  product system of pointed Hilbert spaces, in the sense of Tsirelson
  \cite{Tsir04a,Koes04b}.
\item[(iv)] under the assumptions of (iii), the family
  $(\bimod{\cA_{[0,t]}})_{0 < t < \infty}$ defines a continuous tensor
  product system of Hilbert spaces, now in the sense of Arveson
  \cite{Arve03a,Tsir04a,Koes04b}.
\end{enumerate}
Tsirelson-Vershik's noises provide a rich source of Arveson's product
systems, in particular of type~II. Consequently,  so do
$\Cset$-expected commutative  continuous Bernoulli shifts. The
relationship between these three approaches will be further explored
in \cite{Koes04b}.

\subsection*{Section \ref{section:cocycles}:}
Like stationary $G$-flows  play a central role for a
Tsirelson-Vershik's noise, so do cocycles for a 
continuous Bernoulli shift or its GNS representation. These
cocycles are \emph{adapted to the filtration} of the
continuous Bernoulli shift and satisfy either additive
or multiplicative cocycle equations. The multiplicative cocycles come
in two kinds: \emph{unitary cocycles} for the continuous Bernoulli
shift itself (Definition \ref{unitary-cocycle}) and \emph{unital
cocycles} for its GNS representation (Definition
\ref{unit-cocycle}). The additive cocycles will always be given in the
GNS representation (Definition \ref{add-cocycle}). These cocycles
provide non-commutative versions of L\'evy processes, as applications
of Junge-Pisier-Xu's  non-commutative martingale inequalities show
\cite{PiXu97a,Koes00a,JuXu03a,Koes03a,Koes04a}. Throughout this
work we will consider only additive cocycles with a uniformly
bounded variance operator. Also, the unital cocycles are such that
their compression to $\Ao$ is a uniformly continuous semigroup. 

Let us now present our main results from Section
\ref{section:cocycles}. We will assume here for simplicity that the
$\Ao$-expected  continuous (GNS) Bernoulli shift enjoys
$\dim \Ao < \infty$ (see Theorems \ref{theorem:ln-exp} and
\ref{main-theorem} for a more general case).

\begin{Thm}\label{theorem:intro-main-result} 
  There exists a bijective correspondence between unital cocycles
  and additive cocycles, where the latter ones satisfy some
  structure equation.
\end{Thm} 

Since every unitary cocycle defines a unital cocycle, we
conclude immediately from Theorem \ref{theorem:intro-main-result} 
(see also Theorem \ref{main-theorem-unitary}):

\begin{Thm} \label{theorem:correspondence-unitary}
  There exists a bijective correspondence between unitary cocycles
  and additive cocycles, where now the latter ones satisfy a
  stronger version of the structure equation.
\end{Thm} 

To reveal this stronger structure equation will be the topic of sequel
publications (see Subsection \ref{subsection:CBS-WN} and \cite{Koes03a}
for the case of a tracial state). Theorem
\ref{theorem:correspondence-unitary} resp.\ \ref{main-theorem-unitary}
can be regarded as a non-commutative version of Tsirelson's result that
every stationary $\cU(\Ho)$-flow (continuous in probability) is
`classical' (see \cite[Theorem 8a2]{Tsir04a} and also \cite{Tsir98a}).
Here denotes $\cU(\Ho)$ the unitary operators on the (separable)
Hilbert space $\Ho$ (corresponding to $\Ao$). (Notice that our result
does not cover fully Tsirelson's results, since we will stipulate
stronger continuity conditions.)
    
We emphasize that Theorem \ref{theorem:intro-main-result} is only based
on the infrastructure of continuous Bernoulli shifts and cocycles.
This puts us into the position to introduce non-commutative white
noises, ensuring that they are a non-commutative version (see
Definition \ref{white}) of Tsirelson's `classical' noises: 

\begin{Defn} \label{definition:intro-ncwn}
  An $\Ao$-expected (non-commutative) continuous Bernoulli shift
  $(\cA,\psi,S, (\cA_{[r,s]})_{r\le s})$ is called a
  \emph{(non-commutative) white noise} if the filtration
  $(\cA_{0,t})_{t > 0}$ is generated by its unitary cocycles (in an
  adapted manner).
\end{Defn}

Now the correspondence ensures that white noises, as defined above, are
always `generated' by additive cocycles. This is in parallel to
the well-known fact that, for example,  Brownian motion generates the
$\sigma$-algebras of the filtration of the Gaussian white noise (see
e.g. \cite{Tsir98a}). We will show that the  `non-commutative white
noise part' can always be extracted   from a continuous Bernoulli
shift by the compression with a conditional expectation (see
Subsection \ref{subsection:CBS-WN}).

The proof of Theorem \ref{theorem:intro-main-result} relies on a
non-commutative extension of stochastic It\^o integration and the
construction of non-commutative versions of stochastic exponentials
and logarithms, which we will develop in Sections
\ref{section:nc-ito-integration} and \ref{section:nc-ln-exp}.

\subsection*{Section \ref{section:nc-ito-integration}:} 
We will develop a theory of non-commutative It\^o integration which
includes an existence and uniqueness theorem for solutions of
non-commutative  It\^o differential equations. This theory relies only
on the structure  of non-commutative $\Ao$-continuous (GNS) Bernoulli
shifts and its  additive cocycles. Crucial for this approach is
the notion of $\Ao$-independence which allows to transfer the famous
It\^o isometry of Brownian motion to the non-commutative setting. The
starting point of this theory are preliminary results in
\cite{Prin89a} which gave the evidence that the present approach is
promising in its generality.

Our approach to non-commutative It\^o integration applies to all
examples of non-commutative white noises (in the sense of Definition
\ref{definition:intro-ncwn}), in particular fermionic, bosonic, free
and $q$-white noises, including the operator-valued setting. In the
case of scalar-expected gauge invariant bosonic white noise one recovers
early work on non-Fock bosonic quantum stochastic integration by
Hudson, Lindsay and Wilde \cite{HuLi85a,LiWi86a}, which is constructed
through an amplification of Hudson-Parthasarathy quantum stochastic
integration on symmetric Fock spaces \cite{HuPa84a,Part92a}. Moreover,
one meets the pioneering work on quasi-free quantum stochastic
integrals for the CAR and CCR algebra by Barnett, Streater and Wilde
\cite{BSW83a}. Of special importance in applications is the It\^o
integration theory for so-called `squeezed white noises' which are
relevant in modern quantum optics \cite{GaZo00a}. Their It\^o
integration theory is treated in \cite{HHKKR02a}, based on the present
approach. If the underlying non-commutative probability space comes
from Voiculescu's free probability theory, one is precisely in the
$\Cset$-expected setting of Biane-Speicher's free stochastic calculus
\cite{BiSp98a}. In the case of $q$-commutation relations, one obtains a
non-commutative theory of It\^o integration, as contained already in
\cite{HKK98a} and independently much further developed in
\cite{DoMa03a}.

\subsection*{Section \ref{section:nc-ln-exp}:} 
We will develop the theory of non-commutative logarithms and
exponentials for unital cocycles resp.\ additive cocycles, as
far as it is necessary for the proof of Theorem
\ref{theorem:intro-main-result}. We will introduce the mapping $\Exp$
from the set of additive cocycles (with structure equation) to the
set of unital cocycles (Subsection \ref{subsection:nc-exp}) and the
mapping $\Ln$ from the set of unital cocycles to the set of
additive cocycles (Subsection \ref{subsection:nc-ln}). Finally, we
will show in Subsection \ref{subsection:proof-correspondence} that the
mappings $\Exp$ and $\Ln$ are each others inverse. This result
completes the proof of the main theorems on the correspondence of
additive and unital cocycles, as they are stated in Subsection
\ref{subsection:correspondence}.


\section{Preliminaries}                   \label{section:preliminaries}
We will fix the basic mathematical terminology  for a non-commutative
extension of probability theory,  as we will use it throughout this
paper.
 
\subsection{General terminology}      \label{subsection:prelim-general}
Throughout, $\cA$ is a von Neumann algebra in $\cB(\cH)$, the bounded
operators on some fixed Hilbert space $\cH$. We require that $\cA$ has
a separable predual $\cA_{*}$. Beside the norm topology on $\cA$, we
consider the weak* topology $\si(\cA,\cA_{*})$, the strong operator
(stop) topology and the $\si$-strong operator ($\si$-stop) topology
induced by the seminorms $\d_{\xi}(x):=\nm{x\xi}$, $\xi\in\cH$ resp.~
$\d_{\phi}(x):=\ab{\phi(x^{*}x)}^{1/2}$, $\phi\in\cA_{*}$. The unit of
$\cA$ is denoted by $\1_{\cA}$, or simply by $\1$, if no confusion can
arise.

Since throughout $\cA$ is considered in the presence of a fixed
faithful state $\psi\in \cA_{*}$, we assume for our
convenience that $\cH$ is already the GNS Hilbert space $\cH_{\psi}$
corresponding to $\psi$. Thus we have $\psi = \skp{\Om}{\cdot\Om}$ for
some vector $\Omega \in \cH_{\psi}$ which is cyclic and separating for
$\cA$. Moreover, since $\cA_*$ is separable, the GNS Hilbert space
$\cH_{\psi}$ is also separable. Notice that the scalar product is
taken to be linear in the second component. Two elements $x,y \in \cA$
are called $\psi$-orthogonal if $\psi(x^*y) =0$. Finally, the von
Neumann algebra generated by a family $(\cA_{j})_{j\in J}\subset \cA$
is denoted by $\bigvee_{j\in J}\cA_{j}$.
 
As usual, the von Neumann algebra $\cA'$ is the commutant of $\cA$ in
$\cB(\cH)$ and $\cZ(\cA):=\cA\cap\cA'$ is the center of $\cA$. The von
Neumann algebra $\cA$ is called a factor if $\cZ(\cA)\simeq\Cset$. For
a faithful normal state $\psi$ on $\cA$, the associated modular
automorphism group is denoted by $\si^{\psi}$. The centralizer
$\cA^{\psi}:=\set{x\in\cA}{\psi(xy)=\psi(yx) \text{ for all\ }
y\in\cA}$ is the fixed point algebra of $\si^{\psi}$.

We will use the modulus $|a|:= (a^*a)^{1/2}$ and, occasionally, $\RE
a := (a + a^{*})/2$ and $\IM a := (a - a^{*})/2\i$ for $a \in
\cA$. Finally, for any normed linear space $\cN$ we denote by
$\cN_{1}=\set{x\in\cN}{\nm{x}\le 1}$ the  unit ball of $\cN$.

\subsection{Non-commutative probability spaces and their morphisms}
\label{subsection:spaces-morphisms} 
The pair $\Apsi$ will be understood as a \emph{(non-commutative)
probability space} consisting of a von Neumann algebra $\cA$ which is
equipped with a faithful normal state $\psi$. A \emph{subalgebra $\cB$
of $\Apsi$} is a von Neumann subalgebra $\cB$ of $\cA$ such that the
conditional expectation $E_\cB$ from $\cA$ onto $\cB$ exists and
leaves $\psi$ invariant. We remind that such a conditional expectation
exists (uniquely) if and only if $\cB$ is globally invariant under the
action of the modular automorphism group of $\Apsi$  \cite{Takesaki2}.
 
The morphisms of $\Apsi$ are completely positive unital maps $T$ on
$\cA$ such that $\psi$ is $T$-invariant. They are automatically normal
(see Lemma \ref{morph-cont}) and we denote them by
$\Mor(\cA,\psi)$. Similarly, $\Aut(\cA,\psi)$ denotes the automorphisms
of $\Apsi$. Conditional expectations, as we will consider them
throughout this paper, are always morphisms. The identity map on $\cA$
is denoted by $\id_\cA$ or just by $\id$.
 
An \emph{$\Ao$-expected (non-commutative) probability space
  $(\cA,\psi)$} is a probability space with a distinguished subalgebra
$\Ao$ of $(\cA,\psi)$. Occasionally, such a space will also be
denoted as the triple $(\cA,\psi,\Ao)$. The terminology is motivated
from the fact that, given some von Neumann algebra $\cB$, an injective
*-homomorphism $\iota\colon \cB \to \cA$ with
$\Ao:=\iota(\cB)$ is a
non-commutative random variable (compare e.g.~\cite{Kuem88a}). Here we
will always identify $\Ao$ and $\cB$. Finally,
two elements $x,y \in \cA$ are called $\Ao$-orthogonal if
$E_{\Ao}(y^*x) =0$.
 
Typical examples of $\Ao$-expected probability spaces are the
following: Let $(\Ao,\psi_0)$ and $(\cB,\phi)$ be two probability
spaces. Then an ($\Ao \otimes \id_\cB$)-expected probability space
$(\cA,\psi)$ is defined by the von Neumann algebraic tensor product
$\cA = \Ao \otimes \cB$ and the tensor product state $\psi= \psi_0
\otimes \phi$. If $\cB$ is commutative, then we are in the
context of an operator-valued probability theory.

\subsection{Filtrations}                 \label{subsection:filtrations}
Let $\Apsi$ be an $\Ao$-expected probability space and
$(\cA_{I})_{I\in \cI}$ a family of subalgebras of $\Apsi$ which is
indexed by the set $\cI$ of possibly degenerated or possibly unbounded
intervals $I \subseteq \Rset$ and the empty set
$\emptyset$. Degenerated  intervals are points and will also be
written as $[t,t]$, $t\in \Rset$. It is called a \emph{filtration of
$\Apsi$} if $I \subseteq J$ implies $\cA_I \subseteq \cA_J$ for any
intervals $I,J \in \cI$ (\emph{monotony}). The filtration is
\emph{minimal} if $\bigvee\set{\cA_I}{I \in \cI  \text{\ bounded}}
=\cA$. In particular, a minimal filtration enjoys
$\cA_{\Rset}=\cA$. Notice that monotony is equivalent to
$\cA_{I}\vee\cA_{J} \subseteq \cA_{K}$ whenever $I\cup J=K$.
Finally, the sub-filtrations  $(\cA_{(-\infty,t]})_{t\in \Rset}$ and
$(\cA_{[t,\infty)})_{t\in \Rset}$ are, respectively, called the
\emph{past} and \emph{future filtrations}.
 
A filtration $(\cA_I)_{I\in \cI}$ is \emph{continuous downwards} if
$\bigcap_{\ve >0} \cA_{[s-\ve, t+\ve]}=\cA_{[s,t]}$ for any $s\le
t$. It is sufficient to check the downward continuity of a filtration
for the intersection of closed intervals: from $\cA_{[s,t]}\subseteq
\cA_{(s-\ve, t+\ve)}$ for any $\ve >0$ follows $ \cA_{[s,t]} \subseteq
\bigcap_{\ve >0} \cA_{(s-\ve, t+\ve)} \subseteq \bigcap_{\ve >0}
\cA_{[s-\ve, t+\ve]}$. A filtration $(\cA_I)_{I \in \cI}$ is
\emph{continuous upwards}  if $\bigvee_{\ve>0} \cA_{[s+\ve, t-\ve]} =
\cA_{[s,t]}$ for any $s<t$ (for notational simplicity, the evident
condition $t-s > 2\ve$ is always suppressed). Since $\bigvee_{\ve
>0}\cA_{[s+\ve, t-\ve]} \subseteq \cA_{(s,t)}$, the upward continuity
of a filtration implies immediately $\cA_{(s,t)}= \cA_{[s,t]}$ for any
$s<t$. A filtration is called \emph{continuous} if it is continuous
downwards and upwards. The continuity properties of the past or future
filtration are understood similarly. Notice that the (downward resp.\
upward) continuity of a filtration $(\cA_I)_{I \in \cI}$ is equivalent
to the (downward resp.\ upward) continuity of the associated family
$(E_I)_{I\in \cI}$ of conditional expectations $E_I\colon (\cA, \psi)
\to \cA_I$ in the pointwise stop topology \cite{Hell01a,Koes00a}. In
particular, the past filtration $(\cA_{(-\infty,t]})_{t \in \Rset}$ is
continuous if and only if the family $(E_{(-\infty,t]})_{t \in \Rset}$
is continuous in the pointwise stop topology. A similar equivalence is
valid for the future filtration.

\section{Non-commutative independence}     \label{section:independence}

We will introduce \emph{$\Ao$-independence} in Definition
\ref{definition-independence}  as a non-commutative analogue of
(amalgamated) stochastic independence. Such a notion of non-commutative
independence is in parallel to that 
of commuting squares in subfactor theory \cite{Popa83a,GHJ89a}.
It  will be crucial for the introduction of
\emph{continuous Bernoulli shifts} in Section \ref{section:CBSI}, as
well as non-commutative It\^o integration in Section
\ref{section:nc-ito-integration}.

Subsection \ref{subsection:independence-Popa} provides the definition
of \emph{$\Ao$-independence} and states its elementary properties
which are well-known in subfactor theory. In Subsection
\ref{subsection:independence-examples} we will relate
$\Cset$-independence to classical stochastic independence and to CCR
(or bosonic) independence. Subsection
\ref{subsection:C2A0-independence} provides elementary tools how
$\Cset$-independence is upgraded to $\Ao$-independence by tensor
products of probability spaces. Finally, we illustrate
$\Ao$-independence by further examples coming from non-commutative
probability theory in Subsection
\ref{subsection:nc-independence-examples}. This list of examples
includes CAR (or fermionic) independence, Voiculescu's free independence
and $\Ao$-independence in the context of $q$-commutation relations.

\subsection{$\Ao$-independence and Popa's commuting squares} 
\label{subsection:independence-Popa} 
We remind that in our terminology the phrasing `$\cB$ is subalgebra of
$(\cA,\psi)$' always means: $\cB$ is a von Neumann subalgebra of $\cA$
such that the conditional expectation $E_\cB\colon (\cA,\psi)\to\cB$
exists (see Subsection \ref{subsection:spaces-morphisms}).

\begin{Defn}\label{definition-independence}
  Let $\cB$ and $\cC$ be two subalgebras of the \Ao-expected probability
  space $(\cA,\psi)$ such that $\Ao \subseteq \cB \cap \cC$. The
  algebras $\cB$ and $\cC$ are called $\Ao$-independent, if for any
  $x\in\cB$ and $y\in\cC$
  \begin{align}\label{factor}
         E_{0}(x y)
      =  E_{0}(x)E_{0}(y)\,. 
  \end{align} 
  Here denotes $E_0\in \Mor\Apsi$ the conditional expectation from
  $\cA$ onto $\Ao$. Two families $(x_i)_{i\in I}$ and $(y_j)_{j \in J}$
  in $\cA$ are $\Ao$-independent if $\bigvee_{i \in I}\{x_i\}$ and
  $\bigvee_{j \in J}\{y_j\}$ are $\Ao$-independent.
\end{Defn}

Such a structure was introduced by Popa as a `commuting square' in
subfactor theory \cite{Popa83a}. Thus the statements `$\cB$ and $\cC$
are $\Ao$-independent' and `$\cB$ and $\cC$ form a commuting square
over $\Ao$' are essentially the same. They will used both, depending
on whether we want to emphasize the probabilistic or more the algebraic
aspect.
 
If the von Neumann algebra $\Ao$ is one-dimensional, i.e., $\Ao
\sim \Cset$, then $E_0 = \psi(\cdot)\1_\cA$ is verified immediately.
Whenever it is convenient and does not produce confusion,
$\Ao$-independence is also called $\tilde{\cA}_0$-independence, if
$\tilde{\cA}_0$ and $\Ao$ are isomorphic as von Neumann
algebras. For example, this convention simplifies `$\Cset
\1_\cA$-independence' to `$\Cset$-independence'.
 
$\Ao$-independence is equivalent to other properties of the involved
von Neumann algebras, as it is well-known for commuting squares in
subfactor theory \cite{Popa83a,GHJ89a,JoSu97a}. We will 
make frequently use of this fact.

\begin{Prop}\label{independence-properties}
  Under the assumptions of Definition \ref{definition-independence}
  the following conditions are equivalent:
  \begin{enumerate} 
  \item\label{independence-properties-i} $\cB$ and $\cC$ are
    $\Ao$-independent;
  \item\label{independence-properties-ii} $E_{0}(x_{1}yx_{2}) =
    E_{0}(x_{1}E_{0}(y)x_{2})$ for any $x_{1},x_{2}\in\cB$,
    $y\in\cC$;\,
  \item\label{independence-properties-iii} $E_\cB(\cC) = \Ao$;
  \item\label{independence-properties-iv} $E_{\cB}E_{\cC} =E_{0};$
  \item\label{independence-properties-v} $E_{\cB}E_{\cC}=
    E_{\cC}E_{\cB}$ and $\cB\cap\cC=\Ao$.
  \end{enumerate}
\end{Prop}

\begin{proof} The equivalences follow from the proof given in   
  \cite[Prop. 4.2.1]{GHJ89a}, after some elementary modifications.
  \qedopt
\end{proof} 

\begin{Rem}\normalfont 
  (i) In general, Definition \ref{definition-independence} does not
  incorporate computational rules for expressions like $E_0(xyxy)$ or
  $E_0(yxyx)$. But expressions like $E_0(xyyx)$ or $E_0(yxxy)$ are
  pyramidally ordered and can be simplified to $E_0(xE_0(yy)x)$ resp.\ 
  $E_0(yE_0(xx)y)$ by the module property of conditional expectations.
  Whether enough information for the calculation of non-pyramidally
  ordered expressions is present, this depends on the additional
  algebraic structure of an example.
  
  (ii) $\Ao$-independence applies, in particular, to von Neumann
  algebras of type~III (see Example \ref{CCR-independence}). Notice
  also that $\cB\ \vee \cC$ may be contained  properly in $\cA$. If
  $\cA$ is the weakly* closed linear span of $\set{xy}{x \in \cB, y\in
    \cC}$, then the corresponding commuting square is said to be
  `degenerated' \cite{Popa83a,JoSu97a}. Such a situation appears if the
  von Neumann algebras carry enough (commutation) relations. But
  already \emph{free} probability leads to an important example of
  $\Ao$-independence with $\cA= \cB\vee\cC$, where the corresponding
  commuting squares are not `degenerated' (see Example
  \ref{free-independence}).
 \end{Rem} 

\subsection{Commuting subalgebras and $\Cset$-independence}   
\label{subsection:independence-examples}
A probability space that includes a pair of commuting von Neumann
subalgebras  produces $\Cset$-independence of the two subalgebras
which is characterized  algebraically as tensor product
independence. It comprises `classical independence' and `bosonic or
CCR independence'.

Let $\cB_1$ and $\cB_2$ be two (von Neumann) subalgebras of the
$\Cset$-expected probability space $(\cB,\phi)$. Suppose that  $\cB_1$
and $\cB_2$ commute, i.e., $xy= yx$ for $x \in \cB_1$,  $y \in \cB_2$,
and, for simplicity, that $\cB = \cB_1 \vee \cB_2$. Then are
equivalent:
\begin{enumerate}
\item\label{ind-ex-i}  $\cB_1$ and $\cB_2$ are $\Cset$-independent.
\item\label{ind-ex-ii} $(\cB,\phi)$ is canonically isomorphic to
  $(\cB_1 \otimes \cB_2, \phi_1\otimes\phi_2)$ with $\phi_i :=
  \phi_{|\cB_i}$ for $i=1,2$, where we identify $\cB_1$ with $\cB_1
  \otimes \1$ and $\cB_2$ with $\1 \otimes \cB_2$.
\end{enumerate}
Obviously, (\ref{ind-ex-i}) implies (\ref{ind-ex-ii}) and we are left
to prove the inverse implication. Let $x_i \in \cB_1$ and  $y_i \in
\cB_2$ ($i = 1,2,\ldots,n$). Since  $\phi(\sum_{i=1}^n x_iy_i)= \phi_1
\otimes \phi_2(\sum_{i=1}^n x_i\otimes y_i)$, the map $\sum_{i=1}^n
x_iy_i \mapsto \sum_{i=1}^n x_i\otimes y_i$ is well-defined and
extends to an isomorphism from  $(\cB,\phi)$ onto $(\cB_1 \otimes
\cB_2, \phi_1 \otimes \phi_2)$ which is implemented unitarily on the
corresponding GNS Hilbert spaces.

\begin{Exmp}[Classical independence]
  \normalfont\label{classical-independence}
  In the case of a commutative von Neumann algebra $\cB$, the notion of
  $\Cset$-independence is equivalent to the classical notion of
  independence. Let $\cB = L^\infty(\Omega,\Sigma,\mu)$ and $\psi(f)=
  \int f \d\mu$. Then $\cB_1$ and $\cB_2$ are independent if and only
  if the sub-$\sigma$-algebras of $\Sigma$ generated by $\cB_1$ and
  $\cB_2$ are independent.
\end{Exmp}

\begin{Exmp}[Bosonic or CCR independence]
  \normalfont\label{CCR-independence} The canonical commutation
  relations (CCR) lead to the first non-commutative example of
  $\Cset$-independence. It is convenient to introduce them in their
  Weyl form (see \cite{BrRo2,Petz90a} and cited literature therein).
  These relations are given by
  \begin{align*}
        W(f)W(g) 
   & =  \exp{(-\tfrac{\i}{2} \IM \skp{f}{g})} W(f+g),                \\
        W(f)W(f)^* 
   & =  W(f)^*W(f)=\1,
  \end{align*} 
  where $f,g$ are elements of the Hilbert space $\cK$ with scalar
  product $\skp{\cdot}{\cdot}$. They generate the C*-algebra
  $\operatorname{CCR}(\cK,\IM \skp{\cdot}{\cdot})$. 
  Consider on this C*-algebra the quasi-free (gauge invariant) state
  \[
         \psi_\lambda(W(f))
      =  \exp(-\tfrac{1}{4}(2\lambda+1)\,\nm{f}^2)        
  \] 
  for some fixed $\lambda >0$. Let $\cB$ denote the von Neumann algebra
  which is generated by $\set{W(f)}{f\in \cK}$ in the GNS
  representation associated to $\psi_\lambda$. Furthermore, let
  $\cK_1$ and $\cK_2$ be two orthogonal closed subspaces in $\cK$. Then
  the corresponding von Neumann algebras $\cB_i$, generated by
  $\set{W(f)}{f \in \cK_i}$ ($i=1,2$) in the GNS representation,
  commute and are $\Cset$-independent.\\ Above construction works also
  for more general quasi-free states on a CCR algebra, as we will see
  them in Example \ref{squeezed-wn}. Finally, let us remind that the
  condition $\lambda>0$ ensures that $\psi_\lambda$ extends to a
  \emph{faithful} normal state on $\cB$.
\end{Exmp}

\subsection{From $\Cset$-independence to $\Ao$-independence}
\label{subsection:C2A0-independence}
Examples of {$\Ao$-}inde\-pendence are  produced canonically from
examples of $\Cset$-independence. Given in addition to $(\cA,\psi)$
the probability space $(\cB_0,\phi_0)$, we let $\cA :=\cB_0 \otimes
\cB$, $\psi :=\phi_0 \otimes \phi$ and $\cA_i := \cB_0 \otimes \cB_i$
for $i=1,2$, and $\Ao := \cB_0 \otimes \1$. Then $(\cA,\psi)$ is an
$\Ao$-expected probability space where $E_0 = \id \otimes
\phi(\cdot)\1$ is the conditional expectation onto $\Ao$. Moreover,
$\cA_1$ and $\cA_2$ are $\Ao$-independent if and only if $\cB_1$ and
$\cB_2$ are $\Cset$-independent.

\begin{Exmp}\normalfont
  Following the above (widespread) construction, classical independence
  (Example \ref{classical-independence})  leads immediately to examples
  of so-called amalgamated (or operator-valued) independence.
  Similarly, CCR-independence (Example \ref{CCR-independence}) leads to
  examples of $\Ao$-independence.
\end{Exmp}

The tensor product construction of an $\Ao$-expected probability space
from a $\Cset$-expected probability space looks very specific. It is
worthwhile to point out that, if $\Ao$ is isomorphic to the complex
$n \times n$-matrices $M_n$ with $2 \le n \le \infty$, this
construction captures already the general situation:

\begin{Prop}\label{matrizen-wertig}
  Let $(\cA,\psi)$ be an $\Ao$-expected probability space with $\Ao
  \simeq M_n$. Then there exists a probability space $(\cB,\phi)$ such
  that $\cA \simeq M_n \otimes \cB$ and, under this isomorphism, $\psi
  = \psi_{|\Ao} \otimes \phi$.
\end{Prop}

\begin{proof}
  We define $\cB$ as the relative commutant of $\Ao$ in $\cA$. Then
  $\cA$ splits  canonically into the tensor product $M_n \otimes \cB$,
  \cite[11.4.11]{KaRi2} and we may assume $\cA= M_n \otimes \cB$. Now
  $\phi(x):= \psi(\1\otimes x)$ defines a normal state on $\cB$. It is
  checked  immediately that the conditional expectation $E_0$ from
  $(\cA,\psi)$ onto $\Ao$ acts as $E_0(x \otimes y)= x \otimes
  \phi(y)$ for any $x\in M_n$ and $y\in \cB$. We conclude $\psi(x
  \otimes y)= \psi(E_0(x \otimes y))=\psi(x \otimes \1) \phi(y)$ and
  thus $\psi = \psi_{|\Ao} \otimes \phi$, \cite[11.4.11,
  11.2.7]{KaRi2}.\qedopt
\end{proof}

\subsection{Non-commutative examples of $\Ao$-independence} 
\label{subsection:nc-independence-examples} 
In the remaining part of this section we present further examples of
$\Ao$-independence which, in particular, illustrate that
$\Ao$-independent von Neumann algebras may not commute.

\begin{Exmp}[Fermionic or CAR independence]
  \normalfont\label{CAR-independence} We will consider the canonical
  anticommutation relations (CAR) (see for example
  \cite[5.2.5]{BrRo2}). Let $\cK$ be a Hilbert space and let
  $\operatorname{CAR}(\cK)$ denote the C*-algebra, generated by the
  elements $\set{a(f)}{f\in\cK}$, satisfying for all $f,g \in \cK$: $f
  \mapsto a(f)$ is antilinear and
  \begin{align*} 
        a(f)a(g) + a(g) a(f) 
   & =  0,\\
        a(f)a(g)^* + a(g)^*a(f)
   & =  \skp{f}{g}\1.
  \end{align*}
  Consider on $\operatorname{CAR}(\cK)$ the
  quasi-free (gauge-invariant) state $\psi_{\la}$, defined by
  \[
        \psi_\lambda(a^*(f)a(g))
     =  \la \skp{g}{f}
  \] 
  for some fixed $\lambda$ with $0 < \lambda <1$
  \cite{Arak70a,Arak87a}. Let $(\cB, \psi_{\la})$ be the probability
  space obtained as the weak closure of $\operatorname{CAR}(\cK)$ in
  the GNS representation associated to $\psi_{\la}$. Let $\cK_i$,
  $i=0,1,2$, be mutually pairwise orthogonal closed subspaces in
  $\cK$ and denote by $\cB_i$ the von Neumann subalgebras generated
  by $\set{a(f)}{f\in \cK_i}$ ($i=0,1,2$) in  the GNS
  representation. One verifies immediately that  
  $\cB_0 \vee \cB_1$ and $\cB_0 \vee
  \cB_2$ are  $\cB_0$-independent. In particular, 
  $\cB_1$ and  $\cB_2$ are $\Cset$-independent.
\end{Exmp}

\begin{Exmp}[Free independence]\normalfont\label{free-independence}
  An important example of {$\Ao$-}indepen\-dence is given by
  Voiculescu's (amalgamated) free independence \cite{VDN92a}. Let
  $\cB_1, \cB_2$ be two subalgebras of the $\cB_0$-expected probability
  space $(\cB, \phi)$ such that $\cB_0 \subseteq \cB_1 \cap
  \cB_2$. Let $E_0\in \Mor(\cB,\phi)$ denote the conditional
  expectation onto $\cB_0$. The algebras $\cB_1$ and $\cB_2$ are
  $\cB_0$-freely independent if
  \[
        E_0(b_1b_2\ldots b_n)
     =  0
  \]
  whenever $E(b_i) =0$, $1 \le i \le n$ and $b_i \in \cB_{j(i)}$ with
  $j(i) \ne j(i+1)$, $1 \le i \le n-1$. Notice that for $\cB_0 \simeq
  \Cset$ this definition reduces to free independence with respect to
  the state $\phi$. It is elementary to check that $\cB_0$-free
  independence implies $\cB_0$-independence.
\end{Exmp}

\begin{Exmp}[$q$-Gaussian processes and
  $\Ao$-independence] \normalfont\label{q-independence} 
  A further example originates from the
  construction of $q$-Fock spaces ($-1 <q <1$) by Bo{\.z}ejko and
  Speicher \cite{BoSp91a}. Let $\cK_\Rset$ be a real Hilbert space and
  $\cK = \cK_\Rset\oplus\i\cK_\Rset$ its complexification. Then the
  family $\set{a(f)}{f\in\cK_\Rset}$, satisfying
  \begin{align*}
        a(f)a(g) + a(g)a(f) 
     =  0,   \qquad
        a(f)a(g)^* - q a(g)^*a(f) 
     =  \skp{f}{g}\1
  \end{align*}
  for all $f,g \in \cK_\Rset$, is realized as bounded linear operators
  on the $q$-Fock space $\cF_q(\cK)$. Let $\cB$ denote the von Neumann
  algebra generated by $q$-Gaussian processes $\set{\Phi(f):=a(f)+
  a(f)^*}{f\in\cK_\Rset}$, or equivalently in the case
  $\cK=L^2_\Rset(\Rset)$, generated by all increments of $q$-Brownian
  motions $(\Phi(\chi_{[s,t]}))_{s<t}$ (see \cite{BKS97a}). Then the
  vacuum vector $\Omega \in \cF_q(\cK)$ defines a tracial faithful
  normal state $\tau$ on $\cB$. If $\cK_1$ and $\cK_2$ are two
  orthogonal closed subspaces in $\cK_\Rset$, then the von Neumann
  subalgebras $\cB_i :=\bigvee\set{a(f)+a(f)^*}{f \in \cK_i}$ ($i=1,2$)
  are $\Cset$-independent. If $\cK_0$ is a third closed subspace,
  orthogonal to $\cK_1$ and $\cK_2$, which generates the von Neumann
  subalgebra $\cB_0$, then it is again elementary to verify that
  $\cB_0\vee \cB_1$ and $\cB_0\vee \cB_2$ are $\cB_0$-independent.
 \end{Exmp}

This list of examples can be continued  easily. More examples of
$\Cset$- or $\Ao$-independence arise from von Neumann
algebras generated by so-called generalized Brownian  motions on
deformed Fock spaces \cite{BoSp94a,BoGu02a,GuMa02a,Krol02a}, which
contain $q$-Brownian motions as a simple case. All related
constructions, necessary to provide these further examples,  are
captured by so-called functors of white noise, as  introduced in
\cite{Kuem85a} and further considered in \cite{GuMa02a}. In view of
Anshelevich's results on $q$-L\'evy processes \cite{Ansh04aTA}, it
arises the question whether they provide also examples of
$\Cset$-independence.

Aside of these quantum probabilistic approaches to construct new
examples, it is worthwhile to remind a second rich source for
$\Ao$-independence:  subfactor theory with all its commuting
squares.

\section{Continuous Bernoulli shifts I}            \label{section:CBSI}
This section is devoted to the introduction of a
\emph{(non-commutative) continuous Bernoulli shift}. Its discrete time
versions are (non-commutative) Bernoulli shifts on towers of von
Neumann algebras, as they appear in subfactor theory
(\cite{GHJ89a,Rupp95a}). On the other hand provide continuous Bernoulli
shifts a non-commutative extension of \emph{noises} in the sense of
Tsirelson \cite[Definition 2d1]{Tsir04a}. More technically speaking, a
\emph{noise} is a homogeneous continuous product system of probability
spaces. Thus we may call a \emph{continuous Bernoulli shift} also a
`homogeneous continuous commuting square system of non-commutative
probability spaces'. 

Let us outline the contents of this section. We begin in Subsection
\ref{subsection:cbs-basic} with the definition of a \emph{continuous
Bernoulli shift}, emphasizing the probabilistic point of view. We
will study some of its elementary  properties. These properties will
in particular justify its name. In Subsection \ref{subsection:type}
we prove that the von Neumann algebra of a $\Cset$-expected continuous
Bernoulli shift is either finite or of type~III. An immediate
consequence of this result is that such a continuous Bernoulli shift
has a properly infinite von Neumann algebra if and only if the state
of the continuous Bernoulli shift is non-tracial. After noting in
Subsection \ref{subsection:composition} that continuous Bernoulli
shifts are stable with respect to compositions by tensor products and
direct sums, we turn our attention to the question whether they are
also stable with respect to decompositions. We show in Subsection
\ref{subsection:decomposition} that this is indeed the case for
compressions given by conditional expectations or  orthogonal
projections, both subject to some natural conditions. This opens the
door to compressions onto the relative commutant which allows to
introduce `derived continuous Bernoulli shifts'. Further research in
this direction is planned and should provide classification results,
similar to subfactor theory.

A continuous Bernoulli shift, as stated in Definition \ref{def:cbs},
provides already the sufficient infrastructure for Theorem
\ref{main-theorem}, our main result on the correspondence between shift
cocycles. Most up to the present known examples of continuous Bernoulli
shifts enjoy additional algebraic structures. The study of these
additional algebraic and continuity features is postponed to Section
\ref{section:CBSII}. At this place the reader will also find more
detailed information on the relationship of continuous Bernoulli shifts
and Tsirelson's noises. Some examples of continuous Bernoulli shifts
will also be provided there.

Finally, we want to bring to the reader's attention that  the class of
continuous Bernoulli shifts is richer than those of
\emph{non-commutative white noises} (in the sense of our Definition
\ref{def-white}). This follows already from the surprising result of
Tsirelson and  Vershik \cite{TsVe98a} on the existence of \emph{black
noises} which are `non-classical' (in the terminology of
\cite[Definition 5c4]{Tsir04a}). The relation of continuous Bernoulli
shifts and non-commutative white noises  will be further specified in
Subsection \ref{subsection:CBS-WN}.

\subsection{Continuous Bernoulli shifts and their basic properties}
\label{subsection:cbs-basic}
We start with some notation and remind that our notion of a filtration
does not stipulate continuity properties (see Subsection
\ref{subsection:filtrations}).

\begin{Notation}\normalfont
  The set of all intervals $I$ in $\Rset$ and the
  empty set $\emptyset$ is denoted by $\cI$. Furthermore we let
  $I+t := \set{s+t}{s\in I}$ and $\emptyset + t := \emptyset$. The set
  $\Int{I}$ is the interior of $I$.
\end{Notation}

\begin{Defn}\label{def:cbs}
  Let $\Apsi$ be an $\Ao$-expected probability space, equipped with a
  pointwise weakly* continuous group $S = (S_t)_{t \in \Rset} \subset
  \Aut\Apsi$ and a minimal filtration $(\cA_{I})_{I\in \cI}$. The
  quadruple \CBS is called an \emph{$\Ao$-expected (non-commutative)
  continuous Bernoulli shift}  if it enjoys the following properties:
  \begin{enumerate}
  \item\label{WN-i} $\Ao$ is the fixed point algebra of $S$;
  \item\label{WN-ii} $S$ acts covariantly on the filtration:
    $S_{t}\cA_{I}=\cA_{I +t}$ for any $t\in \Rset$, $I \in \cI$;
  \item\label{WN-iii} $\cA_{I}$ and $\cA_{J}$ are $\Ao$-independent
    whenever $I \cap J = \emptyset$.
  \end{enumerate}
  The $\Ao$-expected continuous Bernoulli shift is said to be
  \emph{trivial} if $S = \id$.
\end{Defn}

For shortness, and if there is no chance of confusion, $\CBS$ as well
as its automorphism group $S$ will be both just called a shift. If the
shift is trivial, then $\CBS$ can be  identified canonically with a
non-commutative probability space.

\begin{Notation}\normalfont
  Throughout, $E_I$ denotes the conditional expectation from
  $(\cA,\psi)$ onto $\cA_I$, where $I \in \cI$.
\end{Notation}

We proceed with elementary results on properties of an $\Ao$-expected
continuous Bernoulli shift which, in particular, will justify its
name.

In Definition \ref{def:cbs} we required that $\cA_I$ and $\cA_J$ are
$\Ao$-independent if $I\cap J = \emptyset$. But boundary points of such
intervals don't matter for $\Ao$-independence. Moreover, we will see
that $\Ao$ equals $\cA_{\emptyset}$, from which we will 
take advantage occasionally in proofs.

\begin{Lem} \label{lemma:independenceI}
  For an $\Ao$-expected continuous Bernoulli shift $\CBS$ is $\Ao =
  \cA_{\emptyset}$.  Moreover, the following are equivalent:
  \begin{enumerate}
  \item[(iii)] $\cA_{I}$ and $\cA_{J}$ are $\Ao$-independent whenever
    $I \cap J = \emptyset$.
  \item[(iii')] $\cA_{I}$ and $\cA_{J}$ are $\Ao$-independent whenever
    $\Int{I} \cap \Int{J} = \emptyset$.
  \end{enumerate}
\end{Lem}

\begin{proof}
  We conclude $\cA_{\emptyset}=\cA_{\emptyset}\cap\cA_{\emptyset}=\Ao$
  from Definition \ref{def:cbs}\spc(\ref{WN-iii}) and Proposition
  \ref{independence-properties}\spc(\ref{independence-properties-v}).
  
  It is obvious that (iii') implies (iii). We are left to prove the
  inverse. For $\Int{I} \cap \Int{J} = \emptyset$ exists $t\in \Rset$
  such that, without loss of generality, $I \subset (-\infty,t]$ and $J
  \subset [t,\infty)$. From the continuity of the past and future
  filtration (see Lemma \ref{P-cont}\spc(\ref{P-cont-ii}) below) we
  conclude $E_{I}E_{J} = E_I E_{(-\infty,t]}E_{[t,\infty)}E_J = E_I
  E_{(-\infty,t)}E_{(t,\infty)}E_J = E_I E_{\emptyset} E_J =
  E_{\emptyset}$.  \qedopt
\end{proof}

We collect further, frequently used properties of a continuous
Bernoulli shift.

\begin{Lem} \label{P-cont} 
  A shift, as stated in Definition \ref{def:cbs}, enjoys the following
  properties:
  \begin{enumerate}
  \item\label{P-cont-i} The shift acts covariantly for any $t \in
    \Rset, I \in \cI$:
    \[          
           S_t E_I
       =   E_{I+t} S_t;
    \]
  \item\label{P-cont-ii} the past filtration $(\cA_{(-\infty,t]})_{t
    \in \Rset}$ and the future filtration $(\cA_{[t, \infty)})_{t \in
    \Rset}$ are continuous, or equivalently, the families of
    conditional expectations $(E_{(-\infty,t]})_{t \in \Rset}$ and
    $(E_{[t,\infty)})_{t \in \Rset}$ are pointwise weakly* continuous.
    In particular, $\cA_{(-\infty, t]} = \cA_{(-\infty, t)}$ and
    $\cA_{[t, \infty)}=\cA_{(t,\infty)}$;
  \item\label{P-cont-iii} $\cA_{[t,t]}= \Ao =\cA_{\emptyset}$ for any
    $t\in \Rset$;
  \item\label{P-cont-iv} $\Ao \subset \cA_I$ for any $I \in \cI$;
  \item\label{P-cont-v} a shift is \emph{tail trivial}: $\bigcap_{t \in
    \Rset} \cA_{(-\infty, t]} = \Ao = \bigcap_{t\in \Rset}
    \cA_{[t,\infty)}$;
  \item\label{P-cont-vi} a shift is \emph{locally trivial}, i.e.,
    $\bigcap_{\ve >0} \cA_{[t-\ve, t+\ve]}= \Ao$ for any $t\in \Rset$.
  \end{enumerate}
\end{Lem}

Notice that the filtration of a shift may not be continuous downwards
or upwards. These and additional properties will be discussed in more
detail in  Section \ref{section:CBSII}.

\begin{proof}
  From the covariant action of the shift $S$ we get
  \begin{align*}
        \psi(x S_t E_I(y))
   & =  \psi ( E_{I+t}(x) S_t E_I (y)) 
     =  \psi (S_{-t} E_{I+t}(x) E_I(y))                              \\
   & =  \psi (S_{-t}E_{I+t}(x) y) 
     =  \psi (E_{I+t}(x) S_t(y)) 
     =  \psi (x E_{I+t}S_t(y))
  \end{align*}
  for any $x,y \in \cA$. This implies (\ref{P-cont-i}).
  
  (\ref{P-cont-ii}) follows from (\ref{P-cont-i}) and the pointwise
  continuity of the shift $S$ in the weak* topology, since
  $E_{(-\infty, t-\ve ]} = S_{-\ve} E_{(-\infty,t]} S_{\ve}$. (The
  equivalence of the two formulations is shown by routine  arguments.)
  
  (\ref{P-cont-iii}) $\cA_{[t,t]}$ and $\cA_{[t,t]}$ are
  $\Ao$-independent by Lemma \ref{lemma:independenceI}. But this
  implies $\cA_{[t,t]}\cap \cA_{[t,t]}= \Ao$ (see Proposition
  \ref{independence-properties} (\ref{independence-properties-v})).
  The equality $\Ao = \cA_{\emptyset}$ is already shown in Lemma
  \ref{lemma:independenceI}.

  (\ref{P-cont-iv}) is part of the definition of $\Ao$-independence
  (or follows directly from (\ref{P-cont-iii}) by the monotony of the
  filtration).
  
  The first equality of (\ref{P-cont-v}) follows from the observation
  that $x \in \bigcap_{t \in \Rset} \cA_{(-\infty,t]}$ and $y \in
  \bigcup_{n \in \Nset} \cA_{[-n, \infty)}$ are $\Ao$-independent
  elements. It follows $\psi((x-\Eo(x))y)=\psi((x-\Eo(x)))\psi(y)
  = 0$. From the weak* density of $\bigcup_{n \in \Nset} \cA_{[-n,
  \infty)}$ in $\cA$, we conclude $x=\Eo(x)$, hence $\bigcap_{t \in
  \Rset} \cA_{(-\infty, t]} = \Ao$. The second equality of
  (\ref{P-cont-v}) is shown by the same arguments.
  
  We are left to prove (\ref{P-cont-vi}). From
  \begin{align*}
        \bigcap_{\ve >0} \cA_{[t-\ve,t+\ve]}
    \subseteq 
        \bigcap_{\ve >0} (\cA_{(-\infty,t+\ve]}
        \cap \cA_{[t-\ve, \infty)})  
     =  \Big[\bigcap_{\ve>0}\cA_{(-\infty,t+\ve]}\Big]
        \cap 
        \Big[\bigcap_{\ve>0}\cA_{[t-\ve,\infty)}\Big]
  \end{align*}
  we conclude with the continuity of the past and future filtration,
  and finally the $\Ao$-independence,
  \begin{align*}
        \Ao 
    \subseteq 
        \bigcap_{\ve>0} \cA_{[t-\ve,t+\ve]} 
   &\subseteq 
        \cA_{(-\infty,t]}\cap \cA_{[t,\infty)} 
     =  \Ao\,. \qedopteqn
  \end{align*}
\end{proof}

The following result states that the shift $S$ is strongly mixing.
It will be crucial in the proof of Theorem
\ref{type-scalar}.

\begin{Lem}\label{ws-lim-Po}
  Let \CBS be an $\Ao$-expected  shift. For any $x\in\cA$ it holds
  \begin{align*}
        \lim_{\ab{t}\to\infty} S_{t}(x) 
     =  \Eo(x)
  \end{align*}
  in the weak* topology.
\end{Lem}

\begin{proof}
  For bounded intervals $I$ and $J$ and $x\in\cA_{I}$, $y\in\cA_{J}$
  one calculates
  \begin{align*}
        \lim_{\ab{t}\to\infty}\psi(y S_{t}(x)) 
    =   \lim_{\ab{t}\to\infty}\psi(\Eo(y S_{t}(x))) 
    =   \psi(\Eo(y)\Eo(x))=\psi(y \Eo(x))\,.
  \end{align*}
  Since the
  filtration is minimal, these identities extend to arbitrary
  $x,y\in\cA$ by standard arguments. Now the assertion follows from the
  norm density of the functionals $\set{\psi(y\,\cdot)}{y\in\cA}$ in
  $\cA_{*}$ and the boundedness of the set $\set{S_{t}(x)}{t\in\Rset}$.
  \qedopt
\end{proof}

Let us mention a subtle fact in the above proof: the mixing property of
the shift hinges on the notion of the minimality of a filtration, as
introduced in Subsection \ref{subsection:filtrations}. It includes that
$\cA$ is approximated by the `local' net $(\cA_{[s,t]})_{-\infty <
s \le t< \infty}$. Thus the `global' structure of a continuous
Bernoulli shift is already determined by its `local' structure.

\subsection{The type of a $\Cset$-expected continuous Bernoulli shift}
\label{subsection:type}
We proceed with a result on the type of the von Neumann algebra of a
continuous Bernoulli shift. In particular, it shows that neither
type~I$_\infty$ nor type~II$_\infty$ can occur for the von Neumann
algebra of a scalar-expected shift, as stated in Definition
\ref{def:cbs}.

\begin{Thm}\label{type-scalar}
  Let $\CBS$ be a non-trivial $\Cset$-expected continuous Bernoulli
  shift. Then $\cA$ is either finite or of type~III. Moreover,
  $\cA$ is finite if and only if $\psi$ is a trace.
\end{Thm}

This result generalizes immediately to a factor-expected continuous
Bernoulli shift if one considers, instead of $\cA$, the relative
commutant ${\Ao}^\prime \cap \cA$.

Before we start the proof of the theorem, we note an immediate
conclusion.

\begin{Cor} 
  Let the $\Cset$-expected shift be given as stated in Theorem
  \ref{type-scalar}. If $\cA$ is a factor, then $\cA$ is either
  of type~II$_1$ or of type~III.
\end{Cor}

\begin{proof} 
  Choose $x\in\cA_{[0,1]}$ with $\Eo(x)\neq 0$ and $x_{n}:=S_{n}(x)$,
  $n\in\Nset$. Then, by $\Cset$-independence, $\set{x_n-\Eo(x_n)}{n
  \in \Nset}$ is a mutually $\psi$-orthogonal family in $\cA$, thus
  $\cA$ is infinite-dimensional. Consequently, the factor $\cA$
  is not of finite type~I, which proves the corollary. \qedopt
\end{proof}

\begin{proof}[Proof of Theorem \ref{type-scalar}.]   
  Let $z\in\cZ(\cA)$ be the maximal central semi-finite projection.
  Since for any automorphism $\al$ on $\cA$, the projection $\al(z)$ is
  again a central semi-finite projection, we have $\al(z)\le z$.
  Interchanging $\al$ with $\al^{-1}$ yields $\al(z)=z$ and thus in
  particular $S_{t}(z)=z$ for any $t \in \Rset$. By the mixing property
  of the shift $S$, as stated in Proposition~\ref{ws-lim-Po}, it
  follows $z\in \Cset\cdot\1$. If $z=0$ then $\cA$ is of type~III and
  we are done. Therefore we may assume in the following that $\cA$ is
  semi-finite. It follows that there exists a stop-continuous unitary
  group $(u_{s})_{s\in\Rset}\subset\cA$ with the property
  $\si_{s}^{\psi}(x) = u_{s}^{*}x u_{s}^{}$ for any $x\in\cA$ (see
  \cite[Prop. 8.14.13]{Pedersen}). Since $S_{t}$ and
  $\sigma_{s}^{\psi}$ commute, one concludes $S_{t}(u_{s}^{*})x
  S_{t}(u_{s}^{}) = u_{s}^{*}x u_{s}^{}$ for all $x\in\cA$ and
  $s,t\in\Rset$. Thus $S_{t}(u_{s}^{})u_{s}^{*} \in \cZ(\cA)$ for any
  $t\in\Rset$. From this follows, again by Proposition \ref{ws-lim-Po},
  that the weak* limit $\Eo(u_{s}^{})u_{s}^{*}$ of this sequence is an
  element of $\cZ(\cA)$. Moreover, one has $\Eo(u_{s}^{}) =
  \psi(u_s)\1$. Since $s\mapsto u_{s}$ is continuous in the stop
  topology and $u_0=\1$, it is $\psi(u_{s}^{})\ne 0$ for all $s$ in
  some small $0$-neighborhood. Consequently, $u_{s}^{*} \in \cZ(\cA)$
  for any $s$ in this $0$-neighborhood. But this implies
  $\si_{s}^{\psi}=\id$ for any $s\in\Rset$ by the group property of
  $\sigma^{\psi}_{}$. Thus $\psi$ is a normal trace on $\cA$ \cite[Lem.
  8.14.6]{Pedersen} and, since $\psi$ is faithful, we conclude that
  $\cA$ is a finite von Neumann algebra,
  \cite[Thm. V.2.4]{Take03a}.\\
  By the above arguments, we have in particular proven that the
  finiteness of $\cA$ implies that $\psi$ is a trace. The converse is
  obvious.\qedopt
\end{proof}

\begin{Rem}\normalfont \label{remark:classification}
  (i) Up to the present we know examples for $\Cset$-expected continuous
  Bernoulli shifts with a von Neumann algebra $\cA$ of type~I$_1$, 
  type~II$_1$ and type~III$_{\lambda}$ ($0 < \lambda \le 1$). We 
  conjecture that for a $\Cset$-expected continuous Bernoulli shift 
  $\cA$ cannot be of type~I$_n$ for $n \ne 1$, but we have yet not 
  been successful to establish this. Also, it is of interest to 
  investigate further the case of type~III$_0$.
  
  (ii) A general result for the type of ${\Ao}^\prime \cap\cA$ is
  obtained by disintegration theory and will be presented elsewhere. It
  relies on the fact that the structure of continuous Bernoulli shifts
  is stable with respect to relative commutants. In particular, we
  introduce the notion of a `derived' continuous Bernoulli shift and
  start to investigate its properties, similar as it is done in
  subfactor theory \cite{HeKo04}.
\end{Rem}

\subsection{Composition of continuous Bernoulli shifts}
\label{subsection:composition}
It is easy to see that shifts, as introduced in Definition
\ref{def:cbs}, are closed under tensor product and direct sum
compositions.

\begin{Prop}\label{composition}
  Let $(\cA^{(i)},\psi^{(i)},S^{(i)},(\cA_I^{(i)})_{I\in \cI})$ be
  $\Ao^{(i)}$-expected continuous Bernoulli shift ($i=1,2$).
  \begin{enumerate} 
  \item $\big(\cA^{(1)}\otimes \cA^{(2)},\psi^{(1)}\otimes
    \psi^{(2)},S^{(1)}\otimes S^{(2)}, (\cA_I^{(1)}\otimes
    \cA_I^{(2)})_{I\in \cI}\big)$ is an $(\Ao^{(1)}\otimes
    \Ao^{(2)})$-expected continuous Bernoulli shift.
  \item $\big(\cA^{(1)}\oplus \cA^{(2)},\mu \psi^{(1)}\oplus
    (1-\mu)\psi^{(2)},S^{(1)}\oplus S^{(2)}, (\cA_I^{(1)}\oplus
    \cA_I^{(2)})_{I\in \cI}\big)$ is an $(\Ao^{(1)}\oplus
    \Ao^{(2)})$-expected continuous Bernoulli shift whenever $0 < \mu <
    1$.
  \end{enumerate}
\end{Prop}

The tensor product of a $\Cset$-expected and a trivial $\Ao$-expected
shift gives rise to an $\Ao$-expected shift.

\begin{Defn}
  The tensor product of an $\Ao$-expected continuous Bernoulli shift
  $\CBS$ and a trivial $\cB$-expected continuous Bernoulli shift is
  called the \emph{amplification of $\CBS$ by $\cB$}.
\end{Defn}

\begin{Rem}\normalfont
  The composition by tensor product extends to infinite tensor
  products. With some more technical efforts the procedure of direct
  sum compositions carries over to direct integrals of
  $\Ao^{(\gamma)}$-expected shifts, $\gamma \in \Gamma$, with respect to
  some standard probability space $(\Gamma,\mu)$.
\end{Rem}

\subsection{Decomposition of continuous Bernoulli shifts}
\label{subsection:decomposition}
In the following we focus onto decompositions of continuous Bernoulli
shifts. They are closed under compression by conditional expectations
and by orthogonal projections, subject to some further conditions.
These compressions provide, roughly speaking, the inverse procedures
to the compositions as stated in Proposition \ref{composition}.

\begin{Prop}\label{E-compression}
  Let $E \in \Mor(\cA,\psi)$ be a conditional expectation and $\CBS$ be
  an $\Ao$-expected continuous Bernoulli shift. If $E S_t = S_t E$ and
  $E E_I = E_I E$ for any $t\in \Rset$, $I \in \cI$, then $(E(\cA),
  \psi_{|E(\cA)}, S_{|E(\cA)}, (E(\cA_I))_{I\in \cI})$ is an
  $E(\Ao)$-expected continuous Bernoulli shift.
\end{Prop}

\begin{Defn}\label{definition:compression}
  The shift $(E(\cA), \psi_{|E(\cA)}, S_{|E(\cA)}, (E(\cA_I))_{I\in
  \cI})$ is called the \emph{compression} of $\CBS$ by the
  conditional expectation $E$.
\end{Defn}

If $E=\Eo$ then the compression leads to a trivial shift.

\begin{proof}[Proof of Proposition \ref{E-compression}.] 
  The restriction $\psi_{|E(\cA)}$ is a faithful normal state on
  $E(\cA)$. Since the conditional expectation $E$ and the shift $S$
  commute, $E(\cA)$ is globally invariant under the action of $S$. Thus
  the restriction $S_{|E(\cA)}$ is well-defined and has $E(\Ao)$ as
  fixed point algebra. Moreover, it is $S_t E(\cA_I) = E S_t(\cA_I) = E
  (\cA_{I +t})$ for any $t \in \Rset$, $I \in \cI$. Since $E(\cA_I)
  \subseteq E(\cA_J)$ whenever $I \subseteq J$, the family
  $(E(\cA_I))_{I\in \cI}$ defines a filtration of $(E(\cA), \psi_{|E
    (\cA)})$. From the minimality of $(\cA_I)_{I\in \cI}$ and the
  normality of $E$ we conclude that $\bigcup\set{E(A_I)}{I \in \cI
    \text{ is bounded}}$ is weakly* dense in $E(\cA)$. The minimality
  of the filtration $(E(\cA_I))_{I\in \cI}$ follows now by the double
  commutation theorem.
  Finally, the independence of $E(\cA_I)$ and E$(\cA_J)$ for $I \cap J
  = \emptyset$ is an immediate consequence of $E E_I E E_J = E E_I E_J
  = E \Eo$. \qedopt
\end{proof}

The following compression will be needed in Subsection
\ref{subsection:local-min-comp}. Recall that $\cA^\psi$  is the
centralizer of $(\cA,\psi)$.

\begin{Cor}
  $(\cA^\psi, \psi_{|\cA^\psi}, S_{|\cA^\psi}, (\cA^\psi \cap
  A_I)_{I\in\cI})$ is an $(\cA^\psi \cap \Ao)$-expected continuous
  Bernoulli shift.
\end{Cor}

\begin{proof}
  We will prove that the conditional expectation $E$ from $\cA$ onto
  the centralizer $\cA^\psi$ commutes with the shift $S$ and the
  conditional expectations $E_I$. Since $\cA^\psi$ is the fixed point
  algebra of the modular automorphism group $\sigma^\psi$, 
  we have $E(x) = \lim_{T\to\infty} \frac{1}{2T}
  \int_{-T}^{T} \sigma_t^\psi(x) \d t$ in the stop topology for any
  $x\in\cA$. Since the modular automorphism group $\sigma^\psi$
  commutes with morphisms of $(\cA,\psi)$, we conclude that $E$
  commutes with $S_t$ and $E_I$ for any $t\in \Rset$, $I\in
  \cI$.\qedopt
\end{proof}

In Proposition \ref{E-compression} is required that the conditional
expectation $E$ commutes with the conditional expectations $E_I$ of
the filtration. Dropping this condition leads to more general
compressions.

\begin{Prop} \label{proposition:EE-compression}
  Let $\CBS$ be an $\Ao$-expected shift and $E\in \Mor(\cA,\psi)$ a
  conditional expectation with $S_t E = E S_t$ and $E \Eo = \Eo E$
  for all $t\in \Rset$. Then  $(\cB, \psi_{|E(\cA)}, S_{|E(\cA)},
  (E(\cA)\cap \cA_I)_{I \in \cI})$ is an $(E(\cA)\cap \Ao)$-expected
  continuous Bernoulli shift, where  $\cB := \bigvee_{I \in \cI}
  E(\cA)\cap \cA_I$.
\end{Prop}
Notice that $\bigvee_{I \in \cI} E(\cA)\cap \cA_I$ can be much
smaller than $E(\cA)$.

\begin{proof}
  All arguments in
  the proof of Proposition \ref{E-compression} are valid, aside those
  for the covariant action of the shift, for the minimality and the
  independence. We will alter them as follows. From $S_t(E(\cA) \cap
  \cA_I) = E S_t(\cA) \cap S_t(\cA_I) = E(\cA)\cap \cA_{I+t}$ it
  follows the covariant action of the shift. The minimality of the
  filtration $(E(\cA)\cap \cA_I))_{I\in \cI}$ is evident. Finally, let
  $\widetilde{E}_I$ denote the conditional expectation from $(\cA,
  \psi)$ onto $E(\cA)\cap \cA_I$. 
  Notice that $\widetilde{E}_I = \widetilde{E}_I
  E_I = E_I \widetilde{E}_I $, 
  $\widetilde{E}_I = \widetilde{E}_I E = E \widetilde{E}_I$, since 
  conditional expectations are $\psi$-selfadjoint (compare Theorem
  \ref{psi-adj}), and $\widetilde{E}_0 = E \Eo = \Eo E$, 
  since $E$ and $\Eo$ commute. Consequently, we obtain
  $\widetilde{E}_I \widetilde{E}_J = \widetilde{E}_I E_I E_J
  \widetilde{E}_J = \widetilde{E}_I \Eo \widetilde{E}_J =
  \widetilde{E}_I E \Eo \widetilde{E}_J = \widetilde{E}_I
  \widetilde{E}_\emptyset \widetilde{E}_J = \widetilde{E}_0$
  for all $I,J \in \cI$. Since $\psi$ and $S$ restrict to $E(\cA)$, the
  compression defines a continuous Bernoulli shift as spelled out in
  the proposition. \qedopt
\end{proof}

Proposition \ref{composition} states that shifts can be composed by
direct sums. In the following we present compressions which, roughly
speaking, provide the reverse procedure. For a non-zero orthogonal
projection $e\in \cA$ we define $\psi_e(x)=\psi(exe)/\psi(e)$ for any
$x \in \cA$.

\begin{Prop}\label{e-e-compression} 
  Let $e\in \Ao \cap \cA^\psi$ be a non-zero orthogonal projection.
  Then $(e\cA e, \psi_e, S_{|e\cA e}, (e\cA_Ie)_{I\in \cI})$ is an
  $e\Ao e$-expected continuous Bernoulli shift.
\end{Prop}

\begin{proof}
  Clearly, $\psi_e$ is a faithful normal state on the von Neumann
  algebra $e\cA e$. Moreover, we observe $\sigma^{\psi_e}_t (e \cA_I
  e)= e\sigma^\psi_t(\cA_I)e = e\cA_I e$. Thus, there exist uniquely
  conditional expectations $Q_I$ from $(e\cA e,\psi_e)$ onto $e\cA_I e$
  for $I\in \cI$ such that $Q_I(exe)= eE_I(x)e$ for any $x \in \cA$.
  Since $e$ is a fixed point of the shift $S$, the restriction of
  $S_{|e\cA e}$ is well-defined. It leaves the state $\psi_e$ invariant
  and has the fixed point algebra $e\Ao e$. Since $S_t (e \cA_I e)= e
  S_t (\cA_I) e = e \cA_{I + t} e$ for any $t\in \Rset$, $I \in \cI$,
  the shift acts covariantly on $(e\cA_I e)_{I\in \cI}$. The family
  $(e\cA_I e)_{I\in \cI}$ defines a filtration on $e\cA e$ since
  $e\cA_I e \subseteq e\cA_J e$ whenever $I \subseteq J$. From the
  minimality of $(\cA_I)_{I\in \cI}$ we conclude that $\bigcup_I \set{e
    \cA_I e}{I \in \cI \text{ is bounded}}$ is weakly* dense in $e\cA
  e$. This ensures the minimality of $(e\cA_I e)_{I\in \cI}$. Finally,
  the $(e\Ao e)$-independence of $e \cA_I e$ and $e \cA_J e$ for $I
  \cap J =\emptyset$ follows from $Q_I Q_J (exe)= e E_J E_I(x)e = 
  e\Eo(x)e= eE_{\emptyset}(x)e = Q_{\emptyset}(x)$. \qedopt
\end{proof}

\section{Continuous Bernoulli shifts II}          \label{section:CBSII}

A continuous Bernoulli shift is a family of von Neumann subalgebras
which enjoys  essentially two further structures: $\Ao$-independence
and shift covariance. We did not stipulate further algebraic
structure  elements in Definition \ref{def:cbs}. Here we will
approach systematically some of these additional features which a
continuous Bernoulli shift may carry. Wide parts of the terminology
will be in analogy to \cite{Arve03a} and/or \cite{Tsir04a}, in the
attempt to highlight common and different grounds. Let us outline
this section's contents.

In Subsection \ref{subsection:local-min-max} we will show that the
algebraic properties of \emph{local minimality/maximality} imply the
\emph{(downward/upward) continuity} of the filtration $(\cA_I)_{I\in
\cI}$. The subject of Subsection
\ref{subsection:enriched-independence} is an  \emph{enriched
independence structure}, as it appears up to the present in all known
examples. We proceed in Subsection
\ref{subsection:commuting-past/future}  with the situation of a
\emph{commuting past/future}, as it appears in the setting of
continuous tensor product systems. Next we discuss the case of
commutative  continuous Bernoulli shifts and relate them to
Tsirelson-Vershik's \emph{noises} in Subsection
\ref{subsection:commutative-BS}. For these additional algebraic
structures will also be commented briefly whether they are stable under
compressions  by conditional expectations (see Subsection
\ref{subsection:decomposition}). In Subsection
\ref{subsection:local-min-comp}  we will present an example of a
compression which obstructs the local  minimality of a shift. Finally,
examples of continuous Bernoulli shifts are provided in Subsection
\ref{subsection:CBS-examples}:  Gaussian and Poisson white noise, CCR
and CAR white noises and $q$-Gaussian white noises ($-1< q < 1$).

\subsection{Local minimality and local maximality}
\label{subsection:local-min-max}
The following two algebraic properties ensure the  continuity of a
filtration.

\begin{Defn} 
  The filtration $(\cA_I)_{I\in\cI}$ of an $\Ao$-expected continuous 
 Bernoulli shift is
  \begin{enumerate}
  \item[(i)] \emph{locally minimal} if $\cA_I \vee \cA_J = \cA_K$ for
    any $I,J,K\in \cI$ with $I \cup J =K$;
  \item[(ii)]\emph{locally maximal} if $\cA_{(-\infty,t]}\cap \cA_{[s,
    \infty)} = \cA_{[s,t]}$ for any $s\le t$.
  \end{enumerate}
\end{Defn}

\begin{Lem} \label{lemma-local-minmax} 
  Let $\CBS$ be an $\Ao$-expected continuous Bernoulli shift. If the
  filtration $(\cA_I)_{I\in\cI}$ is locally minimal (maximal), then it
  is  continuous upwards (downwards).
\end{Lem}

Consequently, such a filtration is continuous if it is locally minimal
and locally maximal.

\begin{proof}
  We first prove that local minimality implies upward continuity. For
  a closed interval $K$ with non-empty interior we choose two
  intervals $I$, $J$ and some $\ve >0$ with $I \cup J = K$ and $I+\ve,
  J -\ve \subset K$. The local minimality guarantees $\cA_I \vee \cA_J
  = \cA_K$. If any $x \in \cA_I$ and $y \in \cA_J$ can be approximated
  by sequences $(x_n)_{n\in \Nset}, (y_n)_{n \in \Nset}\subset \bigcup
  \set{\cA_{K_0}}{K_0=\overline{K_0} \subset \Int{K}}$,
  then $\cA_I \vee \cA_J\subseteq\bigvee \set{ \cA_{K_0} }{
  K_0=\overline{K_0} \subset \Int{K}} \subseteq \cA_{K}$
  implies the upward continuity. But such sequences are given by $x_n
  := S_{\ve/n}(x)$ resp.\ $y_n := S_{-\ve/n}(y)$, since they
  approximate $x$ resp.\ $y$ due to the pointwise weak* continuity of
  $S$.
  
  For a locally maximal filtration $(\cA_I)_{I\in\cI}$ holds
  \begin{align*}
        \bigcap_{\ve >0} \cA_{[s-\ve, t+\ve]}
     =  \bigcap_{\ve > 0}
        (\cA_{(-\infty,t+\ve]}  \cap \cA_{[s-\ve,\infty)})
     =  \Big[\bigcap_{\ve > 0}  \cA_{(-\infty,t+\ve]}\Big] \cap
        \Big[\bigcap_{\ve > 0}  \cA_{[s-\ve,\infty)}\Big].
  \end{align*}
  Now, by Lemma \ref{P-cont}\spc(\ref{P-cont-ii}), the continuity of
  the past filtration and future filtration establishes the downward
  continuity. \qedopt
\end{proof}

\begin{Cor}
  A locally minimal filtration $(\cA_I)_{I\in\cI}$ of a shift enjoys
  $\cA_I = \cA_{\overline{I}}$ for any $I\in\cI$.
\end{Cor}

\begin{proof}
  Lemma \ref{lemma-local-minmax} insures the upward continuity. Since
  $\bigvee_{\ve >0}\cA_{[s+\ve, t-\ve]} \subseteq \cA_{(s,t)}$, the
  upward continuity of a filtration implies immediately $\cA_{(s,t)}=
  \cA_{[s,t]}$ for any $s<t$. Similar arguments prove the two
  remaining cases of unbounded intervals. \qedopt
\end{proof}

\begin{Rem}\normalfont
  It is tempting to stipulate local minimality and/or local maximality
  in Definition \ref{def:cbs} of an $\Ao$-expected continuous Bernoulli
  shift. But local minimality is obstructed by compressions, as it will
  be shown in Subsection \ref{subsection:local-min-comp}. In the other
  case, it is elementary to see that local maximality is stable under
  the compression with conditional expectations. Local maximality is
  always present for continuous Bernoulli shifts, which are
  non-commutative white noises (see Definition \ref{def-white}). But we
  do not know whether every continuous Bernoulli shift is locally
  maximal.
\end{Rem}

\subsection{Enriched independence} 
\label{subsection:enriched-independence}
All examples of $\Ao$-expected shifts enjoy a richer independence
structure than $\Ao$-independence, at least as  known presently by the
authors. Many of these examples come in particular from  `functors of
white noise' (see \cite{Kuem85a, GuMa02a} for their  notion).
 
\begin{Defn}
  An $\Ao$-expected continuous Bernoulli shift or its filtration
  $(\cA_I)_{I\in \cI}$ has an \emph{enriched independence}  if $\cA_I$
  and $\cA_J$ are $\cA_{I \cap J}$-independent for any $I,J \in \cI$.
\end{Defn}

\begin{Prop}\label{enriched-cs}
  For an $\Ao$-expected shift $\CBS$ are equivalent:
  \begin{enumerate}
  \item\label{enriched-cs-iii} $\cA_I$ and $\cA_J$ are $\cA_{I\cap
    J}$-independent for any $I, J \in \cI$;
  \item\label{enriched-cs-ii} $\cA_{I}\cap\cA_{J}=\cA_{I\cap J}$ and
    $E_{I}E_{J}=E_{J}E_{I}$ for any $I,J \in \cI$;
  \item\label{enriched-cs-i} $(\cA_I)_{I\in \cI}$ is locally maximal
    and $E_{(-\infty,t]}E_{[s,\infty)} = E_{[s,\infty)}
    E_{(-\infty,t]}$ for any $s,t \in \Rset$.
  \end{enumerate}  
  A filtration $(\cA_I)_{I\in \cI}$ with enriched independence is
  continuous and enjoys $\cA_I = \cA_{\overline{I}}$ for all $I\in
  \cI$.
\end{Prop}

Notice that the enriched independence of a filtration does not imply
its local minimality (see Proposition \ref{compression-obstruction}
for a counterexample).

\begin{proof}
  The equivalence of (\ref{enriched-cs-iii}) and (\ref{enriched-cs-ii})
  is evident by Proposition
  \ref{independence-properties}\spc(\ref{independence-properties-v}).
  Also it is clear that that (\ref{enriched-cs-ii}) implies
  (\ref{enriched-cs-i}). We are left to prove the converse. Local
  maximality guarantees $\cA_{(-\infty,t]}\cap \cA_{[s,\infty)} =
  \cA_{[s,t]}$ for $s\le t$ and consequently, by the equivalence of
  (\ref{independence-properties-v}) and
  (\ref{independence-properties-iv}) in Proposition
  \ref{independence-properties},
  $E_{(-\infty,t]}E_{[s,\infty)}=E_{[s,t]}$. This yields $E_{(r,t]}
  E_{[s,u)} = E_{(r,t]} E_{(-\infty,t]} E_{[s,\infty)} E_{[s,u)} =
  E_{(r,t]} E_{[s,t]} E_{[s,u)} = E_{[s,t]} = E_{[s,u)} E_{[s,t]}
  E_{(r,t]} = E_{[s,u)} E_{(r,t]}$ for any $-\infty \le r < s < t < u
  \le \infty$. Aside of the cases $I \cap J = \emptyset$ or $I
  \subseteq I \cap J$ (which follow easily from the $\Ao$-independence
  resp.\ the monotony of the filtration), the proof (iii) $\Rightarrow$
  (ii) is completed if the filtration $(\cA_I)_{I\in \cI}$ is
  continuous upwards, entailing $\cA_{I}= \cA_{\overline{I}}$ for any
  $I\in \cI$.  But $E_{[s+\ve,t-\ve]} =
  E_{(-\infty,t-\ve]}E_{[s+\ve,\infty)}
  \xrightarrow{{\scriptstyle\ve\searrow 0}}
  E_{(-\infty,t]}E_{[s,\infty)} = E_{[s,t]}$ from the continuity of
  $t\mapsto E_{(-\infty,t]}$ resp.\ $s\mapsto E_{[s,\infty)}$ in the
  pointwise stop topology. Thus the filtration is continuous upwards.
  Finally, the local minimality implies the downward continuity.
  Consequently, $(\cA_I)_{I\in \cI}$ is continuous. \qedopt
\end{proof}

\begin{Cor}\label{cor-comp-enriched}
  If a continuous Bernoulli shift has an enriched independence, then
  also its compression by a conditional expectation $E$, as stated in
  Proposition \ref{E-compression}.
\end{Cor}

\begin{proof}
  This follows immediately from Proposition
  \ref{enriched-cs}\spc(\ref{enriched-cs-ii}) and $E E_I = E_I E$ for
  any interval $I\subset \Rset$, since $\big(E(\cA) \cap \cA_I\big)
  \cap \big(E(\cA) \cap \cA_J\big) = E(\cA) \cap \cA_{I \cap J}$ and
  $E_I E E_J E = E_J E E_I E$. \qedopt
\end{proof}

The stability of shifts with enriched independence under compressions
suggests to stipulate this structure in Definition \ref{def:cbs}.
Nevertheless, it is not needed for the proofs of our main results in
this paper. Moreover, for more general compressions (see Proposition
\ref{proposition:EE-compression}) the proof of Corollary
\ref{cor-comp-enriched} breaks down, since the conditional expectations
$E$ and $E_I$ may no longer commute.

\begin{Rem}\normalfont
  (i) It is an open problem to give examples of continuous Bernoulli
  shifts without enriched independence.
  
  (ii) Further structure can be added to the index set $\cI$ on $\Rset$
  to provide other enriched forms of $\Ao$-independence. Here we focus
  on (possible degenerated and unbounded) intervals as set $\cI$.
  Stipulating a Boolean algebra structure for $\cI$ (similar as done
  for example in \cite{TsVe98a}) leads to even more enriched forms of
  $\Ao$-independence.
\end{Rem}

\subsection{Commuting past and future}
\label{subsection:commuting-past/future}
Shifts with a commuting past/future form an interesting class on their
own, as they stem from tensor product independence (see Subsection
\ref{subsection:independence-examples}). Roughly speaking, such
commuting structures are present in Hudson-Parthasarathy's approach to
quantum probability theory \cite{Part92a}, in Arveson's approach to
continuous product systems of Hilbert spaces \cite{Arve03a} and in
Tsirelson-Vershik's approach to continuous product systems of
probability spaces \cite{Tsir04a}.

\begin{Defn}\label{definition:cp-f}
  A $\Cset$-expected shift $\CBS$ (or its filtration) is said to have a
  \emph{commuting past/future} if $\cA_{(-\infty,0]}$ and
  $\cA_{[0,\infty)}$ commute.
\end{Defn}

\begin{Lem}\label{PFcommute}     
  Suppose that a $\Cset$-expected shift $\CBS$ satisfies one of the
  following (equivalent) additional conditions:
  \begin{enumerate}
  \item\label{PFcommute-i} the past $\cA_{(-\infty,0]}$ and the future
    $\cA_{[0,\infty)}$ commute;
  \item\label{PFcommute-ii} $\cA_{(-\infty,t]}$ and $\cA_{[t,\infty)}$
    commute for any $t\in \Rset$;
  \item\label{PFcommute-iii} $\cA_I$ and $\cA_J$ commute for all
    $I,J\in \cI$ with $\Int{I} \cap \Int{J} = \emptyset$.
  \end{enumerate}
  If in addition the filtration is locally minimal, 
  then it is locally maximal, continuous and enjoys 
  an enriched independence.
\end{Lem}

It is well-known that local maximality does not imply local
minimality, even when adding commutativity of the von Neumann algebras
to the assumptions of this lemma (see for example
\cite[Rem. 3.9]{Tsir03a}).

\begin{proof}
  The equivalence of (\ref{PFcommute-i}) and (\ref{PFcommute-ii}) is
  obvious by stationarity, (\ref{PFcommute-iii})  clearly implies
  (\ref{PFcommute-ii}), the inverse implication follows from the
  monotony of the filtration and the $\Cset$-independence. For a
  locally minimal filtration with commuting past and future is
  $\cA_{(-\infty,t]}= \cA_{(-\infty,s]}\vee \cA_{[s,t]}$ isomorphic to
  $\cA_{(-\infty,s]}\otimes \cA_{[s,t]} \otimes \1_{\cA_{[t,\infty)}}$
  with $s<t$. Similarly, one decomposes $\cA_{[s,\infty)}$ and concludes
  $ \cA_{(-\infty,t]}\cap \cA_{[s,\infty)} \simeq (\cA_{(-\infty,s]}
  \otimes \cA_{[s,t]}\otimes \1_{\cA_{[t,\infty)}}) \cap
  (\1_{\cA_{(-\infty,s]}} \otimes \cA_{[s,t]} \otimes\cA_{[t,\infty)})
  \simeq \cA_{[s,t]}$. This shows the local maximality of the
  filtration. The continuity of the filtration follows directly from
  Lemma \ref{lemma-local-minmax}.  Finally, using ideas from Subsection
  \ref{subsection:C2A0-independence}, it is easy to see that
  $E_{(-\infty,t]}$ and $E_{[s,\infty)}$ commute. This ensures, by
  Lemma \ref{enriched-cs}, the enriched independence structure.\qedopt
\end{proof}

\begin{Cor}\label{cor-commuting-p/f}
  If the $\Cset$-expected shift $\CBS$ has a commuting past/future,
  then also its compression by a conditional expectation $E$, as stated
  in Proposition \ref{E-compression}.
\end{Cor}

\begin{proof}
  This is an elementary conclusion from Definitions
  \ref{definition:compression} and \ref{definition:cp-f}. \qedopt
\end{proof}

\begin{Rem}\normalfont
  $\Cset$-expected continuous Bernoulli shift with a locally minimal
  filtration and a commuting past/future lead to examples of
  continuous tensor product systems of W*-algebras (see
  \cite[Sec. 7.3, Def. 7.1]{Liebscher1}).
\end{Rem}

\subsection{Commutative von Neumann algebras}
\label{subsection:commutative-BS}
A rich source for $\Cset$-expected continuous Bernoulli systems  with a
commutative von Neumann algebra is provided  by probability
theory. Here we focus on the connection to  noises, as they appear in
the work of Tsirelson and Vershik  \cite{TsVe98a, Tsir98a, Tsir04a}.

Let us start with the measure theoretic notion of a continuous
Bernoulli shift. We assume throughout that the probability  spaces are
Lebesgue spaces and that the $\si$-algebras are completed.

\begin{Defn}\label{definition:CBSPS}
  A \emph{continuous Bernoulli shift (on a probability space)}
  consists of a probability space $(\Omega, \Sigma, \mu)$, a measure
  preserving Borel-measurable group $(s_t)_{t\in \Rset}$ on $\Omega$
  and a family of sub-$\sigma$-algebras $(\Sigma_I)_{I\in \cI}\subset
  \Sigma$, such that  for all $I,J,K  \in \cI$
  \begin{enumerate}
  \item[(o)]   $\Si$ is generated by the family 
               $\set{\Si_{I}}{I\in\cI\text{\normalfont{\ bounded}}}$;  
  \item[(i)]   $s_t$ maps $\Si_I$ onto $\Si_{I+t}$ for any 
               $t\in \Rset$;
  \item[(ii)]  $\Si_I$ and $\Si_J$ are independent
               whenever $\Int{I} \cap \Int{J} = \emptyset$;
  \item[(iii)] $\Si_K$ contains the sub-$\sigma$-algebra 
               generated by the union of $\Si_I$ and 
               $\Si_J$ whenever $K = I \cup J$.
  \end{enumerate}
  $(\Sigma_{I})_{I\in \cI}$ (or the continuous Bernoulli shift) is
  called \emph{locally minimal} if $\Si_K$ is generated by the union of
  $\Si_I$ and $\Si_J$ whenever $K = I \cup J$.
  
  A \emph{locally minimal} continuous Bernoulli shift on a probability
  space is also called \emph{Tsirelson-Vershik's noise} or a
  \emph{homogeneous continuous product (system) of probability spaces}
  (see \cite{Tsir98a} or \cite[Definition 2d1]{Tsir04a}).
\end{Defn}

It is elementary to turn a measure theoretic shift into an algebraic
shift in the sense of Definition \ref{def:cbs}. The other direction is
less elementary, but also a very familiar fact. Let us state in the
following result without proof and in our terminology. Actually, it is
an immediate corollary of Mackey's paper \cite{Mac}.

\begin{Defn}
  An $\Cset$-expected continuous Bernoulli shift $\CBS$ is called to be
  \emph{commutative} if $\cA$ is commutative.
\end{Defn}

\begin{Thm}\label{theorem:TVN-CBS}
  There is one-to-one correspondence between (isomorphism classes of)
  \begin{enumerate}
  \item (locally minimal) continuous Bernoulli shifts 
       $(\Omega, \Sigma, \mu)$, $(s_r)_{r\in
    \Rset}$, $(\Sigma_I)_{I\in \cI}$;
  \item (locally minimal) $\Cset$-expected commutative continuous 
         Bernoulli shifts $\CBS$.
  \end{enumerate}
  The correspondence is given by $\cA = L^\infty(\Omega, \Sigma, \mu)$,
  $\psi = \int \cdot\, \d\mu$, $S_r(f)= f\circ s_r$ for all $f \in
  L^\infty(\Omega, \Sigma, \mu)$ and $\cA_I = L^\infty(\Omega,
  \Sigma_I, \mu_{|\Sigma_I})$ for all $I \in \cI$.
\end{Thm}

The terminology `isomorphism' means in the present context `isomorphism
between the dynamical systems which preserves the filtration
structure'. Since we will use this equivalence only occasionally, we
have omitted the (evident) notion of an isomorphism between algebraic 
and measure theoretic continuous Bernoulli shifts.

\begin{Cor}\label{corollary:TVN}
  Tsirelson-Vershik's noise corresponds to a $\Cset$-expected continuous
  Bernoulli shift which is commutative, locally minimal, locally
  maximal and has a continuous filtration. Moreover it enjoys an
  enriched independence and a commuting past/future.
\end{Cor}

\begin{proof} 
  The isomorphism of Theorem \ref{theorem:TVN-CBS} respects local
  minimality, consequently Lemma \ref{PFcommute} applies.
\end{proof}

It is clear that all these properties can be translated back into
properties of the underlying probability spaces, for example like
upwards and downwards continuity in \cite[2d2, 2d4]{Tsir04a}. Let us
close this section with a remark on the notation. Since
Tsirelson-Vershik's noises are continuous, one always has $\cA_{(s,t)}
= \cA_{[s,t]} =: \cA_{s,t}$. The latter notation is used for example in
\cite{Tsir04a}.

\subsection{Local minimality and compressions}
\label{subsection:local-min-comp}
Local minimality of a shift is not aggregated to Definition
\ref{def:cbs}, since this property is not stable with respect to
compressions. This is already a well-known phenomenon in probability
theory (see Remark \ref{remark:obstruction-classical}). But the
situation is even more dramatic for probability spaces based on
properly infinite von Neumann algebras. We will describe a
class of shifts, for which a quite natural compression by a conditional
expectation destroys the local minimality. The von Neumann algebra of
these shifts is of type~III$_\lambda$ ($0 < \lambda < 1$), equipped
with a periodic state, and its compression onto the centralizer yields
a von Neumann algebra of type~II$_1$. The example also shows that this
obstruction cannot be removed or controlled by stipulating
algebraic structures as an enriched $\Ao$-independence and/or a
commuting past/future.

Recall that a state $\psi$ on $\cA$ is called periodic if there exists
$T>0$ such that $\sigma^\psi_T = \id$. The smallest such $T$ is called
the period of $\psi$  \cite{Take73a}.

\begin{Lem}\label{compression-obstruction}
  Let $\CBS$ be a non-trivial $\Cset$-expected continuous Bernoulli
  shift with a locally minimal filtration and a periodic state $\psi$.
  If the filtration carries a commuting past/future, then the 
  compressed filtration $(\cA_I \cap \cA^{\psi})_{I\in \cI}$ is not 
  locally minimal.
\end{Lem}

\begin{proof}
  Since the shift is $\Cset$-expected, $S \subset \Aut(\cA,\psi)$ acts
  ergodically on $\cA$. Thus $\psi$ is a homogeneous periodic state
  with period $T>0$  \cite{Take73a}. Put $\kappa := e^{-2\pi/T}$ and
  let $\cB_n := \set{x \in \cA}{\sigma^\psi_t(x) = \kappa^{\i nt}x}$
  for $n \in \Zset$. Note that $\cB_0 =\cA^\psi$. Furthermore, $\cB_n$
  and $\cB_m$ are $\psi$-orthogonal whenever $n\neq m$. Recall also
  $\cB_n \cB_m \subseteq \cB_{n+m}$ for any $n,m \in \Zset$. Finally,
  we will need from \cite{Take73a} that $\ve_n(x): =
  \frac{1}{T}\int_{0}^{T} \kappa^{-\i nt}\sigma_t^\psi(x) \d t$ is a
  projection from $\cA$ onto $\cB_n$ and $x\Om =
  \sum_{n\in\Zset}\ve_{n}(x)\Om$ for any $x\in\cA$.
  
  Since $\sigma_t^\psi$ commutes with $S_t$ and $E_I$, also $\ve_n$
  does so. Thus, one has $S_t (\cB_n) = \cB_n$ and $E_I(\cB_n)\subseteq
  \cB_n$ for any $t\in \Rset$, $I\subseteq\Rset$ and $n \in \Zset$.
  Now, let $I:=[0,1]$, $J:=[1,2] =S_{1}(I)$ and $K:=[0,2]$. First,
  we show that there exists a non-zero $x \in \cA_I\cap\cB_n$ for some
  $n \neq 0$. If $\cA_I \cap \cB_n = \{0\}$ for all $n \neq 0$, then
  $\cA_I \subset \cA^\psi$ and furthermore $\cA = \bigvee_{n\in\Zset}
  S_n\cA_{I} \subseteq \cA^\psi$ by local minimality. But this
  contradicts the periodicity of $\psi$. Hence there exists a non-zero
  $x \in \cA_I \cap \cB_n$ for some $n \neq 0$. Since $S_1(x^*) \in
  \cA_J \cap \cB_{-n}$ is also non-zero, the $\Cset$-independence of
  $x$ and $S_1(x^*)$ implies that $x S_1(x^*) \in \cA_K \cap \cB_0$ is
  non-zero: $\psi(S_{1}(x)x^{*}xS_{1}(x^{*})) =
    \psi(xx^{*})\psi(x^{*}x) \neq 0$.  We are left to prove that $x
  S_1(x^*) \in \cA_K \cap \cB_0$ is $\psi$-orthogonal to $(\cA_I \cap
  \cB_0)\vee (\cA_J \cap \cB_0)$. It holds
  \begin{align}\label{obstruction-eqn}
        \psi(S_1(b)a xS_1(x^*))
     =  \psi(ax)\psi(S_1(b)S_1(x))
     =  0
  \end{align}
  for any $a,b \in \cA_I \cap \cB_0$, since $a\in \cB_0$ and $x\in
  \cB_n$ are $\psi$-orthogonal. Now, the commuting past/future
  structure ensures that $c \in (\cA_I\cap \cB_0) \vee (\cA_J\cap
  \cB_0)$ is approximated in the weak* topology by some sequence
  $(c_i)_{i \in \Nset}$ with terms of the form $c_i = \sum_l
  S_1(b_{i,l})a_{i,l}$, where $b_{i,l},a_{i,l} \in \cA_I\cap \cA^\psi$.
  Thus, (\ref{obstruction-eqn}) extends to $\psi(c xS_1(x^*))=0$ for
  any $c \in (\cA_I \cap \cA^\psi) \vee (\cA_J \cap \cA^\psi)$. This
  shows that $(\cA_I \cap \cA^\psi)_I$ is not locally minimal.\qedopt
\end{proof}

A concrete example for such a continuous Bernoulli shift is provided
by the CCR white noise in Example \ref{CCR-white-noise}.

\begin{Rem}\normalfont\label{remark:obstruction-classical}
  It is well-known in probability theory that multi-dimensional
  stochastic processes may lead to filtrations without local
  minimality. Let us next sketch the construction of such an example,
  starting from Gaussian white noise, as introduced in Example 
  \ref{Gaussian-wn}. Consider two independent
  generalized stochastic processes $(X_f^{(1)})_{f \in \srS}$ and
  $(X_g^{(2)})_{g \in \srS}$, realized on $(\srS^\prime \times
  \srS^\prime, \Sigma \times \Sigma, \mu \times \mu)$ with
  characteristic functional $C(f,g):= \e^{-1/2(\nm{f}^2+ \nm{g}^2)}$.
  A slight modification of the construction in Example
  \ref{Gaussian-wn} shows that the corresponding $\Cset$-expected 
  continuous Bernoulli shift is given by $(\cC \otimes \cC, 
  \psi_\mu \otimes \psi_\mu, (S_t \otimes S_t)_{t\in \Rset}, 
  (\cC_I \otimes \cC_I)_{I \in \cI})$, the tensor product of two
  Gaussian white noises. Now let 
  $\widetilde{E}$ be the conditional
  expectation associated to the sub-$\sigma$-algebra $\widetilde{\Sigma}$
  generated by $\set{X_f^{(1)}\cdot X_f^{(2)}}{f \in \srS}$. Then
  it can be checked that 
  $(\widetilde{E}(\cC\otimes \cC) \cap (\cC_I\otimes \cC_I))_{I \in \cI}$
  is a filtration without local minimality. From this follows that 
  the tensor product of two Gaussian white noises is compressed by
  the conditional expectation $\tilde{E}$ to a $\Cset$-centred
  continuous Bernoulli shift which fails to be locally minimal. 
\end{Rem}

\subsection{Examples from probability theory}
\label{subsection:CBS-examples} 
Tsirelson-Vershik's noises \cite[Definition 2d1]{Tsir04a} correspond
to $\Cset$-expected continuous Bernoulli shift which are locally
minimal and commutative (see Subsection
\ref{subsection:commutative-BS}).  These noises divide into two
classes: \emph{classical} and \emph{non-classical} (in the sense of
\cite[Definition 5c4]{Tsir04a}), or alternatively phrased:
\emph{type~I} and \emph{non-type~I} (following the analogy with
Arveson's product systems \cite{Arve03a}).

`Type I' examples  in probability theory correspond  to L\'evy
processes described by the L\'evy-Khinchin formula. Essentially, these
processes are combinations of Brownian motion and Poisson
processes. We will address the corresponding $\Cset$-expected
continuous Bernoulli shifts as white noises. The attribute `white'
emphasizes that the spectral density of L\'evy processes is constant
(see also Subsection \ref{subsection:CBS-WN} for the general
case). Notice that our usage of `white noise' differs from
\cite{Tsir04a} where it is reserved for `Gaussian white noise'(as
stated in Example \ref{Gaussian-wn}). Thus a `classical noise' therein
corresponds to a `white noise' herein.

\begin{Exmp}[Gaussian white noise]\normalfont\label{Gaussian-wn}
  Let $\srS\subset L_{\Rset}^{2}(\Rset)$ denote the space of all
  smooth rapidly decreasing real valued functions and $\srS'$ its
  dual, the space of tempered distributions. Consider the generalized
  stochastic process $(X_{f})_{f\in\srS}$ with $X_{f}(x'):=\langle
  f,x'\rangle$, $x'\in\srS'$, on the probability space
  $(\srS',\Si,\mu)$, where the measure $\mu$ is determined by the
  characteristic functional $C(f) := \e^{-\frac{1}{2}\nm{f}^{2}} =
  \int_{\srS'} \e^{\i X_{f}}\d\mu$ \cite{Hida,GeVi}. Let $\Si_I$
  be the $\si$-algebra generated by the functions $\e^{\i X_{f}}$ with
  $\text{supp}\,f\subseteq I$, where $I \in \cI$. From
  $C(f+g)=C(f)C(g)$ for all
  functions $f,g\in\srS$ with support in the disjoint intervals $I$
  and $J$ we obtain the independence of the random variables $X_{f}$
  and $X_{g}$ and hence the independence of $\Si_{I}$ and $\Si_{J}$.
  The characteristic functional is invariant under the right shift
  $\si_{t}$ on $\srS$. Consequently, $\mu$ is invariant under the dual
  action $s_{t}:= \si_{t}^*$ on $\srS'$. Now properties (o) to (iii)
  of Definition \ref{definition:CBSPS} can be checked. Hence, by
  Theorem \ref{theorem:TVN-CBS}, we get a locally minimal
  $\Cset$-expected shift $(\cC,\psi_{\mu},S,(\cC_{I})_{I\in \cI})$
  which is called \emph{Gaussian white noise}. Notice that Brownian
  motion $B_t\in L^2(\srS^\prime, \Sigma, \mu)$ is approximated by
  $(X_{f_n})_{n \in \Nset}$ with $f_n \to \chi_{[0,t]}$ in the
  $L^2$-norm and generates the sub-$\sigma$-algebra $\Sigma_{[0,t]}$
  for any $t >0$. 
\end{Exmp}

\begin{Exmp}[Poisson white noises]\normalfont\label{Poisson-wn}
  Let $N:=(N_{t})_{t\ge0}$ be the Poisson process with intensity
  $\la>0$. Then $N$ can be realized on a probability space
  $(\Om,\Si,\mu)$, where $\Om$ is the set of paths
  $\om\colon\Rset\to\Zset$, $\om(0)=0$, which are
  increasing and right continuous, and with left limits
  \cite{Prot95a}. We extend $N$ to negative times by $N_{-t}(\om): =
  \om(-t)$. Thus $\om(t)$ and $-\om(-t)$ count the jumps of $\om$ in
  $[0,t]$ resp. $(-t,0]$. The $\si$-algebras $\Si$ and
  $\Si_{[s,t]}$  are generated by the sets
  $\Om_{n}((s,t]) := \{N_{t}-N_{s}=n\}$, $s<t\in\Rset$,
  $n\in\Nset_{0}$ resp.\ $\Om_{n}((s',t'])$ with $s\le s'$ and $t'\le
  t$, $n\in\Nset_{0}$. The measure $\mu$ is given by
  $\mu(\Om_{n}((s,t])) := \la^{n}(t-s)^{n} \e^{-\la(t-s)}/n!$. The
  $\sigma$-algebras $\Si_{[s,t]}$ and $\Si_{[u,v]}$ are
  independent for disjoint intervals $[s,t]$ and $[u,v]$. Finally, a
  measure preserving shift $(s_{r})_{r\in\Rset}$ on $\Om$ is defined
  by $(s_{r}(\om))(t):= \om(t+r)-\om(r)$. Now, by Theorem
  \ref{theorem:TVN-CBS}, one can associate  canonically a locally
  minimal $\Cset$-expected shift $(\cP,\psi_{\mu},S,(\cP_{I})_{I\in
  \cI})$ to the Poisson process. It is called \emph{Poisson white
  noise}. Notice that the sub-$\sigma$-algebras of the filtration 
  are generated by increments of the Poisson process.
\end{Exmp}

\begin{Rem}\normalfont
  Tensor products of Poisson white noises and Gaussian white noises,
  compressed to the von Neumann subalgebra generated by a specified
  linear combination of the underlying Brownian motion and Poisson
  processes, give first examples of white noises coming from L\'evy
  processes. Moreover, (countable many) tensor products of Gaussian 
  white noise lead again to Gaussian white noises, now with 
  multiplicities. Furthermore, the amplification with a
  non-commutative probability space gives operator-expected white
  noises. We leave the details to the reader.
\end{Rem}

Tsirelson-Vershik's noises with a `non-classical' or `non-type~I' part
have no representation in  Fock spaces. The existence of such
intrinsically non-linear  random fields was revealed by A.~Vershik and
B.~Tsirelson \cite{TsVe98a}. For a detailed  survey on recent
developments and examples, as well as the close connection to
Arveson's non-type~I product systems, we refer the interested reader
to \cite{Tsir04a}. Prominent examples among these `non-classical'
noises are `black noises', as they are named by B.~Tsirelson in 
\cite{Tsir98a}. These `black noises' lead to Arveson's continuous 
product systems of Hilbert spaces of type~II$_0$ \cite{Arve03a}.

\begin{Exmp}[Black noises] \normalfont\label{black-noise}
  Examples of `black' noises, as stated in \cite{Tsir04a}, are locally
  minimal continuous Bernoulli shifts on probability spaces. They
  promote to $\Cset$-expected locally minimal continuous Bernoulli shifts
  according to Corollary \ref{corollary:TVN} and will also be addressed
  as $\Cset$-expected black noises. They are called `black' because they
  have only trivial additive shift cocycles (see Definition
  \ref{add-cocycle}). Thus, no `linear sensors' (as phrased by
  Tsirelson) exist to detect their color (see also the discussion of
  `whiteness' at the end of Subsection \ref{subsection:CBS-WN}). For
  further details on black noises and their construction we refer to
  \cite{Tsir04a} and the literature cited therein.
\end{Exmp}

\subsection{Examples from quantum probability theory}
\label{subsection:nCBS-examples} 
We continue with examples of continuous Bernoulli shifts coming from
quantum probability. In analogy to probability theory, we address the
`type~I' examples as `quantum white noises'.  Heuristically and
justified up to now by all known examples, these examples come from
quantum L\'evy processes which generate the  filtration of the quantum
white noise (see Definition \ref{def-white}). These quantum L\'evy
processes are provided by additive shift cocycles  (see Definition
\ref{add-cocycle}).  Already a rich source for the construction of
quantum white noises is provided by generalized Brownian motions which
are realized on (deformed) Fock spaces \cite{BoSp91a,BoSp94a,
BKS97a,GuMa02a, BoGu02a}.  Moreover, out the work of Anshelevich on
$q$-L\'evy processes \cite{Ansh04aTA} appear promising candidates for
further examples of quantum white noises.


\begin{Exmp}[CCR white noises]\normalfont \label{CCR-white-noise}  
  We continue the discussion of Example \ref{CCR-independence}. Let
  $\cK:= L^2(\Rset)$ and let $\cB_I$ be the von Neumann algebra
  generated by functions $f\in \cK$ with support in the interval $I$.
  Second quantization of the right shift on $L^2(\Rset)$ provides the
  shift $S$. It is elementary to check that $(\cB, \psi_\lambda, S,
  (\cB_I)_{I\in \cI})$ is a locally minimal $\Cset$-expected continuous
  Bernoulli shift, called CCR white noise. It is well-known that the
  von Neumann algebra of such a shift is a factor of type
  III$_{\la/(1+\la)}$  Notice that this shift has a
  commuting past/future and an enriched independence. 
\end{Exmp}

Multi-dimensional CCR white noises are just tensor products of
$\Cset$-expected CCR white noises.

\begin{Exmp}[Squeezed CCR white noises]
  \normalfont\label{squeezed-wn} More generally as done in Examples
  \ref{CCR-independence} and \ref{CCR-white-noise}, consider now the
  quasi-free state $\psi_{\lambda,c}$ on
  $\operatorname{CCR}(L^2(\Rset), \IM\skp{\cdot}{\cdot})$, given by
  $\psi_{\lambda,c}(W(f))= \exp{(-\frac14 q_{\lambda,c}(f)})$ with
  \[
        q_{\lambda,c}(f)
    :=  (2\lambda+1)\nm{f}^2 + 2 \RE(c\skp{f}{Jf}),
  \]
  where $J$ is the complex conjugation in $L^2(\Rset)$ and $\lambda >
  \sqrt{|c|^2+1/4} -1/2$ for some $c\in \Cset$. This state is non-gauge
  invariant for $c\neq 0$ and the corresponding Araki-Woods
  representation leads again to a $\Cset$-expected quantum white noise.
  This noise has important applications in quantum optics, where it is
  referred to as `squeezed white noise'. We refer the reader to
  \cite{HHKKR02a} for its construction and further references, moreover
  to \cite{GaZo00a} for its applications.
\end{Exmp}

\begin{Exmp}[CAR white noise]\normalfont\label{CAR-wn}
  Let $(\cB,\psi_{\la})$ be the probability space introduced in Example
  \ref{CAR-independence} with $\cK:=L^{2}(\Rset)$. For any interval
  $I\subseteq\Rset$ let $\cB_{I}$ be the subalgebra generated by the
  functions in $L^{2}(\Rset)$ with support in $I$. Second quantization
  of the right shift on $L^{2}(\Rset)$ provides a shift $S$, which
   fulfills obviously $S_{t}\cB_{I}=\cB_{I+t}$. This gives the
  $\Cset$-expected locally minimal shift
  $(\cB,\psi_{\la},S,(\cB_{I})_{I\in \cI})$, called
  \emph{$\operatorname{CAR}$ white noise}. Note that these shifts have
  an enriched independence, but they do not have a commuting past and
  future. Moreover, the von Neumann algebra of these shifts is a
  type~III$_{\la/(1-\la)}$ factor for $0 < \la <1/2$
  and a type~II$_1$ factor in the case $\la = 1/2$.
 
  Aside of amplification with a non-commutative probability space, the
  above construction can also be promoted to a $\cB_0$-expected shift as
  follows. Let $\cK := \cK_0 \oplus L^2(\Rset)$ ($\cK_0$ separable)
  and let $\cB_I$ be the von Neumann algebra generated by the closed
  subspace $\cK_I := \cK_0 \oplus L^2(I) \subseteq \cK$. Define the
  shift $S$ by second quantization of $\id \oplus s$. Here denotes $s$
  the right shift on $L^2(\Rset)$. It is again elementary to verify
  that $(\cB,\psi_{\la},S,(\cB_{I})_{I\in \cI})$ is a $\cB_0$-expected
  shift.
\end{Exmp}

\begin{Rem}\normalfont
  Let two $\Cset$-expected continuous Bernoulli shifts have von Neumann
  algebras of type~III$_\kappa$ resp.\ III$_{\tilde{\kappa}}$ ($0 <
  \kappa, \tilde{\kappa}<1$). If $\ln(\kappa)/\ln(\tilde{\kappa})$ is
  irrational then the tensor product of these two shifts leads to a
  $\Cset$-expected white noise with a von Neumann algebra of
  type~III$_1$ (see also \cite[Theorem 4.16]{Take03c} for approximately
  finite dimensional von Neumann algebras). Such examples are in
  particular provided by the CAR white noises. The authors doubt that
  examples of $\Cset$-expected continuous Bernoulli shifts with a von
  Neumann algebra of type~III$_0$ exist, but have yet not been
  successful in clarifying this point.
\end{Rem}

The following class of examples contains quantum white noise from free
probability theory, including a special case of amalgamated free
independence.

\begin{Exmp}[$q$-Gaussian white noises]
  \normalfont\label{q-white-noise} We continue the discussion of
  Example \ref{q-independence}. Let $\cK_\Rset:=
    L^2_\Rset(\Rset)$ and let $\cB_I$ be
  the von Neumann algebra, generated by the $q$-Gaussian processes
  $\set{\Phi(f)}{f\in\cK_I}$ with $\cK_I:= L^2_\Rset(I) \subseteq
  \cK_\Rset$. The shift $S$ is again obtained by second quantization of
  the right shift on $L^2(\Rset)$. From this one verifies that $(\cB,
  \tau, S, (\cB_I)_{I\in \cI})$ is a locally minimal $\Cset$-expected
  shift with enriched independence, but without commuting past and
  future. Such shifts are called \emph{$q$-Gaussian white noises} and
  their von Neumann algebras are factors of type~II$_1$ \cite{BKS97a}.
  
  Operator-expected shifts are produced by amplification. Similar to
  Example \ref{CAR-wn}, let us here present an alternative way in the
  multi-dimensional case. Put $\cK_I := \cK_0 \oplus L^2_\Rset(I,
  \cK_1)$, where $\cK_0, \cK_1$ are two real separable Hilbert spaces.
  The dimension of $\cK_1$ will give the multiplicity of the shift. As
  usual, $\cK_I$ is  identified canonically as a subspace of
  $\cK_\Rset$. Let $\cC_I :=\bigvee \set{\Phi(f)}{f\in\cK_{I}}$. The
  shift $S$ arises again through the second quantization of
  $\id_{\cK_0} \oplus s_t$, where $(s_t)_{t\in \Rset}$ is the right
  shift on $L^2(\Rset, \cK_1)$. Now one obtains the
  $\cC_{0}$-expected shift $(\cC_\Rset, \psi, S, (\cC_I)_{I\in
  \cI})$. Notice in the case $q=0$ that $\cC_\Rset$ is isomorphic to
  the free product of $\cB_{0}$ with the $\dim(\cK_1)$-fold
  free product of $\cB$, as stated at the beginning of the present
  example.
\end{Exmp}

\begin{Rem}\normalfont
  Further examples of $\Ao$-expected shifts arise from generalized
  Brownian motions by functors of white noise, as indicated already at
  the end of Section \ref{section:independence}.
\end{Rem}

\begin{Rem}\normalfont
  We close this section with a digression on `non-classical' (or
  `non-type~I') examples of continuous Bernoulli shifts which come from
  quantum probability. Presently, the existence of such examples is an
  open problem, if one insists on the local minimality of the
  filtration (see Subsection \ref{subsection:local-min-comp} for an
  example lacking local minimality). The analogy between
  Tsirelson-Vershik's noises and continuous Bernoulli shifts is
  evident, even more for the `classical or type~I' parts. Thus it is
  tempting to conjecture the existence of `quantum black noises', in
  the sense of continuous Bernoulli shifts. Such objects would enjoy a
  probabilistic interpretation and could provide the notion of an
  `intrinsically non-linear random quantum field'. Moreover, such
  examples would serve as a quantum probabilistic source for Arveson's
  continuous tensor product systems of type~II \cite{Arve03a}.
\end{Rem}

\section{Continuous GNS Bernoulli shifts}        \label{section:CGNSBS}

In this section we will develop the GNS representation theory of an
$\Ao$-expected continuous Bernoulli shift, as far as it will be needed
throughout this paper. The goal of this section is Definition
\ref{HMBS}, which provides our notion of an $\Ao$-expected continuous GNS
Bernoulli shift. It will be integral for the remaining Sections
\ref{section:cocycles} to \ref{section:nc-ln-exp}. We will start now
with a motivation.

Already a Gaussian white noise shows clearly that some of its most
interesting processes  will not sit in the  von Neumann algebra
$L^\infty(\srS',\Si,\mu)$:  Brownian motion, its generating stochastic
process, is contained in  the GNS Hilbert space
$L^2(\srS',\Si,\mu)$. Thus, starting with the Gaussian white noise
\[
      \big(L^\infty(\srS',\Si,\mu), \psi_\mu, S,
      \big(L^\infty(\srS',\Si_I,\mu_{|\Sigma_I})\big)_{
      I \in \cI}\big),
\]
we are led to the notion of a \emph{$\Cset$-expected continuous GNS
Bernoulli shift}, hereafter formulated in the example of Gaussian
white noise:
\begin{align}\label{GNSCBS-motivation}
      \big(L^2(\srS',\Si,\mu), \1, \overline{S},
      \big(L^2(\srS',\Si_I,\mu_{|\Sigma_I})\big)_{I \in \cI}\big).
\end{align} 
It captures all the necessary structure on the Hilbert space
level. Here is  $\1\in L^2(\srS',\Si,\mu)$ the cyclic and separating
vector coming from the GNS representation associated to
$\psi_\mu$. Moreover denotes  $\overline{S}=(\overline{S}_t)_{t\in
\Rset}$ the strongly continuous unitary group coming from the GNS
representation of the shift $S$.

Looking at this most familiar example, one realizes that
\eqref{GNSCBS-motivation} represents data which can also be captured
by a  \emph{homogeneous continuous product (system) of pointed Hilbert
spaces}, as it is introduced in \cite[Definitions 3c1 and
6d6]{Tsir04a}. It is \emph{pointed} because it has the \emph{unit}
$\1$.

In the general setting of an $\Ao$-expected continuous Bernoulli shift
$\CBS$ we can neither rely on the tensor product (independence), nor on
Hilbert spaces. But due to the notion of $\Ao$-independence (Definition
\ref{definition-independence}) and the rigid infrastructure of an
$\Ao$-expected continuous Bernoulli shift (Definition \ref{def:cbs}),
all means are at hand to construct a proper analogue. Technically
speaking and picking up the analogy, a \emph{continuous GNS Bernoulli
shift} will be a `homogeneous continuous commuting square system of
pointed Hilbert bimodules'. Here we will take a `lay-man's approach'
to construct the family of Hilbert bimodules $(\bimod{\cA_I})_{I\in
\cI}$ in $\cB(\cH_\psi)$, as we will need them to extend the
structure of an $\Ao$-expected continuous Bernoulli shift $\CBS$. Let
us outline the contents of the present section.

Subsection \ref{subsection:hm-spaces} introduces on an elementary level
the GNS construction of the Hilbert $\cA$-$\Ao$-bimodule $\bimod{\cA}$,
starting from an $\Ao$-expected probability space $(\cA,\psi)$. In
Subsection \ref{subsection:hm-morphisms} we will identify the morphisms
of an $\Ao$-expected probability space which extend to adjointable
bounded linear operators on the Hilbert bimodule $\bimod{\cA}$. In
Subsection \ref{subsection:hm-mult} we will introduce the product of
two $\Ao$-independent Hilbert bimodule elements. This product will be
crucial for the introduction of multiplicative shift cocycles in
Section \ref{section:cocycles}, as well as the development of a theory
of non-commutative It\^o integration in Section
\ref{section:nc-ito-integration}. Next we will extend in Subsection
\ref{subsection:CBS-GNS} the structure of an $\Ao$-expected continuous
Bernoulli shift $\CBS$ to its GNS Hilbert $\cA$-$\Ao$-bimodule
$\bimod{\cA}$ which leads to the notion of an \emph{$\Ao$-expected
  continuous GNS Bernoulli shift} \CBSE, finally stated Definition
\ref{HMBS} .

There are various approaches to Hilbert bimodules. For the
convenience of the reader, we provide in
Appendix~\ref{sec-appendix-hm} a short survey on Hilbert W*-modules,
as far as we will need them. Further references on  Hilbert modules
can also be found there.

\subsection{Hilbert bimodules of $\Ao$-expected probability spaces}
\label{subsection:hm-spaces}
In the following we will present the concrete construction of the
Hilbert $\cA$-$\Ao$ bimodule $\bimod{\cA}$, starting from the
$\Ao$-expected probability space $(\cA,\psi)$, occasionally also
denoted by the triple $(\cA,\psi,\Ao)$. It will be realized in
$\cB(\cH_\psi)$, the bounded linear operators on the GNS Hilbert space
$\cH_\psi$. We remind that $\cA \subseteq \cB(\cH_\psi)$ can be assumed
(see Section \ref{section:preliminaries}). Notice also that the
constructed Hilbert bimodule will be (isomorphic to) the GNS Hilbert
space $\cH_\psi$ if $\Ao \simeq \Cset$.

Let $\Eo$ be the conditional expectation from $(\cA,\psi)$ onto
$\Ao$. Its GNS representation defines an orthogonal projection $e_{0}$
such that $\Eo(x)e_{0}=e_{0} x e_{0}$ for any $x \in \cA$. Furthermore,
let $\Ho :=e_{0}\cH_{\psi}$. The vector space $\cA
e_{0}\subset\cB(\cH_\psi)$ is an \cA-\Ao-bimodule with left
multiplication
\begin{align} \label{left-mult}
       \cA\times\cA e_{0} \ni (y,xe_{0})
   \mapsto
       y x e_{0} \in \cA e_{0}
\end{align}
and right multiplication
\begin{align}\label{right-mult}
       \cA e_{0}\times\Ao\ni(x e_{0},a)
   \mapsto
       x e_{0} a
    =  xa e_{0} \in \cA e_{0}\,.
\end{align}
We introduce the $\Ao$-valued inner product
\begin{align}\label{inner-product}
       \cA e_{0}\times\cA e_{0} \ni (x,y) \mapsto \skpo{x}{y}
  &:=  x^{*}y
       \in\Ao
\end{align}
which turns $\cA e_{0}$ into a pre-Hilbert \Ao-module. From now on,
$\Ao e_{0}$ and $\Ao$ will be  identified canonically.

\begin{Defn}
  $\bimod{\cA}$ is the closure of $\cA e_{0}$ in the stop topology of
  $\cB(\cH_{\psi})$.
\end{Defn}

This closure can be identified as $\bimod{\cA} =\{\cA,e_{0}\}''e_{0}$,
thus it is elementary to see that the left and right multiplication
\eqref{left-mult} resp.\ \eqref{right-mult}, as well as the inner
product \eqref{inner-product}, extend to $\bimod{\cA}$. The
$\Ao$-valued `norm' and its induced norm on $\bimod{\cA}$ are denoted
by
\[
       \abo{x}
   :=  \skpO{x}{x}^{1/2}\quad \text{resp.}
       \quad
       \nmo{x}
   :=  \nm{\abo{x}}.
\]
In order to see the connection with the definition of a Hilbert
W*-module in Appendix \ref{sec-appendix-hm}, we note, that
$\cB(\cH_{\psi})e_{0}$ is canonically isomorphic to
$\cB(\Ho,\cH_{\psi})\subset\cB(\Ho\oplus\cH_{\psi})$.  Under this
isomorphism $\bimod{\cA}$ becomes a Hilbert W*-module in the sense of
Definition \ref{definition:Hilbert-W-module}. Identifying $\bimod{\cA}$
with its isomorphic image, we arrive at the following Lemma.

\begin{Lem}
  $\bimod{\cA}$ is an $\cA$-$\Ao$-bimodule and a Hilbert W*-module
  (over $\Ao$).
\end{Lem}

Notice that the left multiplication with an element in $\cA$ defines a
bounded linear operator on $\bimod{\cA}$ which is adjointable.

The strong operator (stop) topology on $\bimod{\cA}$ is generated
by the seminorms 
\[
      \d_{\xi}(x)
  :=  \nm{\abo{x}\xi},
      \qquad \xi\in\Ho
\]
and the $\si$-stop topology by
\[
      \d_{\phi}(x)
  :=  \phi(\skpo{x}{x})^{1/2},
      \quad \phi\in\cA_{0*}^{+}.
\] 
Moreover are the following continuity properties valid:
\begin{enumerate}
\item The map $ \bimod{\cA}\ni x\mapsto \skpo{x}{y} $ is
  weakly*-weakly* continuous for any $y\in\bimod{\cA}$;
\item The map $\bimod{\cA}\times \bimod{\cA}\ni (x,y)\mapsto
  \skpo{x}{y}$ is jointly continuous on bounded sets in the stop
  topology on $\bimod{\cA}$ and in the weak* topology on \Ao.
\end{enumerate}
Notice that $\cA$ is a dense subspace of $\bimod{\cA}$ in the stop
topology, since the algebra $\cA$ embeds contractively into
$\bimod{\cA}$ by the strongly continuous mapping $\cA\ni x\mapsto
x e_{0}$. From the separability of $\cA_{*}$ follows that $\bimod{\cA}$
has a separable predual. (see also Theorem \ref{self-dual}). Finally,
we note that Kaplansky's density Theorem ~\ref{Kaplansky} ensures that
elements in $\bimod{\cA}$ can be approximated in the stop topology by
bounded sequences in $\cA$.

\subsection{GNS representation of morphisms} 
\label{subsection:hm-morphisms}
We will characterize the elements of $\Mor\Apsi$ which extend to
adjointable bounded linear operators on $\bimod{\cA}$.
The following definition is needed for the formulation of Theorem
\ref{lap-morph}.

\begin{Defn}
  A morphism $T^*\in\Mor(\cA,\psi)$ is called the \emph{$\psi$-adjoint
    of $T\in\Mor(\cA,\psi)$} if $\psi(T^*(x)y)=\psi(xT(y))$ for any
  $x,y\in\cA$.
\end{Defn}

Notice that $T^{*}$  exists (uniquely) if and only if $T$ commutes with
the modular automorphism group $\sigma^\psi$ (see Theorem
\ref{psi-adj}).

\begin{Thm}\label{lap-morph}
  Let $T \in \Mor(\cA,\psi)$ commute with the modular automorphism
  group $\si^{\psi}$ and leave $\Ao$ pointwise fixed. Then the morphism
  $T$ has a unique extension to an adjointable bounded linear operator
  $\overline{T}$ on $\bimod{\cA}$ such that $\overline{\mbox{$T^{*}$}} =
  \overline{T}^{\,*}$. Moreover, $T\mapsto\overline{T}$ is pointwise
  weakly*-weakly* continuous and pointwise strongly-strongly
  continuous.
\end{Thm}

\begin{proof}[Proof of Theorem \ref{lap-morph}.] 
  The morphism $T$ defines by $\overline{T}x\Om:=T(x)\Om$, $x\in\cA$, a
  contraction which extends to the GNS Hilbert space $\cH_\psi$.
  Moreover, we know from Theorem \ref{psi-adj} that its $\psi$-adjoint
  $T^*\in \Mor(\cA,\psi)$ exists uniquely  and thus also extends to a
  contraction $\overline{T^*}$ such that $\overline{T^*}=
  \overline{T}^*$.
 
  We will need that $\Ao$ is contained in the fixed point algebra of
  $T^{*}$: Since $\Ao$ is contained in the fixed point algebra of $T$,
  we conclude $T(xa)=T(x)a\,$ and $T(ax)=aT(x)\,$ for any $x\in\cA$ and
  $a\in\Ao$ (see also \cite{Kuem84aUP, Robinson1}). Furthermore is
  $\psi(T^{*}(a)x) = \psi(aT(x)) = \psi(T(ax)) = \psi(ax)$ for any
  $x\in\cA$, $a\in\Ao$. Thus $\Ao$ is also contained in the fixed point
  algebra of $T^{*}$. We notice for later arguments that $T^{*}$
  satisfies also the conditions of the theorem.
  
  We conclude next for any $x,y\in\cA$
  \begin{align*}
        \overline{T}x e_0 y\Om
     =  \overline{T}xE_{0}(y)\Om
     =  T(xE_{0}(y))\Om
     =  T(x)E_{0}(y)\Om
     =  T(x) e_0 y\Om\,,
  \end{align*}
  and consequently $\overline{T}x e_0=T(x) e_0$. But this implies $
  \overline{T}{\bimod{\cA}}\subseteq{\bimod{\cA}}$ and
  \[
        \abO{\overline{T}x e_0}^{2}
     =   e_0(x^{*}\overline{T}^{\,*}\overline{T}x)e_0 
   \le  \nm{\overline{T}}^{2}\abO{x}^{2}
  \] 
  for any $x\in\cA$. From inequality \eqref{A-linear} follows now
  $\overline{T} \in \cB({\bimod{\cA}})$. Furthermore we conclude with
  Corollary \ref{BE-LE} that $\overline{T}$ is adjointable. The
  morphism $T^{*}$ satisfies again all conditions of the theorem and
  therefore extends also to an adjointable operator
  $\overline{\mbox{$T^{*}$}}$. Thus we can conclude
  \begin{align*}
        \psi(a\skpo{x}{\overline{T}y}) 
   & =  \psi(aE_{0}(x^{*}T(y)))
     =  \psi(ax^{*}T(y))                                             \\
   & =  \psi(T^{*}(ax^{*})y)
     =  \psi(aE_{0}(T^{*}(x^{*})y))
     =  \psi(a\skpo{\overline{\mbox{$T^{*}$}}x}{y})
  \end{align*}
  for any $a\in\Ao$, $x,y\in\cA$. But this implies 
  $\overline{\mbox{$T^{*}$}}=\overline{T}^{\,*}$.
  
  Assume that $(T_{\al})_{\al\in I}\in\Mor\Apsi$ converges pointwise to
  $T\in\Mor\Apsi$ in the weak* topology. Moreover, assume that each
  $T_{\al}$ and $T$ satisfy the assumptions of the theorem. Then it
  follows that $(\overline{T}_{\al})_{\al\in I}$ maps the unit ball
  of $\bimod{\cA}$ into itself.
  
  The map $T\mapsto T^{*}$ is pointwise weakly* continuous, since the
  weak* topology on $\cA_{1}$ is induced by the family of seminorms
  $\set{\ab{\psi(y\,\cdot\,)}}{y\in\cA}$. It follows
  \[
        \phi(\skpo{y}{(\overline{T}_{\al}-\overline{T})x})
     =  \phi(\skpo{(T_{\al}^{*}-T_{}^{*})(y)}{x})
        \xrightarrow{\al} 0
  \]
  for any $\phi\in\cA_{0*}^{+}$, $y\in\cA_{1}$ and
  $x\in{\bimod{\cA}}_{1}$, $\nmo{x}\le1$.
  This implies the pointwise weak*-weak*
  continuity of the extension $\overline{T}$, since the weak* topology
  on the unit ball of $\bimod{\cA}$
  is induced by the family of seminorms
  $\set{\ab{\phi(\skpo{y}{\cdot})}}{\phi\in\cA_{0*}^{+}, y\in\cA_{1}}$.
  
  Now, let us assume that $F_{\al}^{}:=T_{\al}^{}-T$ converges to $0$
  in the pointwise stop-topology. Then it follows for any $x,y\in\cA$
  that 
  \[
         \psi(y(T_{\al}^{*}\circ T_{\al}^{} - T^{*}\circ T)(x))
     =   \psi(F_{\al}^{}(y)T_{\al}^{}(x)) + \psi(T(y)F_{\al}^{}(x)).
  \]
  In other words, $T_{\al}^{*}\circ T_{\al}^{}$ converges to
  $T^{*}_{}\circ T$ in the pointwise weak* topology. Since $T^{*}$ is
  weakly*-weakly* continuous, we conclude that $F_{\al}^{\,*}\circ
  F_{\al}^{}$ converges to $0$ in the pointwise weak* topology on \cA.
  With $x\in{\bimod{\cA}}$ and $\phi\in\cA_{0*}^{+}$, we conclude
  further for any $y\in\cA$ that
  \[  
        \phi(\skpo{y}{\overline{F}_{\al}^{\,*}\,
        \overline{F}_{\al}^{}x})
     =  \phi(\skpo{F_{\al}^{*}\circ F_{\al}^{}(y)}{x})
        \xrightarrow{\al} 0.
  \]
  Thus one has $\lim_{\al} \overline{F}_{\al}^{\,*} \,
  \overline{F}_{\al}^{} x = 0$ in the weak* topology. The pointwise
  strong convergence of $\overline{T}_{\al}$ to $\overline{T}$ is an
  immediate consequence. \qedopt
\end{proof} 

\begin{Notation}\normalfont
  The argument $x$ of a morphism $T \in \Mor\Apsi$ will always be put
  in parenthesis, in contrast to the argument of $\overline{T}$, its
  extension to an adjointable bounded linear operator on
  ${\bimod{\cA}}$. This distinguishes morphisms acting on $\cA$ and
  bounded linear operators on ${\bimod{\cA}}$ in most cases
  sufficiently; thus we will denote them from now on by the same symbol
  $T$ to lighten the notation.
\end{Notation}

Notice finally that a conditional expectation $E \in \Mor(\cA,\psi)$
extends to an orthogonal projection on ${\bimod{\cA}}$. In particular,
the extension of the conditional expectation $\Eo$ onto $\Ao$
satisfies $\Eo x = \skpo{\1}{x}$ $(= e_{0} x e_0)$, where
$x\in{\bimod{\cA}}$.

\subsection{The product of $\Ao$-independent elements}
\label{subsection:hm-mult}
Let $\cB$ be a subalgebra of $(\cA,\psi,\Ao)$ with $\Ao\subseteq\cB$,
such that the conditional expectation $E_{\cB}\in \Mor(\cA,\psi)$
exists. Since $E_{\cB}$ is a morphism which commutes with the modular
automorphism group and which leaves $\Ao$ pointwise fixed, it extends
to an orthogonal projection on ${\bimod{\cA}}$. One verifies  easily 
that $E_{\cB}{\bimod{\cA}}= \bimod{\cB}$, where $\bimod{\cB}$
corresponds to $(\cB,\psi,\Ao)$.

\begin{Prop}\label{new-product}
  Let \cB and \cC be two von Neumann subalgebras of the probability
  space $(\cA,\psi,\Ao)$ such that the conditional expectations
  $E_{\cB}$ resp.\ $E_{\cC}$ exist.
  
  If \cB and \cC are $\Ao$-independent, then
  \begin{align}\label{new-pr}
        \bimod{\cB}\times \bimod{\cC} \ni (x,y)
   \mapsto 
        x y \in \bimod{\cB\vee\cC}
  \end{align}
  defines a product which extends the left multiplication 
  \[
         \cB \times {\bimod{\cA}} \ni  (x,y) 
    \mapsto 
         xy \in{\bimod{\cA}}. 
  \]
  This product is jointly \nmot-continuous. If its first component is
  \nmot-bounded, it is also jointly $\si$-stop - $\si$-stop-continuous.
  The product satisfies
  \begin{align}\label{factor-lap}
        \skpo{x_{1}y_{1}}{x_{2}y_{2}} 
     =  \skpo{y_{1}}{\skpo{x_{1}}{x_{2}}\,y_{2}}\,,
  \end{align}
  for any $x_{1},x_{2}\in \bimod{\cB}$ and $y_{1},y_{2}\in
  \bimod{\cC}$. Moreover, the module property of the conditional
  expectations $E_{\cB}$ resp.\ $E_{\cC}$ extends to
  \begin{align}\label{P-new-product}
        E_{\cB}\,xy
   & =  x\,\Eo y 
        \qquad\text{and}\qquad
        E_{\cC}\,xy 
     =  (\Eo x)y \,,
  \end{align}
  for any $x\in \bimod{\cB}$ and $y\in \bimod{\cC}$.
\end{Prop}
 
We emphasize that in general the product \eqref{new-pr} is not minimal:
The $\nmot$-closure of the product $\bimod{\cB}\, \bimod{\cC}$, where
$\cB$ and $\cC$ are $\Ao$-independent, may be contained properly in
$\bimod{\cB\vee\cC}$. Consider for example $\Cset$-free independence as
stated in Example \ref{free-independence}.

\begin{proof}
  The equation \eqref{factor-lap} reduces for $x_{i}\in\cB$ and
  $y_{i}\in\cC$ ($i=1,2$) simply to the equation in Proposition
  \ref{independence-properties}\spc(\ref{independence-properties-ii}).
  It is still valid for $y_{i}\in \bimod{\cC}$. This is 
  concluded immediately from the approximation of $y_{i}$ by bounded 
  sequences in
  $\cC$, from the continuity of the left multiplication by elements in
  $\cB$, and finally the continuity of the inner product.
  Next, we approximate $x\in \bimod{\cB}$ by a sequence 
  $({x_n})_{n\in\Nset}$, 
  with $x_n \in \cB$, in the strong (operator)
  topology on $\bimod{\cB}$
  and show that for $y\in \bimod{\cC}$ the element
  \begin{align}\label{def-new-product}
        xy  
    :=  \lim_{n\to\infty} x_{n}y\,
  \end{align}
  is well-defined. We verify the claimed convergence in the strong
  topology and the independence from the choice of the approximating
  sequence $(x_{n})_{n\in\Nset}$ as follows: for another approximating
  sequence $(\tilde{x}_{n})_{n\in\Nset} \subset \cB$ of $x$ in the
  strong topology on $\bimod{\cA}$ and for any $\xi\in\Ho$ it is:
  \begin{align*}
        \nm{\abo{(x_{n}-\tilde{x}_{m})y}\xi}^{2} 
   & =  \bskp{(x_{n}-\tilde{x}_{m})y\xi}{
        (x_{n}-\tilde{x}_{m})y\xi}                                   \\
   & =  \bskp{y\xi}{ \abO{x_{n}-\tilde{x}_{m}}^{2}y\xi}
     =  \d_{\phi}(x_{n}-\tilde{x}_{m})^{2}
        \xrightarrow{m,n\,\to\,\infty} 0\,.
  \end{align*}
  with $\phi:=\skp{y\xi}{\cdot y\xi}\in {\Ao}_{*}^{+}$. Here we used 
  that stop-convergent sequences are also $\si$-stop convergent.
  One observes that $xy$ is just the usual
  left multiplication whenever $x$ is in $\cB$.  Consequently, the
  equation \eqref{factor-lap} follows, since the inner product is
  continuous.
  
  We are left to prove the continuity of the product. Let
  $x=\lim_{\al}x_{\al}$ in the stop topology, where $x_{\al}\in
  \bimod{\cB}$ and $\nmo{x_{\al}}\le M$. Moreover, let
  $y=\lim_{\be}y_{\be}$ in the stop topology, where $y_{\be}\in
  \bimod{\cC}$. We observe
  \begin{align*}
        x_{\al}y_{\be}-xy 
     =  (x_{\al}-x)(y_{\be}-y) + (x_{\al}-x)y + x(y_{\be}-y)\,.
  \end{align*}
  For the first summand we conclude for any 
  $\phi\in\cA_{0*}^{+}$
  \begin{align*}
        \d_{\phi}((x_{\al}-x)(y_{\be}-y))^{2}
   & =  \phi(\skpo{y_{\be}-y}{\abO{x_{\al}-x}^{2}\,(y_{\al}-y)})     \\
   &\le \nmO{x_{\al}-x}^{2}\, \d_{\phi}(y_{\be}-y)^{2}
    \le 4M^{2}\d_{\phi}(y_{\be}-y)^{2}
        \xrightarrow{\be} 0\,.
  \end{align*}
  Similarly, we proceed with the third summand. The second summand
  delivers the expression $\phi(\skpo{y}{\abO{x_{\al}-x}^{2}\,y})$
  which converges to zero, since the map ${\bimod{\cA}}\ni
  z\mapsto\phi(\skpo{y}{\abO{z}^{2}\,y})$ is continuous in the stop
  topology. The \nmot-continuity is shown easily. Finally, we conclude
  \eqref{P-new-product} from the defining equation
  \eqref{def-new-product} of the product $xy$ and the strong continuity
  of the conditional expectations $E_{\cB}$ resp.\ $E_{\cC}$.\qedopt
\end{proof}

\begin{Rem}\normalfont 
  If \cB and \cC are $\Ao$-independent, we obtain for $x\in
  \bimod{\cB}$ and $y\in \bimod{\cC}$
  \begin{align}\label{factor-2}
        \skpo{x}{y}
     =  \skpo{x}{\1}\skpo{\1}{y}\,.
  \end{align}
  This is evident by choosing $x_{1}=y_{2}=\1$ and $y_{1}=x$, $x_{2}=y$
  in equation \eqref{factor-lap}. Similarly,
  \begin{align}\label{factor-2a}
        \Eo xy 
     =  \skpo{\1}{xy} 
     =  \skpo{\1}{x} \skpo{\1}{y} 
     =  \Eo x \Eo y.   
  \end{align}
  Notice that the module property of $\Eo$ insures $\Eo(x \Eo y) =
  (\Eo x)(\Eo y)$. The identities \eqref{factor-2} and
  \eqref{factor-2a} feature the factorization property \eqref{factor}
  in the language of Hilbert modules.
\end{Rem}

\begin{Rem}\normalfont
  The $n$-tuple product $x_1x_2\cdots x_n$ of elements $x_i \in
  \bimod{\cB_i}\subset \bimod{\cA}$ ($i=1,2,\ldots,n$) is
  well-defined and associative whenever $\cB_j$ and
  $\bigvee_{i=1}^{j-1}\cB_i$ are $\Ao$-independent for all
  $j=2,\ldots,n$. Moreover, it is $\Eo x_1x_2x_3\cdots x_n
  = \Eo x_{1}\Eo x_{2}\Eo x_{3}\cdots \Eo x_n$.
\end{Rem}

The following notion of $\Ao$-independence will be used for
${\bimod{\cA}}$.

\begin{Defn}
  $\bimod{\cB}$ and $\bimod{\cC}$ in $\bimod{\cA}$ 
  are called \emph{$\Ao$-independent} if $\cB$ and $\cC$ are
  $\Ao$-independent. Two elements $x, y \in \bimod{\cA}$ are
  $\Ao$-independent if they are, respectively, elements of two
  $\Ao$-independent $\bimod{\cB}$ and
  $\bimod{\cC}$.
\end{Defn}

All properties of commuting squares, as stated in Proposition
\ref{independence-properties}, carry over to ${\bimod{\cA}}$.

\begin{Prop}
  Under the assumptions of Proposition \ref{new-product}, the
  following conditions are equivalent:
  \begin{enumerate}
  \item $\bimod{\cB}$ and $\bimod{\cC}$ are $\Ao$-independent;
  \item $ \skpo{x_{1}y_{1}}{x_{2}y_{2}} =
    \skpo{y_{1}}{\skpo{x_{1}}{x_{2}}\,y_{2}}$ for any $x_1,x_2 \in
    \bimod{\cB}$ and $y_1,y_2 \in \bimod{\cC}$;
  \item $E_\cB \bimod{\cC} = \Ao$;
  \item $E_\cB E_\cC = \Eo$.
  \end{enumerate}
\end{Prop}

\begin{proof}
  The stated equivalences follow from Proposition
  \ref{independence-properties} and \eqref{def-new-product}, if one
  approximates elements in $\bimod{\cA}$ by sequences in the underlying
  von Neumann algebra $\cA$.\qedopt
\end{proof}

\subsection{Continuous GNS Bernoulli shifts}
\label{subsection:CBS-GNS}
We will use the results of the previous subsections to extend the
structure of  an $\Ao$-expected  continuous Bernoulli shift \CBS.

\begin{Notation}\normalfont
We lighten the notation of Hilbert bimodules and let  $\Lapneu :=
\bimod{\cA}$ and $\LapIneu := \bimod{\cA_I}$.
\end{Notation}

\begin{Prop} \label{proposition:CGNSBSa}
  Let the $\Ao$-expected continuous Bernoulli shift $\CBS$ be given. The
  filtration $(\cA_{I})_{I\in \cI}$ induces the family of Hilbert
  bimodules $(\LapIneu)_{I\in \cI}$ such that
  \begin{enumerate}
  \item[(o)] $\Lapneu$ is the closure of $\set{\LapIneu}{I \in \cI
    \text{ bounded}}$ in the stop topology.
  \item \label{proposition:CGNSBSai} $\LapIneu \subseteq \LapJneu $
    whenever $I \subseteq J$;
  \item \label{proposition:CGNSBSaii} $\bigcup\set{\LapIneu}{I \in
    \cI\text{\ bounded}}$ is strongly dense in ${\Lapneu}$.
  \item \label{proposition:CGNSBSaiii} $\LapIneu$ and $\LapJneu$
    are $\Ao$-independent whenever $I\cap J =\emptyset$;
  \end{enumerate}
\end{Prop}

\begin{proof}
  This follows from the results in the present section.
\end{proof}

Notice that due to Proposition \ref{new-product}, the product $xy$ of
independent elements $x$ and $y$ is well-defined in $\Lapneu$.

The shift $S_t$ satisfies the assumptions of Theorem \ref{lap-morph}.
Thus $S_t$ extends to an identically denoted operator $S_t$ on
${\Lapneu}$ for any $t\in \Rset$. 

\begin{Prop}\label{proposition:CGNSBSb}
  Let the $\Ao$-expected continuous Bernoulli shift $\CBS$ be given.
  Then the shift $S$ on the Hilbert bimodule $\Lapneu$ enjoys the
  following properties:
  \begin{enumerate}
  \item \label{proposition:CGNSBSbi} $S=(S_{t})_{t\in\Rset}$ is a
    unitary group on ${\Lapneu}$, which is pointwise continuous in the
    stop topology;
  \item \label{proposition:CGNSBSbii} the fixed point space of $S$
    in ${\Lapneu}$ is $\Ao$, in particular it is $S_t \Eo = \Eo$
    for any $t \in \Rset$;
  \item \label{proposition:CGNSBSbiii}
    $S_t$ maps $\LapIneu$ onto $\cE_{I+t}$ for     
    any interval $I \subseteq \Rset$ and $t\in\Rset$.
  \end{enumerate}
\end{Prop}

\begin{proof}
  This is clear.
\end{proof}

We will  use frequently the following result.

\begin{Cor}\label{Corollary:cexp-pastfuture}
  Let a $\Ao$-continuous Bernoulli shift be given.
  $(E_{(-\infty,t]})_{t\in \Rset}$ and $(E_{[t,\infty)})_{t\in \Rset}$
  are pointwise continuous in the stop topology on $\Lapneu$.
\end{Cor}

\begin{proof}
  This is concluded from Proposition
  \ref{proposition:CGNSBSa}\spc(\ref{proposition:CGNSBSaiii}) and
  \ref{proposition:CGNSBSb}\spc(\ref{proposition:CGNSBSbi}), similar as
  in Lemma \ref{P-cont}.
\end{proof}

In this paper we refrain from an axiomatization of these structures in
the language of Hilbert bimodules and work  with the following
definition.

\begin{Defn} \label{HMBS}
  The quadruple $(\Lapneu, \1, S, (\LapIneu)_{I\in \cI})$ is called an
  \emph{$\Ao$-expected continuous GNS Bernoulli shift}, whenever it is
  constructed from an $\Ao$-expected continuous Bernoulli shift $\CBS$
  as stated above.
\end{Defn}

For simplicity, we will  refer occasionally to continuous GNS Bernoulli
shifts just as shifts.

\begin{Rem}\normalfont
  (i) Since a $\Cset$-expected locally minimal commutative continuous
  Bernoulli shift corresponds to a Tsirelson-Vershik noise, it is
  evident that both cover the same class of of continuous product
  systems of pointed Hilbert spaces (see \cite[Definitions 3c1 and
  6d6]{Tsir04a} for a definition).
  
  (ii) $\Ao$-expected locally minimal continuous GNS Bernoulli shifts
  with a commuting past/future are in close contact with product
  systems of Hilbert modules, as they are considered in \cite{MuSo02a}
  and \cite{BhSk00a,BBLS04a}.
\end{Rem}

\section{Cocycles of continuous (GNS) Bernoulli shifts}
\label{section:cocycles}

In this section we will introduce and investigate additive and
multiplicative cocycles for $\Ao$-expected (non-commutative) continuous
Bernoulli shift. These cocycles are a non-commutative version of
operator-valued L\'evy processes and thus enjoy a probabilistic
interpretation. The \emph {unital} (or \emph{unitary}) \emph{shift
cocycles} will take over the role of stationary multiplicative
$G$-flow in Tsirelson's theory \cite{Tsir04a}, or of units in
Arveson's product systems of Hilbert spaces \cite{Arve03a}. As to
expect from this parallel, an \emph{additive shift cocycle} will play
in the present operator algebraic setting the part of a stationary
additive $G$-flows \cite{Tsir04a} resp.\ of the `logarithm' of a unit
\cite{Arve03a}, or of an `addit' in Bhat-Srinivasan's sum systems
\cite{BhSr04a}. Notice that the present approach does not cover
general $G$-flows of Tsirelson' theory; our (semi-)group $G$ is 
restricted to be the additive (resp.\ multiplicative semi-)group $\Ao$.

Our main result is a \emph{bijective correspondence} between
\emph{unital shift cocycles} and \emph{additive shift cocycles},
subject to a continuity condition resp.\ a structure equation  (see
Theorems \ref{theorem:ln-exp} and \ref{main-theorem}).  This
correspondence is based on the construction  of \emph{non-commutative
logarithms} and \emph{exponentials} which is the subject of  Section
\ref{section:nc-ln-exp}. Moreover, it rests on \emph{non-commutative
It\^o integration}, which will be developed in 
Section \ref{section:nc-ito-integration}. 

The correspondence can also be regarded as a non-commutative stochastic
extension of Stone's theorem on unitary groups with bounded generators.
Since our main result is formulated in the language of Hilbert
bimodules, it can moreover be
seen as a correspondence between a concrete Stinespring representation
of uniformly continuous Markovian (or $CP_0$-)semigroups on a von
Neumann algebra and the Christensen-Evans form of its generator. From 
the stochastic point of view it provides a non-commutative analogue of the
correspondence between square integrable L\'evy processes with values
in the (additive group) $\Ao$ and `unital' square integrable
multiplicative L\'evy processes with values in the (multiplicative
semigroup) $\Ao$. Obviously, there exists also a parallel to subfactor
theory in the (one- or two-sided) `discrete time' case. We have not
developed enough material here to point out this connection explicitely
and thus postpone further details in this direction.

Our main result insures already the abstract form of a \emph{bijective
correspondence} between \emph{unitary shift cocycles} and
\emph{additive shift cocycles}. But now the latter one's are satisfying
a stronger structure equation. An explicit rigorous presentation of
this correspondence is beyond our limits in this paper and will be
presented in sequel publications. But already the abstract form of this
correspondence provides enough information to introduce the notion of a
\emph{non-commutative white noise} in Definition \ref{white}. We show
that the `non-commutative white noise part' can always be extracted
from an $\Ao$-expected continuous Bernoulli shift. This is in parallel to
results that one can extract the `classical part' from a
Tsirelson-Vershik noise \cite{Tsir04a}, or that one can extract the
`type~I part' of Arveson's tensor product systems \cite{Arve03a}.
   
Let us summarize the contents of this section. In Subsection
\ref{subsection:uc-cbs} we introduce unitary cocycles for an
$\Ao$-expected continuous Bernoulli shift and show that they lead as
usual to $CP_0$- (or Markovian) semigroups by compressions. This gives
in particular rise to stationary quantum Markov processes as they are
considered in \cite{Kuem85a}, but now in continuous time. Subsections
\ref{subsection:uc-cgnsbs} and \ref{subsection:ac-cgnsbs} contain the
definitions of unital cocycles resp.\ additive cocycles for
$\Ao$-expected continuous GNS Bernoulli shifts and some of their
elementary properties from which we will make frequent use. In
Subsection \ref{subsection:correspondence} we present with Theorems
\ref{theorem:ln-exp} and \ref{main-theorem} the bijective
correspondence between additive and unital cocycles. The first
theorem states this correspondence in an abstract manner, the second in
a concrete manner which allows to identify the non-commutative
logarithms and exponentials. We also include there some immediate
implications of this correspondence: If the cocycles are weakly*
differentiable, one recovers Stone's theorem from the
correspondence. Moreover allows the correspondence to identify the
Christensen-Evans generator of a uniformly continuous $CP_0$-semigroup
in terms of additive cocycles. We proceed in Subsection
\ref{subsection:CBS-WN} with an existence result about another
correspondence, this time between unitary cocycles and additive
cocycles. This result puts us into the position to introduce the
notion of an $\Ao$-expected non-commutative white noise (Definition
\ref{white}) and we will study some of its properties. Finally, we
apply our main result in Subsection
\ref{subsection:examples-correspondence} to the examples of Gaussian
and Poisson white noise, CCR, CAR and q-white noises.
 
Throughout we will assume that a fixed $\Ao$-expected continuous
Bernoulli shift \CBS  and its continuous GNS Bernoulli shift \CBSE are
given.

\subsection{Multiplicative cocycles of continuous Bernoulli shifts}
\label{subsection:uc-cbs}
We introduce unitary cocycles. They are a non-commutative version
of L\'evy processes, taking values in unitary operators. 

\begin{Defn}\label{unitary-cocycle}
  Let $\CBS$ be an $\Ao$-expected shift. A family of unitaries $u =
  (u_t)_{t\ge 0}$ in $\cA$ is called a \emph{unitary cocycle} if
  for any $s,t \ge 0$
  \begin{enumerate} 
  \item\label{unitary-cocycle-i}    $t \mapsto u_t$ is weakly*
                                    continuous; 
  \item\label{unitary-cocycle-ii}   $u_t \in \cA_{[0,t]}$;
  \item\label{unitary-cocycle-iii}  $u_{s+t}= S_t(u_s)u_t$. 
  \end{enumerate}
  The unitary cocycle $u$ is called \emph{trivial} if $u \subseteq
  \Ao$. The set of all ($\nmot$-continuous) unitary cocycles is
  denoted by $\Cm(\cA,\boldsymbol{\cdot}\,)$ (resp.\ 
  $\Cmo(\cA,\boldsymbol{\cdot}\,)$).
\end{Defn}
Putting $u_t = S_t^{}(u_{-t}^*)$ for $t \le 0$, the unitary cocycle
$u$ extends to $\Rset$. For a given $\Ao$-expected shift, a unitary
cocycle $u$ defines via $T_{t}^{}(x) = u_t^*S_t^{}(x)u_t^{}$, $x
\in \cA$ a pointwise weakly* continuous group of automorphism $T =
(T_t)_{t\in \Rset}$ on $\cA$. As usual, the compression $R=\Eo T \Eo$
defines a pointwise weakly* continuous unital semigroup of completely
positive contractions, also called a \emph{Markovian or 
$CP_0$-semigroup on $\Ao$}.

Notice that $T \subset \Aut(\cA,\psi)$ if and only if $u \subset
\cA^\psi$. In such a situation, the quadruple $(\cA, \psi,T; \Ao)$
defines an $\Ao$-valued stationary quantum Markov process in the sense
of \cite{Kuem85a}. The Markovian semigroup $R=(R_t)_{t\ge0}$ is again 
obtained by the compression  $R_t =\Eo T_t \Eo$ and leaves $\psi$ 
invariant.

\subsection{Multiplicative cocycles of continuous GNS Bernoulli shifts}
\label{subsection:uc-cgnsbs}
In the framework of continuous GNS Bernoulli shifts, the role of a
unitary cocycle is filled by a  \emph{unital} cocycle.

\begin{Defn}\label{unit-cocycle} 
 Let $\CBSE$ be an $\Ao$-expected continuous GNS Bernoulli shift. A
  \emph{unital cocycle} is a weakly* continuous family
  $u:=(u_{t})_{t\ge0} \subset{\Lapneu}$ such that for any $s,t \ge 0$
  \begin{enumerate}
  \item\label{unit-cocycle-i} 
    $\abo{u_{t}}=\1$;                                 \hfill(unitality)
  \item\label{unit-cocycle-ii} 
    $u_{t}\in \LapXneu{[0,t]}$;                     \hfill(adaptedness)
  \item\label{unit-cocycle-iii} 
    $u_{t+s}=(S_{t}u_{s})u_{t}$.              \hfill(cocycle identity)
  \end{enumerate}
  The unital cocycle $u$ is called \emph{trivial} if $u \subseteq
  \Ao$. 
  
  The set of all ($\nmot$-continuous) unital cocycles is denoted by
  $\Cm(\cE,\boldsymbol{\cdot}\,)$ (resp.\ 
  $\Cmo(\cE,\boldsymbol{\cdot}\,)$).
\end{Defn}

Notice that every unitary cocycle gives a unital cocycle. We
remark also  that the cocycle identity (\ref{unit-cocycle-iii}) relies
on the fact that  the product of the independent elements
$S_{t}u_{s}\in \LapXneu{[t,t+s]}$ and $u_{t}\in \LapXneu{[0,t]}$ is
well-defined by Proposition \ref{new-product}.

A unital cocycle $u$ is trivial if and only if $u$ is a strongly
continuous semigroup of isometries in $\Ao$. Notice moreover that a
$\nmot$-continuous unital cocycle is trivial if and only if $u$
is a uniformly continuous semigroup of unitaries in $\Ao$. 
  
If $\Ao \simeq \Cset$, then a
unital cocycle is a family of unit vectors in the GNS Hilbert space
$\cH_{\psi}$.

From the definition of a unital cocycle follows already
$u_{0}=\1$. Indeed, it is $u_{0}\in\Ao$ by adaptedness and
$u_{0}^{}=u_{0}^{2}$ by the cocycle identity. Moreover,
$u_{0}^{*}u_{0}^{}=\1$ is satisfied by unitality. From this we
conclude $u_{0}-\1=u_{0}^{*}u_{0}^{}(u_{0}^{}-\1)=0$. As in the case
of unitary cocycles, the $\Ao$-independence ensures that the
compression of a unital cocycle yields contractive semigroups as
follows.

\begin{Prop} \label{semigroup} 
  Let $u$ be a unital cocycle.
 \begin{enumerate}
 \item\label{semigroup-i} The compression $A := (\Eo u_{t})_{t\ge0}$
  is a strongly continuous semigroup of contractions in $\Ao$.  
  Moreover, the following are equivalent:
   \begin{enumerate}
    \item[(a)] the unital cocycle $u$ is \nmot-continuous;
    \item[(b)] the semigroup $A$ is \nmt-continuous;
    \item[(c)] the generator $K$ of $A$ is in $\Ao$, i.e., $A_t =
               \e^{tK}$ for any $t\ge 0$.
   \end{enumerate}
 \item\label{semigroup-ii}
   The compression $R_{t}(a):=\skpo{u_{t}}{a\,u_{t}}$ for $a\in\Ao$ 
   and $t \ge 0$ defines a pointwise weakly* continuous Markovian 
   semigroup $R$ on $\Ao$. If $u$ is
   \nmot-continuous then $R$ is uniformly continuous.
 \end{enumerate}
\end{Prop}

Notice that the state $\psi$ of the continuous Bernoulli shift may not
be $R$-invariant.

\begin{proof}
  (\ref{semigroup-i}): From the cocycle identity follows $A_{s+ t} =
  \Eo  u_{s+t} = \Eo  (S_t u_s) u_t$. Since $S_tu_s$ and $u_t$ are
  $\Ao$-independent and $\Eo  S_t= \Eo $, we conclude $A_{s+t} = A_s
  A_t$ from equation \eqref{factor-2a}. The equivalence of (b) and (c)
  is evident from the theory of semigroups on Banach spaces (see
  e.g. \cite{Davi80a}). The equivalence of (a) and (b) 
  comes essentially from Cauchy-Schwarz's inequality \eqref{CS-ineq}:
  \begin{align}\label{cont-u-A}
        \nm{A_{t}-\1}^{2}
   & =  \nm{\skpo{\1}{u_{t}-\1}}^{2}
    \le \nmO{u_{t}-\1}^{2}
     =  \nm{\1-A_{t}^{*}-A_{t}^{}+\1}\,.
  \end{align}
  (\ref{semigroup-ii}): The semigroup property is a consequence of
  \ref{unit-cocycle}\spc(\ref{unit-cocycle-iii}) and \eqref{factor-lap}:
  \begin{align*}
        R_{s+t}(a) 
     =  \bskpo{(S_{t}u_{s})u_{t}}{a(S_{t}u_{s})u_{t}} 
     =  \bskpo{u_{t}}{\skpo{u_{s}}{a\,u_{s}}u_{t}} 
     =  R_{t}(R_{s}(a)).
  \end{align*}
  Using \eqref{CS-ineq} again we obtain the $\nmt$-continuity of $R$
  from \nmot-continuity of $u$:
  \begin{align}\label{cont-u-R}
        \nm{R_t(a) -a}
   &\le \nm{\skpo{u_{t}-\1}{a(u_{t}-\1)}}
       +\nm{\skpo{u_{t}-\1}{a\1}}
       +\nm{\skpo{a^{*}\1}{u_{t}-\1}}                         \notag \\
   &\le \nmo{u_{t}-\1}\nm{\skpo{u_{t}-\1}{a^{*}a(u_{t}-\1)}}^{1/2}
       +2\nm{a}\nmo{u_{t}-\1}                                 \notag \\
   &\le 4\nm{a}\nmo{u_{t}-\1}\,.
  \end{align}
  All other properties of the semigroup are clear by
  construction.\qedopt
\end{proof}

Throughout this paper we will consider only  $\nmt$-continuous
semigroups on $\Ao$.

\subsection{Additive cocycles of continuous GNS Bernoulli shifts}
\label{subsection:ac-cgnsbs}
We turn our attention to additive cocycles which are
non-commutative versions of L\'evy processes with values in an operator
algebra.

\begin{Defn}\label{add-cocycle}\label{def:add-cocycle}
  Let \CBSE be an $\Ao$-expected
  continuous GNS Bernoulli shift.  An \emph{additive cocycle} is a
  family $b:=(b_{t})_{t\ge0} \subset{\Lapneu}$ such that for any $s,t
  \ge 0$
  \begin{enumerate}
  \item\label{add-cocycle-i}
    $t \mapsto b_{t}$ is \nmot-continuous.
  \item\label{add-cocycle-ii}
    $b_{t}\in \LapXneu{[0,t]}$,                     \hfill(adaptedness)
  \item\label{add-cocycle-iii}
    $b_{t+s}= b_{t} + S_{t}b_{s}$.             \hfill(cocycle identity)
  \end{enumerate}
  The additive cocycle $b$ is \emph{centred} if in addition $\Eo
  b_{t}=0$ for all $t\ge0$. The
  operator $\skpo{b_t}{b_t}$ is called the \emph{variance} and $\Eo
  b_{t}$ the \emph{drift part} of $b$.
  
  $\Ca(\Lapneu,+)$ is the set of all additive cocycles.
  $\Cao(\Lapneu,+)$ is the set of all additive cocycles $b$
  satisfying the structure equation
  \[ 
        \skpo{b-\Eo  b}{b -\Eo b} + \Eo b + (\Eo b)^* 
     =  0.
  \]
\end{Defn}
\begin{Notation} \normalfont
Whenever it is convenient, we will use $\Eo b := (\Eo b_t)_{t\ge 0}$ 
and $\skpo{b}{c}:=(\skpo{b_t}{c_t})_{t\ge 0}$ for additive 
cocycles $b$ and $c$. We have already used this convention 
in above structure equation.    
\end{Notation}
Instead of $\nmot$-continuity in condition (\ref{add-cocycle-i})  one
could require the additive cocycle to be weakly* continuous  (or
measurable). But standard arguments show that all these conditions
are equivalent. Note also that $\nmot$-continuity is a redundant
requirement for a centred additive cocycle $b$. This continuity
property follows already  from the martingale property
\begin{align}\label{martingal-prop}
      E_{(-\infty,t]}b_{t+s}  
 & =  b_{t}\, \quad (s\ge 0).
\end{align} 
and the continuity of the past filtration (see Corollary
\ref{Corollary:cexp-pastfuture}). Thus one needs only  some continuity
or measurability condition to ensure $\Eo b_t= t\,\Eo b_1$ for any $t
\ge 0$ (compare the proof of the next proposition).

\begin{Prop}
  The variance  of a centred additive cocycle $b$ satisfies
  \begin{align} \abO{b_{t}}^{2}  & =  t \abO{b_{1}}^{2}\,.
  \label{bt-abs} \end{align} More generally, the covariance operator
  $\skpo{b_{t}}{\cdot\, c_{t}}$ of two centred additive cocycles
  $b$ and $c$ satisfies for any $a \in \Ao$ \begin{align}
  \skpo{b_{t}}{ac_{t}}  & =  t\skpo{b_{1}}{ac_{1}}\,.
  \label{inner-map} \end{align}
\end{Prop} 

\begin{proof} 
  Equation \eqref{bt-abs} is obtained from equation \eqref{inner-map}
  by putting $b=c$ and $a = \1$. For the proof of the latter one we
  define the linear map $\Ga_{t}:\Ao\ni a\mapsto\skpo{b_{t}}{ac_{t}}$.
  By the continuity of $b$ and $c$ in the stop topology mentioned
  above, $\Ga$ is pointwise weakly* continuous. From
  $\skpo{b_{t}}{aS_{t}c_{s}} = \skpo{b_{t}}{\1}\skpo{\1}{ac_{s}} = 0$
  we conclude $\Ga_{t+s}(a) = \skpo{b_{t}}{ac_{t}} +
  \skpo{S_{t}b_{s}}{S_{t}ac_{s}} = \Ga_{t}(a) + \Ga_{s}(a)$. Thus we
  obtain for any $\phi\in{\Ao}_{*}$ the functional equation
  $\phi(\Ga_{t+s}(a))=\phi(\Ga_{t}(a))+\phi(\Ga_{s}(a))$. Since $t
  \mapsto \phi(\Ga_{t}(a))$ is continuous, this functional equation has
  the unique solution 
  $\phi(\Ga_{t}(a))=t\phi(\Ga_{1}(a))$.  Equation \eqref{inner-map}
  follows now, since the predual of $\Ao$ is separating for
  $\Ao$.\qedopt
\end{proof}

\begin{Rem} \normalfont
  In the case of $q$-Gaussian white noises and the CAR-white noises,
  it is well-known that their additive cocycles are bounded
  in the operator norm and thus contained in the von Neumann algebra
  $\cA$ itself. Nevertheless, we refrained here to define explicitely 
  such additive cocycles for continuous Bernoulli shifts. 
  This additional feature can be exploited in applications, for
  example, if it is necessary or helpful on the computational side.
\end{Rem}

\subsection{The correspondence}       \label{subsection:correspondence}
Recall that $\Cmo(\cE,\boldsymbol{\cdot}\,)$ and $\Cao(\cE,+)$ are 
sets of $\nmot$-continuous unital resp.\ additive cocycles with
structure equation, as introduced in Definitions \ref{unit-cocycle}
and \ref{add-cocycle}. 
The abstract version of our main result is as follows. 

\begin{Thm}\label{theorem:ln-exp} 
  Let $\CBSE$ be an $\Ao$-expected continuous GNS Bernoulli shift. Then
  there exist two bijective mappings   
  \begin{align*}
        \Ln  & \colon \Cmo(\cE,\boldsymbol{\cdot}\,) 
    \to \Cao(\cE,+)                                                  \\
        \Exp & \colon \Cao(\cE,+) 
    \to \Cmo(\cE,\boldsymbol{\cdot}\,) 
  \end{align*}
  such that
  \begin{align*}
        \Ln \circ \Exp 
     =  \id 
     =  \Exp \circ \Ln.  
   \end{align*}
\end{Thm}

Notice that one has $\Cmo(\cE,\boldsymbol{\cdot}\,) = \Cm(\cE,
\boldsymbol{\cdot}\,)$ for $\dim\Ao<\infty$, since in this case every
semigroup $A$ is uniformly continuous and thus Proposition
\ref{semigroup}\spc(\ref{semigroup-i}) applies.

\begin{Notation}\normalfont
The cocycles $(\Exp(b)_t)_{t\ge 0}$ and $(\Ln(u)_t)_{t \ge 0}$
will also be written as $(\Exp(b_t))_{t\ge 0}$ resp.\ 
$(\Ln(u_t))_{t \ge 0}$. Since these two mappings will be constructed
'pointwise', this slight abuse of notation will vanish anyway. 
\end{Notation}

Above result can be viewed as an abstract corollary of Theorem
\ref{main-theorem} stated below. The latter one is formulated from a
(quantum) stochastic perspective and brings much more structure to the
surface. This will allow us to establish the mappings $\Ln$ and $\Exp$
in a constructive manner. The related proof is based on
non-commutative It\^o integration, non-commutative exponentials 
and logarithms. All these tools will be developed in Section
\ref{section:nc-ito-integration} and \ref{section:nc-ln-exp}, where we
will also finish the proof of Theorem \ref{theorem:ln-exp}.
  
\begin{Notation}\normalfont
 $\cZ(t)$ denotes an arbitrary net of
  partitions $Z$ of the interval $[0,t]$, which is partially ordered by
  inclusion, such that its grid tends to zero.
\end{Notation}

\begin{Thm}\label{main-theorem}
  Let \CBSE be an  $\Ao$-expected continuous GNS Bernoulli shift. The
  following are in a bijective correspondence:
  \begin{enumerate}
  \item\label{main-theorem-i}
    \nmot-continuous unital cocycles $u$ in ${\Lapneu}$; 
  \item\label{main-theorem-ii} pairs $(c,K$), where $c\subset\Lapneu$
    is a centred additive cocycles and $K\in\Ao$, satisfying the
    structure equation
    \[
          \abO{c_{t}}^{2}+t(K+K^{*})
       =  0.
    \]
  \end{enumerate}
  The unital cocycle $u$ is obtained from the pair $(c,K)$ as the
  solution of the non-commutative It\^o differential equation (IDE)
  \begin{align}\label{main-t-dgl}
        u_{t} 
   & =  \1+\int_{0}^{t}\!\d c_{s} u_{s} + \int_{0}^{t}\!\d t\,Ku_{s}\,.
  \end{align}
  Conversely, the centred additive cocycle $c$ is obtained from the
  unital cocycle $u$ as the non-commutative logarithm
  \begin{align}\label{main-t-constr}
        c_{t} 
   & =  \nmolim{Z\in\cZ(t)}\sum_{t_{i}\in Z}
        S_{t_{i}}(u_{t_{i+1}-t_{i}} - A_{t_{i+1}-t_{i}})\,,
  \end{align}
  and $K$ is the generator of the semigroup $A:=(\Eo u_{t})_{t\ge0}$.
  \end{Thm}
  
  Let us give some guidelines for the proof strategy. By doing so we
  will present, in particular, the concrete form of the mappings $\Ln$
  and $\Exp$.

\begin{proof}[Outlined proof of Theorem \ref{theorem:ln-exp} 
  and Theorem \ref{main-theorem}.]     
  In Section \ref{section:nc-ito-integration} we will develop a theory
  of non-commutative It\^o integration that is based solely on an
  $\Ao$-expected continuous GNS Bernoulli shift $\CBSE$ (see Definition
  \ref{HMBS}) and its additive cocycles $\Ca(\Lapneu,+)$ (see
  Definition \ref{add-cocycle}). Then
  Proposition \ref{proposition:nc-ito-integral} ensures that the
  non-commutative It\^o integral $\int_{0}^{t}\!\d c_{s} u_{s}$ is
  well-defined. With Theorem \ref{ex-eind} is established that the IDE
  \eqref{main-t-dgl} has a unique solution. Now, all terms are
  introduced and well-defined as they appear in the formulation of
  Theorem \ref{main-theorem}.
  
  In Section \ref{section:nc-ln-exp} we will develop the notion of
  non-commutative exponentials and logarithms for a given
  $\Ao$-expected continuous Bernoulli shift. We start in Subsection
  \ref{subsection:nc-exp} to investigate the IDE \eqref{main-t-dgl}
  for an arbitrary additive cocycle in $\Ca(\Lapneu,+)$. We ensure
  with Theorem \ref{sdgl-u} that an additive cocycle in
  $\Cao(\Lapneu,+)$ gives a $\nmot$-continuous unital cocycle as
  the solution of the IDE \eqref{main-t-dgl}. Now we are in the
  position to introduce non-commutative exponentials (see Definition
  \ref{definition:nc-exp}). From the uniqueness of the solution we
  conclude that the mapping $\Exp\colon \Cao(\Lapneu,+) \to
  \Cmo(\Lapneu,\boldsymbol{\cdot}\,)$ is well-defined by the family
  of IDEs
  \[
        \Exp(b_t) 
     =  \1 + \int_{0}^t \d b_s \Exp(b_s)\quad (b \in
        \Cao(\Lapneu,+), t\ge 0).
  \]
  Here we make use of the convention $\int\d b_s\cdot \:= \int
  \d c_s\cdot + \int\d s K\,\cdot$.
  
  We proceed in Subsection \ref{subsection:nc-ln} with the proof that
  \begin{align*}
        \Ln_0(u_t)
    :=  \nmolim{Z\in\cZ(t)}\sum_{t_{i}\in Z}
        S_{t_{i}}u_{t_{i+1}-t_{i}} - A_{t_{i+1}-t_{i}} 
        \qquad (t \ge 0)\,
  \end{align*}
  and
  \begin{align*}
        \Ln(u_t)
    :=  \nmolim{Z\in\cZ(t)}\sum_{t_{i}\in Z}
        S_{t_{i}}(u_{t_{i+1}-t_{i}} - \1) 
        \qquad (t \ge 0)\,
  \end{align*}
  are well-defined in $\Lapneu$ for every unital cocycle $u$ in
  $\Cmo(\Lapneu,\boldsymbol{\cdot}\,)$. In particular, we show that
  $\Ln(u_t) = \Ln_0(u_t) + K t$, where $K$ is the generator of the
  semigroup $A$. Moreover, we verify in Theorem \ref{unit-add-main}
  that $\Ln(u_t)$ is an additive cocycle in $\Cao(\Lapneu, +)$.
  Thus we have obtained the mapping $\Ln\colon
  \Cmo(\Lapneu,\boldsymbol{\cdot}\,) \to \Cao(\Lapneu, +)$, (see
  Definition \ref{definition:nc-ln}).  At this state all terms are
  well-defined, as they appear in the formulation of Theorem
  \ref{theorem:ln-exp}.
  
  We are left to prove that the correspondence in Theorem
  \ref{main-theorem} is bijective, respectively that the mappings
  $\Exp$ and $\Ln$ in Theorem \ref{theorem:ln-exp} are injective. But
  this is ensured by Proposition \ref{Exp-injectivity} and Proposition
  \ref{Ln-injectivity}. We complete the proof of Theorem
  \ref{theorem:ln-exp} and Theorem \ref{main-theorem} at the end of
  Subsection \ref{subsection:proof-correspondence}.\qedopt
\end{proof} 

If the shift $\CBS$ is trivial, i.e.\ if $\Ao=\cA$ and thus
$S=\id_{\cA}$, then the correspondence reduces to Stone's theorem on
uniformly continuous unitary groups. The following result emphasizes
the stochastic character of the additive and unital cocycles.

We remind that a function $\Rset \ni t \mapsto x_t \in {\Lapneu}$  is
weakly* differentiable if $\frac{\d}{\d t}\phi(x_{t})$ exists for any 
$\phi$ in the predual of ${\Lapneu}$ (see Theorem \ref{self-dual}).
   
\begin{Cor}\label{stone}
  If the unital cocycle $u$ or the additive cocycle $c$,
  respectively as stated in (\ref{main-theorem-i}) or
  (\ref{main-theorem-ii}) of Theorem \ref{main-theorem}, is weakly*
  differentiable, then $u$ is a semigroup of unitaries in $\Ao$ with
  generator $K = \i H$ for some selfadjoint operator $H \in \Ao$ and
  $c=0$.
\end{Cor}

Consequently, $\nmot$-continuous cocycles are
weakly* differentiable if and only if they are trivial.

\begin{proof}
  We will show that a weakly* differentiable unital cocycle $u$
  lies in $\Ao$, the fixed point algebra of the shift $S$. Since
  $u_t-\1 \in \LapXneu{(-\infty,t]} \cap \LapXneu{[0,\infty)}$
  for any $t>0$ and since the filtration $(\cA_{(-\infty,t]})_{t
  \in \Rset}$ is continuous we conclude that the weak* limit
  $u_{0}':=\lim_{t\to 0}\frac{1}{t}(u_t-\1)$ is in $\bigcap_{t>0}
  \Lapneu_{(-\infty,t]}\cap \Lapneu_{[0,\infty)} =
  \Lapneu_{(-\infty,0]}\cap \Lapneu_{[0,\infty)} = \Ao$.
  Furthermore, the cocycle identity implies $u^\prime_t = (S_t
  u^\prime_0) u_t = u^\prime_0 u_t$, where $u_{t}'$ is the weak*
  derivative of $u_{t}$. Consequently, $u_{t}=\e^{Kt}\in\Ao$ with
  $K:=u_{0}'$. Since $u_t = \Eo u_t$, the correspondence of Theorem
  \ref{main-theorem} implies $c_{t}=0$, hence $K+K^{*}=0$ and
  consequently $K = \i H$ for some selfadjoint operator $H \in \Ao$.
 
  Conversely, the centred additive cocycle $c$ is weakly*
  differentiable if and only if it is weakly* differentiable at $t=0$.
  Indeed, by the cocycle property, it is $c_t^\prime = S_t c_0^\prime$.
  It follows with similar arguments as for the unital cocycle that
  $c^\prime_0 \in \Ao$ and thus $c^\prime_t =c^\prime_0$. Hence $c_t =
  c^\prime_0 t =0$ for any $t\ge 0$, since the cocycle is centred.  Now
  $K+K^{*}=0$ entails again $K = \i H$ for some selfadjoint operator $H
  \in \Ao$ and \eqref{main-t-dgl} implies $u_{t}=\e^{\i Ht}$, which
   is weakly* differentiable.\qedopt
\end{proof}

The correspondence allows to identify the generator of the semigroup
$R$, introduced in Proposition \ref{semigroup}\spc(\ref{semigroup-ii}).

\begin{Cor}\label{CE-generator}
  The uniformly continuous semigroup of completely positive
  contractions $R$ on $\Ao$ has the Christensen-Evans generator
  \cite{ChEv79a}
  \begin{align*}
        \cL(a) 
     =  \Lambda(a) + K^*a + a K, \qquad a \in \Ao,
  \end{align*} 
  where $\Lambda := \skpo{c_1}{\cdot\, c_1}\,$.
\end{Cor}

\begin{proof}
  Since the convergence in \eqref{main-t-constr} is independent of the
  chosen net, we consider the sequence of equidistant partitions
  $Z_{n}:=\set{i\de_{n}}{\de_{n}:=2^{-n}t,i=0,\ldots,2^{n}}\in\cZ(t)$
  of $[0,t]$. Together with \eqref{inner-map} we calculate
  \begin{align*}
        \La(a) t
     =  \skpo{c_t}{a c_t}                                              
   & =  \lim_{n\to\infty}\sum_{i,j=0}^{2^n-1}
        \skpo{S_{i\de_n}(u_{\de_n}-A_{\de_n})}{%
              a S_{j\de_n}(u_{\de_n}-A_{\de_n})}                     \\
   & =  \lim_{n\to\infty}\sum_{i=0}^{2^n-1}
        \skpo{u_{\de_n}-A_{\de_n}}{a (u_{\de_n}-A_{\de_n})}          \\
   & =  \lim_{n\to\infty} \frac{1}{\de_n}
        (R_{\de_n}(a) - a + a- A_{\de_n}^*a A_{\de_n}^{}) t          \\
   & =  (\cL(a)-K^*a-aK) t \,.
  \end{align*}
  Thus, $\cL$ is given by $\cL(a) = \La(a) + K^{*}a + a K$. \qedopt
\end{proof}

The structure equation of the additive cocycle gives immediately
conditions, how additive cocycles should be composed for the
construction of unital cocycles. 

\begin{Lem}\label{lem:addition}
 Let $b,c \in \Cao(\Lapneu, +)$. Then it is
  \[
        b+c - \RE\skpo{b-\Eo b}{c-\Eo c}
    \in \Cao(\Lapneu, +). 
  \]
If the centred parts of $b$ and $c$ are, in addition, 
$\Ao$-independent or $\Ao$-orthogonal 
then $b+c \in \Cao(\Lapneu, +)$.  
\end{Lem}   
\begin{proof}
This is an elementary consequence of Definition 
\ref{def:add-cocycle}. \qedopt
\end{proof}
Notice that $\Ao$-independence of two centred additive cocycles
implies their $\Ao$-orthogonality. 

Lemma \ref{lem:addition} allows to generalize Corollary 
\ref{CE-generator} to a countable family of mutually 
$\Ao$-orthogonal additive centred cocycles $(c^{i})_{i\in\Nset}$  
with drifts $(K^{i})_{i\in\Nset}$, satisfying the structure equation.
If the sequences, defined by $\sum_{i=1}^n c^{i}_t + K^{i}t$, are 
stop-convergent for $n\to \infty$ for any $t>0$, then their limits define
an additive cocycle in $\Cao(\Lapneu, +)$. 
In this case, the Christensen-Evans generator of the semigroup $R$, 
associated to $\sum_{i\in\Nset} c^{i}_1 + K^{i}$, is given by
\begin{align*}
      \cL(a) 
   =  \bskpo{\sum_{i\in\Nset} c^{i}_1}{a \sum_{i\in\Nset} c^{i}_1}
     +\sum_{i\in\Nset} {K^{i}}^*a + aK^{i}\,.
\end{align*}
The more familiar form $  \sum_{i\in\Nset}\skpo{c^{i}_1}{a c^{i}_1} +
{K^{i}}^*a + aK^{i} $  of the generator is obtained if the considered
family of additive cocycles is mutually  $\Ao$-independent. Notice
that the sum over infinitely many terms is meant as limit in the weak*
topology on $\cA_0$.

\begin{Rem}\normalfont
  (i) The correspondence in Theorem \ref{main-theorem} provides a
  concrete Stinespring decomposition of uniformly continuous
  $CP_0$-semigroups on a von Neumann algebra. Further details are
  postponed to sequel publications.
 
  (ii) There are many other independent approaches to the dilation of
  uniformly continuous $CP_0$-semigroups on von Neumann algebras. For
  further information on this huge subject we refer the reader to
  \cite{Arve03a, BhSk00a, BBLS04a, GLW01a,GLSW03a, GoSi99a,
  LiWi00a,MuSo02a} and the references therein. In this context we
  want further mention the approach of \cite{Sauv86a,CiSa03a,CiSa03b}.
\end{Rem}

\subsection{Non-commutative white noises}     \label{subsection:CBS-WN}
Our main result, stated in Theorem \ref{main-theorem}, establishes a
bijective correspondence between unital and additive shift cocycles.
Actually, it contains already an (abstract) result for unitary shift
cocycles, as introduced in Definition  \ref{unitary-cocycle}. Recall
from Subsection  \ref{subsection:correspondence}  that
$\operatorname{Ln}$ maps a $\nmot$-continuous cocycle to an
additive cocycle.

\begin{Thm}\label{main-theorem-unitary}
  Let $\CBS$ be a $\Ao$-expected continuous Bernoulli shift. Then there
  exists a bijective correspondence between the set of
  $\nmot$-continuous unitary cocycles
  $\Cmo(\cA,\boldsymbol{\cdot}\,)$ and the set of additive
  cocycles $
  \Ln\big(\Cmo(\cA,\boldsymbol{\cdot}\,)\big)\subset
  \Cao(\Lapneu,+)$.
\end{Thm}

\begin{proof}
  Every unitary cocycle defines a unital cocycle. According to
  Theorem \ref{main-theorem}, the latter one corresponds uniquely to an
  additive cocycle. This ensures the existence of the claimed
  bijection. \qedopt
\end{proof}

What is the structure of the set $\Ln(\Cmo(\cA,\boldsymbol{\cdot}\,))$?
The development of the necessary material for a satisfying
answer goes beyond the limits of the present paper. Due to the
importance of this question, let us at least outline an answer for
$\Cset$-expected shifts with a tracial state $\psi$ (see also the
survey \cite{Koes03a}). Notice that it holds $\Eo = \psi(\cdot)\1$ in
this case.

Up to now we have not further specified the GNS Hilbert space of a
continuous Bernoulli shift. It turns out that the non-commutative
$L^2$-space provides all the further infrastructure which is sufficient
to reveal the structure of $\Ln(\Cmo(\cA,\boldsymbol{\cdot}\,))$. We
have available the whole scale of non-commutative $L^p$-spaces and in
particular the adjoint $b_t^*$ of an additive cocycle $b_{t}$. More
importantly, the sesquilinear quadratic variation
\begin{align*}
      \llbracket b,b \rrbracket_{t}
 &:=  \lim_{n\to \infty}
      \sum_{j=0}^{n-1} S_{jt/n} |b_{t/n}|^2,
\end{align*}
is well-defined as an $L^1$-norm limit \cite{Koes04a}. Now it turns
out that a $\nmot$-continuous unital cocycle is unitary if and
only if the corresponding additive cocycle $b$ satisfies
\begin{align*}
      \llbracket b,b \rrbracket_{t} + b_t^* + b_t^{} 
 & =  0. 
\end{align*}
Notice that the compression of this structure equation by the
conditional expectation $\Eo $ gives back the structure equation as it
appears for the correspondence between unital and additive
cocycles. This conceptual approach works also in the more general
setting of an $\Ao$-expected continuous Bernoulli shift, at least in
the case of a tracial state. Publications on this issue are in
preparation, including the non-tracial case.

We return to the discussion of the structure of continuous Bernoulli
shifts and present an operator algebraic notion of `white
noise'.

\begin{Defn}\label{def-white}
  An $\Ao$-expected continuous Bernoulli shift $\CBS$ is called an
  \emph{$\Ao$-expected (non-commutative) white noise} if for any $t >
  0$
  \begin{align}\label{white}
        \cA_{[0,t]}
     =  \bigvee\set{u_s}{u \in 
        \Cm(\cA,\boldsymbol{\cdot}\,), 0 \le s\le t}.
  \end{align}
  The $\Ao$-expected white noise is said to be \emph{generated by
  $\Cm(\cA,\boldsymbol{\cdot}\,)$}.
\end{Defn}

Notice that $\Ao \subseteq \bigvee\set{u_s}{u \in
\Cm(\cA,\boldsymbol{\cdot}\,), 0 \le s\le t}$ for any $t > 0$,
because each unitary in $\Ao$ defines canonically a trivial unitary
cocycle and a von Neumann algebra is generated by its unitaries. A
plausible explanation of `whiteness' is provided for the reader's
convenience at the end of this subsection.

It is easy to see from the cocycle equation that an  $\Ao$-expected
white noise is always locally minimal. This implies
upward continuity and in particular $\cA_I= \cA_{\overline{I}}$ for
its filtration   (see Subsection \ref{subsection:local-min-max}).

\begin{Prop}
  An $\Ao$-expected white noise with an enriched
  independence is locally minimal, locally maximal and has a continuous
  filtration.
\end{Prop}

\begin{proof}
  Enriched independence implies the local maximality of the filtration
  and thus downward continuity (see Subsection
  \ref{subsection:enriched-independence}).
\end{proof}

We expect that the enriched independence condition can be removed, due
to the continuity of the unitary cocycles, but we didn't yet look
closer at this problem.  It is evident that all continuity properties of
the filtration pass to the corresponding continuous GNS
Bernoulli shift $\CBSE$.

At this point it is worthwhile to remind our convention. We use the
attribute  `non-commutative' always in the sense of `not necessarily
commutative'. If an $\Ao$-expected non-commutative white noise does not
come from (operator-expected) probability theory, we will call it a 
\emph{quantum white noise}.

\begin{Rem}\normalfont 
  (i) If the von Neumann algebra $\cA$ of a $\Cset$-expected white noise
  is commutative then its measure theoretic version is identified as a
  `classical noise' in \cite{Tsir04a}. Note also that unitary
  cocycles play a generating role similar as it do units of
  Arveson's product systems. Thus white noises may also be called
  continuous Bernoulli shifts of `type~I'. The `classical or type~I'
  part is well-understood for continuous product systems of probability
  spaces or Hilbert spaces. Stressing the analogy, we hope to gain a
  better understanding of `continuous commuting square systems of
  operator algebras', starting in the (time-)homogeneous `type~I'
  setting. 
  
  (ii) There exist other notions of `quantum white noise' in the
  literature. Essentially, these approaches have in common that they
  start with generalized or quantum Brownian motions (or L\'evy
  processes). These processes are given explicitely in (deformed) Fock
  spaces and generate their `quantum white noises'. Out of our
  results arises the question whether every $\Cset$-expected 
  non-commutative white noise induces a deformed Fock space such that 
  it can be generated from quantum L\'evy processes on this Fock
  space. Notice also in this context that such a Fock space structure
  is anticipated by multiple non-commutative It\^o integrals which 
  can be formulated easily, starting from the results in Section 
  \ref{section:nc-ito-integration}.     
\end{Rem}

The following result states that one can always extract the `classical'
or `type~I part' of a shift.

\begin{Prop}\label{prop:wn-part}
  Let $\CBS$ be an $\Ao$-expected continuous Bernoulli shift $\CBS$.
  Then there exists a conditional expectation $E$ such that the
  compression of $\CBS$ by $E$ is an $\Ao$-expected non-commutative
  white noise generated by $\Cm(\cA,\boldsymbol{\cdot}\,)$.
\end{Prop}

Notice that a single unitary cocycle 
$u\in\Cm(\cA,\boldsymbol{\cdot}\,)$ (together with
all trivial unitary cocycles) may not generate a non-commutative
white noise. But this is guaranteed if 
$u \subset  \Cm(\cA,\boldsymbol{\cdot}\,) \cap \cA^\psi$.

\begin{proof} 
  Let $\cB_{[r,r+t]}:= \bigvee \set{S_r(u_s)}{u \in
  \Cm(\cA,\boldsymbol{\cdot}\,), r \in \Rset, 0\le s \le t}$ and
  define similarly $\cB_I$ for more general intervals $I \subseteq
  \Rset$. We conclude from the $\sigma^\psi$-invariance of $\cB_I$ that
  the conditional expectation $E_{[r,r+t]}$ from $(\cA,\psi)$ onto
  $\cB_{[r,r+t]}$ exists. We note that $\cB_\Rset$ is $S$-invariant and
  $S_t(\cB_I) = \cB_{I +t}$. The $\Ao$-independence of $\cB_I$ and
  $\cB_J$ for disjoint $I$ and $J$ follows immediately from the
  inclusions $\Ao\subseteq \cB_K \subseteq \cA_K$ for any $K \in \cI$.
  Thus $(\cB, \psi_{|\cB}, S_{|\cB}, (\cB_I)_{I\in \cI})$ is an
  $\Ao$-expected continuous Bernoulli shift and, by construction, an
  $\Ao$-expected non-commutative white noise. \qedopt
\end{proof}

\begin{Prop}
  Any $\Ao$-expected non-commutative white noise $\CBS$ is
  generated already by a finite or countable set of unitary cocycles.
\end{Prop}

\begin{proof}
  Since $\cA$ has a separable predual and is represented with respect
  to the faithful normal state $\psi$, there exists a stop-dense
  sequence $(x_n)_{n \in \Nset}$ \cite[Proposition 3.8.4]{Pedersen}.
  Kaplansky's density theorem (and the separability of the predual)
  ensures that each $x_n$ can be approximated by sequence $(x_{n,k})_{k
  \in \Nset}$ in the algebraic hull of the set $\set{S_t u_s}{u \in
  \Cm(\cA,\boldsymbol{\cdot}\,), t \in \Rset, s\ge 0}$. Clearly there
  are at most countable many cocycles involved to generate all elements
  $(x_{n,k})_{n,k \in \Nset}$. \qedopt
\end{proof}

An immediate consequence of Theorem \ref{main-theorem-unitary} is the
following result.

\begin{Cor}\label{cor-shift2white}
  If the set of centred additive cocycles is $\{0\}$, then the
  $\Ao$-expected continuous Bernoulli shift $\CBS$ restricts to a
  trivial non-commutative white noise.
\end{Cor}

\begin{proof} 
  It follows from Theorem \ref{main-theorem-unitary} that
  $\Cm(\cA,\boldsymbol{\cdot}\,)$ is a subset of $\Ao$. But this
  implies the triviality of the compressed shift stated in Proposition
  \ref{prop:wn-part}.  \qedopt
\end{proof}

From the cocycle identity follows that a continuous Bernoulli shift
without local minimality already fails to be a white noise. An example
of a continuous Bernoulli shifts without locally minimal filtration is
presented in Subsection \ref{subsection:local-min-comp}. We emphasize
that a locally minimal \emph{commutative} ${\Cset}$-expected
continuous Bernoulli shift may not be a white noise. The surprising
existence of such examples emerges from work of Tsirelson and Vershik
on the construction of intrinsically non-linear random fields
\cite{TsVe98a}. We summarize their result, phrased in our terminology,
as follows (see also Example \ref{black-noise}).

\begin{Thm}[Tsirelson-Vershik]
  There exist non-trivial locally minimal commutative 
  $\Cset$-expected continuous Bernoulli shift with $\{0\}$ 
  as the set of centred additive cocycles.
\end{Thm}

We close this subsection with a (widely known) explanation of
`whiteness' which also captures the non-commutative case. By Theorem
\ref{main-theorem-unitary}, a $\nmot$-continuous unitary cocycle $u$
determines uniquely a centred additive cocycle $b$. 
Let $\cK$ be the
Hilbert space given by the closed $\Cset$-linear span of $\set{S_t
b_s}{t \in \Rset, s\ge 0}$. The shift $S$ defines on $\cK$ a 
strongly continuous unitary
group, denoted by the same symbol.  Thus we obtain a Hilbert space
cocycle system $(\cK, S,b)$ which captures the linear theory of a 
white noise (up to multiplicity). Its theory is equivalent to that of 
the triple {$(L^2(\Rset), (\sigma_t)_{t \in \Rset}, (\chi_{[0,t]})_{t\ge 0})$}
(see e.g.\ \cite{Guic71a,Guic72a}).
Here denotes $\sigma_t$ the right
shift on $L^2(\Rset)$ and the family $(\chi_{[0,t]})_{t \ge 0}$
satisfies the cocycle equation $\chi_{[0,t+s]}= \chi_{[0,t]} +
\sigma_t \chi_{[0,s]}$. In the discussion of
`whiteness', it is instructive to think of `white noise' as the 
family of formal derivatives `$db_t/dt$', thus as the family of 
Dirac distributions
$\delta_t$ centred at $t$ in the equivalent picture. The latter one is
better expressed in the spectral representation of the shift. By
Fourier transformation one obtains the triple $(L^2(\Rset), (\e^{\i
t\lambda})_{t\in \Rset}, ((\i\lambda)^{-1} (\e^{\i t
\lambda}-1))_{t\ge 0})$. The derivatives
`$db_t/dt$' are now given by the family of functions $(\e^{\i
t\lambda})_{t\in \Rset}$ (not contained in $L^2(\Rset)$). The
attribute `white' features that these functions (`frequency modes')
evolve independently in time $t$ \emph{and} that they are equally
weighted by the Lebesgue measure. It is worthwhile to remind that
$\Cset$-expected white noises, each of them generated by a single
unitary cocycle, have \emph{all} the same linear theory as sketched
above. It is determined by the second order correlation functions of
the additive cocycle. The differences between $\Cset$-expected
non-commutative white noises appear, aside of multiplicities, by
looking at the correlation functions of higher orders.

\subsection{Examples for the correspondence}
\label{subsection:examples-correspondence}
We illustrate the explicit form of additive cocycles which satisfy
the structure equation of Theorem \ref{main-theorem} and thus lead to
unital  cocycles. Let us remark that all examples presented below
are $\Ao$-expected  non-commutative white noises in the sense of
Definition \ref{def-white}) (we will omit here most arguments about
this fact).
 
\begin{Exmp}[Gaussian white noise]\normalfont
  \label{Gaussian-correspondence}
  We continue the discussion of Example \ref{Gaussian-wn}. The
  $\Cset$-expected GNS Bernoulli shift is given explicitely by 
  $(L^2(\srS^\prime, \Sigma, \mu), \1, S, L^2(\srS^\prime,
      \Sigma_I, \mu_I))$. Here denotes $\mu_I$ the restriction of
  $\mu$ to $\Sigma_I$. The Brownian motion $B_t \in L^2(\srS^\prime,
  \Sigma, \mu)$ is the limit of $(X_{f_n})_{n\in \Nset}$ with $f_n \to
  \chi_{[0,t]}$ in the $L^2$-norm. All additives cocycle $b_t$ are of
  the form $\lambda B_t + K t$, where $\lambda, K \in \Cset$. It is
  elementary to check that $b_t$ satisfies the structure equation as
  stated in Theorem \ref{main-theorem}\spc(\ref{main-theorem-ii}) if
  and only if $\RE{K} = -\ab{\lambda}^2/2$ and $\IM{K}=h$ for some $h
  \in \Rset$. Thus we obtain $b_t = \lambda B_t -(\ab{\lambda}^2/2 -\i
  h)t$. By straightforward calculations the corresponding unital
  cocycle is given by
  \[
        \Exp{b_t}
     =  \exp(\lambda B_t -(\lambda \RE(\lambda) -\i h)t).
  \]
\end{Exmp}

\begin{Exmp}[Poisson white noise]\normalfont
  \label{Poisson-correspondence}
  Consider the Poisson process $N$ with intensity $\lambda >0$ of
  Example \ref{Poisson-wn}. Then $c_t:= \ve(N_t-\la t)$ defines a
  centred additive cocycle. From $\psi_\mu(|c_t|^2) = |\ve|^2 \la
  t$ we conclude that the additive cocycle $b_t := c_t + K t$ satisfies
  the structure equation if and only if $K = - |\ve|^2 \la/2 +\i h$ for
  some $h \in \Rset$. The explicit form of the corresponding unital
  cocycle is found after some calculations to be of the form
  \[
        \Exp{b_t}
     =  (1+\ve)^{N_t} \exp{(-(\tfrac{1}{2}|\ve|^2\la+\ve\la-\i h)t)}.
  \]
\end{Exmp}

\begin{Exmp}[$\operatorname{CCR}$ white noises]\normalfont
  \label{CCR-correspondence}
  We continue Examples \ref{CCR-independence} and \ref{CCR-white-noise}.
  To find the form of additive cocycles, we pass to a concrete GNS
  representation of the C*-algebra $\operatorname{CCR}(L^2(\Rset),\IM
  \skp{\cdot}{\cdot})$ with the state $\psi_\lambda(W(f))=
  \exp{(-(2\lambda+1)/4\,\nm{f}^2)}$ ($\lambda >0$). This
  representation of Araki-Woods type is given on the tensor product of
  two symmetric Fock spaces $\cF_+(L^2(\Rset)) \otimes
  \cF_+(L^2(\Rset))$ with cyclic separating vector $\Omega \otimes
  \Omega$ such that
  \[
        \psi_\lambda(W(f))
     =  \skp{\Omega \otimes \Omega}%
        {W_\cF(\sqrt{\la +1}f)\Omega \otimes W_\cF(\sqrt{\la}Jf)\Omega}
  \]
  Here denotes $W_\cF(f)$ the Weyl operator on the symmetric Fock
  space and $J$ the complex conjugation on $L^2(\Rset)$. One finds
  from the represented Weyl operators, via Stone's Theorem, the
  annihilation operator
  \[
        a_\lambda(f)
    :=  \sqrt{\la+1}\, a(f) \otimes \1 + \sqrt{\la}\, \1 \otimes a^*(Jf) 
  \]
  and the annihilation operator $a_\lambda(f)^*$ as its adjoint, in
  terms of the usual annihilation operator $a(f)$ on
  $\cF_+(L^2(\Rset))$. This gives immediately the general form of a
  centred additive cocycle as a linear combination of
  $a_\lambda(\chi_{[0,t]})$ and $a_\lambda^*(\chi_{[0,t]})$:
  \[
        c_t
    =   a_1 \sqrt{\lambda}\, \Omega \otimes \chi_{[0,t]} 
       +a_2 \sqrt{\lambda+1}\,\chi_{[0,t]}\otimes \Omega
       \qquad (a_1, a_2 \in \Cset). 
  \]
  One verifies  easily from this form that an additive cocycle $b_t
  = c_t + K t$ satisfies the structure equation whenever $K =
  -1/2(\lambda |a_1|^2 + (\lambda+1)|a_2|^2) + \i h$ for some $h \in
  \Rset$. Obviously, the fixed constants $\sqrt{\lambda}$ and
  $\sqrt{\lambda+1}$ are superfluous. Their effect can be compensated
  by rescaling the coefficients, as long as one calculates only
  time-ordered higher moments factorizing into second moments. 
\end{Exmp}

The examples \ref{Gaussian-correspondence} to
\ref{CCR-correspondence} generalize straightforward to the case of 
$\Ao$-expected white noises with infinite multiplicity
\[
      (\Ao \otimes (\bigotimes_{n=1}^\Nset \cC),
      \psi \otimes (\bigotimes_{n=1}^\Nset \psi),
      \id \otimes (\bigotimes_{n=1}^\Nset S), 
      (\Ao \otimes (\bigotimes_{n=1}^\Nset \cC_I))_{I\in \cI}). 
\]
Here denotes $(\cC, \psi , S, (\cC_I)_{I\in \cI})$ a
$\Cset$-expected white noise as considered in the last three examples.
The general form of an additive cocycle is now
\[
       b_t
   =  \sum_n a_n \otimes c_t^{(n)} + (K \otimes \id)t,
\]
where $K \in \Ao$ and $c^{(n)}_t$ denotes the canonical embedding of
the $\Cset$-expected centred additive cocycle with variance
$\skp{c_1}{c_1}=1$ into the infinite tensor product at the $n$th
position. Moreover, it is required for the sequence $(a_n)_{n\in\Nset}
\subset \Ao$ that $\sum_{i=1}^n |a_i|^2$ is stop convergent for $n \to
\infty$. The structure equation is satisfied by $b$ if and only if
$\RE K = -1/2\sum_n |a_n|^2$ and $\IM K= h$ for some selfadjoint
operator $h \in \Ao$. One obtains, as usual, the Markovian semigroup
\[
      \Ao
  \ni x \mapsto R_t(x)
  :=  \skpo{\Exp(b_t)}{x\Exp(b_t)},
\]
where we identify $\Ao$ and  $\Ao \otimes \1$. As it is
well-known, this semigroup has a generator of Lindblad form:
\[
     \cL(x)
   = \sum_n (a_n^*x a_n^{}-\tfrac{1}{2}\{a_n^*a_n^{},x\})+\i[h,x] 
\]
(here denote $\{a,b \}$ and $[a,b]$ the anti-commutator resp.\ the
commutator).

\begin{Exmp}[$\operatorname{CAR}$ white noises]\normalfont
  \label{CAR-correspondence}
  We continue the discussion of Examples \ref{CAR-independence} and
  \ref{CAR-wn}. Similar as done for the $\operatorname{CCR}$-algebra,
  we pass for the C*-algebra $\operatorname{CAR}(L^2(\Rset))$ with the
  quasi-free state $\psi_\lambda(a^*(f)a(g))= \lambda \skp{g}{f}$ ($0
  <\lambda <1$) to an Araki-Woods representation on the tensor product
  of two antisymmetric Fock spaces $\cF_{-}(L^2(\Rset)) \otimes
  \cF_{-}(L^2(\Rset))$ and find that an additive cocycle is always
  of the form
  \[
        b_t
     =  a_1 \Omega \otimes \chi_{[0,t]} 
       +a_2 \chi_{[0,t]} \otimes \Omega \qquad (a_1, a_2 \in \Cset).
  \]
  In the case of $\operatorname{CAR}(K_0 \oplus L^2(\Rset, \cK_1))$,
  where one obtains a $\cB_0$-expected white noise, the general form of
  an additive cocycle is now 
  \[
        b_t
     =  \sum_{n} x_n \Omega \otimes (\chi_{[0,t]}\otimes e_n) 
       +y_n (\chi_{[0,t]}\otimes e_n) \otimes \Omega \qquad 
        (x_n, y_n \in \cB_0).
  \]
  Here, $\{e_n\}_{n}$ is an orthonormal basis of $\cK_1$. Moreover, we
  identified $L^2(\Rset,\cK_1)$ and $L^2(\Rset)\otimes \cK_1$.  If
  $\cK_1$ is infinite dimensional, one needs also that $\sum_i^n
  {x_i^*x_i^{}} + {y_i^*y_i^{}}$ is stop convergent for $n \to \infty$.
  The computation of the Christensen-Evans generator is
  straightforward. We leave these details to the reader.
\end{Exmp}

\begin{Exmp}[$q$-Gaussian white noises]\normalfont
  \label{q-correspondence} 
  We continue the discussion of Examples \ref{free-independence},
  \ref{q-independence} and \ref{q-white-noise}. Let us immediately
  consider the case
  \[
        (\cF_0(\cK_0 \oplus L^2(\Rset)), 
        \tau, S, (\cF_0(\cK_0 \oplus L^2(I))_{I\in \cI}). 
  \]
  The general form of an $\cF_0(\cK_0)$-expected additive cocycle is
  now
  \[ 
        b_t
     =  \sum_{i=1}^n 
        a_i \Phi(0 \oplus \chi_{[0,t]}) \tilde{a}_i + K t,  
  \]
  where $a_i, \tilde{a}_i, K \in \cF_q(\cK_0 \oplus 0)$ ($1 \le i \le
  n$). The structure equation is satisfied if $K =
  \frac12\sum_{i,j=1}^n \tilde{a}_i^* \Gamma_q(q)(a_i^*a_j^{})
  \tilde{a}_j^{} + \i h$ for some selfadjoint operator $h \in \cB_0$.
  Here denotes $\Gamma_q(q)$ the second quantization of the
  multiplication operator $M_q(g) := qg$ with $g \in \cK_0 \oplus
  L^2_\Rset(\Rset)$. Moreover we used the identity
  \[
        \Eo (\tilde{a}^* \Phi(0 \oplus \chi_{[0,t]}) 
        a^*a \Phi(0 \oplus \chi_{[0,t]})\tilde{a})
     =  t \tilde{a}^* \Gamma_q(q)(a^*a)\tilde{a},
  \]
  where $a, \tilde{a} \in \cF_q(\cK_0 \oplus 0)$ (see \cite{DoMa03a}).
  Notice that $\Gamma_0(0)(x) = \tau(x)$ in free probability. The
  generator of the semigroups has now the Christensen-Evans form
  \[
        \cL(x)
     =  \sum_{i,j=1}^n
        \big(\tilde{a}_i^* \Gamma_q(q)(a_i^*x a_j^{}) \tilde{a}_j^{} 
       -\tfrac12\,\{\tilde{a}_i^* \Gamma_q(q)(a_i^*a_j^{}) 
        \tilde{a}_j^{}, x \}\big) + \i[h,x]. 
  \]
\end{Exmp}

\begin{Rem}\normalfont
  This list of examples can be enlarged by examples coming from
  generalized Brownian motions, as soon as white noise functors are
  available and the vacuum vector of the deformed Fock space is
  separating for the von Neumann algebra, generated by these
  generalized Brownian motions.
\end{Rem}

\section{Non-commutative It{\^o} integration}
\label{section:nc-ito-integration}

This section is devoted to the development of operator-valued
non-commutative  It\^o integration, as it is needed for the
correspondence stated in Theorem \ref{main-theorem}. In the case of a
$\Cset$-expected Bernoulli shift with a commutative von Neumann algebra,
our approach reduces to an $L^2$-theory of stochastic It\^o
integration, as it is known for L\'evy processes, in particular
Brownian motion. Since we work with $\Ao$-expected shifts, this
approach covers It\^o integration for 
operator-valued L\'evy processes with uniformly bounded covariance
operators. Notice also that our filtrations may not necessarily be
`generated by the non-commutative L\'evy processes', we will make
use only of the additive cocycle's adaptedness. This adaptedness implies
already that a centred additive cocycle is a non-commutative
martingale with respect to the filtration of the given continuous
Bernoulli shift.

The GNS representation $\CBSE$ of a continuous Bernoulli shift \CBS
and its additive cocycles $\Ca(\Lapneu,+)$  provide in our
approach all the infrastructure which is needed to develop  this
theory of non-commutative It\^o integration. We will restrict our
presentation to results, as far as they are  necessary for the proof
of Theorem \ref{main-theorem}. We emphasize that our approach  to
It\^o integration is purely operator algebraic. A priori, it is not
based on (deformed) Fock space structures. Nevertheless, it can be
utilized for quantum stochastic integration on (deformed) Fock spaces,
reserved to the condition that a continuous Bernoulli shift is present
in the bounded operators of this Fock space.
 
Throughout this section, we work in the presence of a fixed
$\Ao$-expected continuous GNS Bernoulli shift $\CBSE$. Moreover, we
assume that this shift has at least one non-zero centred additive
cocycle. (We remind that Tsirelson's black noises  provide
examples  which have only trivial additive cocycles.)

\subsection{Non-commutative It{\^o} integrals for simple adapted
processes}                                \label{subsection:ito-simple}
In the following, a centred additive cocycle $\cac$ will serve as
the non-commutative generalization of a stochastic process with
stationary independent increments, like Brownian motion or Poisson
process.

We start with the usual notion of adapted processes, as it is known
from stochastic It\^o integration in probability theory. A
\emph{process} is a family $x =(x_t)_{t\ge 0}\subset{\Lapneu}$. It is
called \emph{(locally) adapted} if $x_{t}\in \LapXneu{(-\infty,t]}$
(resp.~$\LapXneu{[0,t]}$) for any $t\ge 0$. A \emph{simple adapted
process} $x =(x_t)_{t\ge 0}\subset {\Lapneu}$ is given by
\begin{align*}
      x 
  :=  \sum_{i\ge0}x_{i}\,\chi_{[s_{i},s_{i+1})}\,, 
      \qquad 
      x_{i}\in \LapXneu{(-\infty,s_{i}]}\,. 
\end{align*}
Here $Z:=\set{s_{i}}{s_{i}<s_{i+1},\, i\in\Nset_{0}}$ defines a
partition of $\Rset^{+}$. For a simple adapted process $x$ (putting
$s_m:=t_0$ and $s_n:=t$ for suitable $0\le m < n$), we introduce the
((left) non-commutative) It\^o integral
\begin{align}\label{simple-integral}
      \int_{t_0}^{t}\! \d \cac_{s} x_{s}
 &:=  \sum_{i=m}^{n-1} (\cac_{s_{i+1}}-\cac_{s_{i}})\,x_{i}\,.
\end{align}
It is sufficient to define the integrals of simple adapted processes
  just for intervals $[t_{0},t]$ with boundaries $t_{0},t\in Z$ (using
  a sub-partition of $Z$ if necessary). Indeed, this expression is
  well-defined, since $\cac_{s_{i+1}}-\cac_{s_{i}} =
  S_{s_{i}}\cac_{s_{i+1}-s_{i}}$ and $x_{i}$ are $\Ao$-independent,
  and thus, by Proposition \ref{new-product}, their product makes
  sense. The following operator identity in $\Ao$, the so-called
  \emph{(non-commutative) It\^o identity}, will be crucial for the
  extension of the integral to a larger class of adapted
  processes.

\begin{Lem}\label{simple-isometry} 
  If $x$ and $y$ are two simple adapted processes, then
  \begin{align*}
        \Bskpo{
        \int_{t_0}^{t}\!\d \cac_{s} x_{s}}%
        {\int_{t_0}^{t}\!\d \cac_{s} y_{s}}
   & =  \int_{t_0}^{t} \skpo{x_{s}}{\La(\1)\, y_{s}}\d s\,,
  \end{align*}
  where $\La=\skpo{\cac_{1}}{\cdot\, \cac_{1}}$ is the uniformly
  bounded covariance operator of $b$.
\end{Lem}

\begin{proof}
  By refinement, we may assume that $x$ and $y$ are simple adapted
  processes with respect to the same partition. From equation
  \eqref{factor-lap} we conclude
  \begin{align*}
        \skpo{(\cac_{s_{i+1}}-\cac_{s_{i}})\, x_{i}}{%
              (\cac_{s_{j+1}}-\cac_{s_{j}})\, y_{j}}
     =  0 
        \qquad \text{\ for\ } i\ne j\,. 
  \end{align*}                             
  Thus,
  \begin{alignat*}{2}
           \lefteqn{\Bskpo{\int_{t_0}^{t} 
           \!\d \cac_{s} x_{s}}{\int_{t_{0}}^{t} \!\d \cac_{s} y_{s}}
   \hspace*{0.5ex} =
           \sum_{i=m}^{n-1}
           \skpo{(\cac_{s_{i+1}}-\cac_{s_{i}})\,x_{i}}{%
           (\cac_{s_{i+1}}-\cac_{s_{i}})\,y_{i}}}\hspace{1.5cm}      \\
   &\hspace*{1.2ex}
    \stackrel{\makebox[0pt]{\scriptsize\eqref{factor-lap}}}{=}
    \hspace*{1.2ex} & &
           \sum_{i=m}^{n -1}
           \skpo{x_{i}}{\abO{\cac_{s_{i+1}-s_{i}}}^{2}\,y_{i}}
   \hspace*{1.2ex}
    \stackrel{\makebox[0pt]{\scriptsize\eqref{bt-abs}}}{=}
    \hspace*{1.2ex} 
           \sum_{i=m}^{n-1}
           \skpo{x_{i}}{\La(\1)\,y_{i}}(s_{i+1}-s_{i})               \\
   &\hspace*{1.2ex} =
    \hspace*{1.2ex} & &
           \int_{t_0}^{t}\skpo{x_{s}}{\La(\1)\,y_{s}}\d s\,. \qedopteqn
 \end{alignat*}
\end{proof}

\begin{Rem}\normalfont 
  In this paper we make  use only of left It\^o
  integrals. Nevertheless, using the same techniques, the right
  non-commutative It\^o integral $\int_{t_0}^{t}x_{s}\d \cac_{s}$ is
  introduced. It goes along with an It\^o identity of the form
  $\bskpo{\int_{t_0}^{t} x_{s}\d \cac_{s}}{\int_{t_0}^{t} y_{s}\d
    \cac_{s}} = \int_{t_0}^{t}\La(\skpo{x_{s}}{y_{s}})\d s$. We
  emphasize that here, in contrast to stochastic It\^o integration or
  Hudson-Parthasarathy quantum stochastic integration, the left and
  right  integral differ in most cases, since the
  past/future structure may not commute, even in the case $\Ao\simeq
  \Cset$. A typical example is  It\^o
  integration in the case of free Brownian motion (see also
  \cite{BiSp98a}).
\end{Rem}

\subsection{An extension of the non-commutative It{\^o} integral}
\label{subsection:ito-extension}
For the purposes of this paper, it will be sufficient to extend the
 It\^o integral, introduced in Subsection
\ref{subsection:ito-simple}, to the vector space $\cV$ of piecewise
stop-continuous adapted processes
\[
      \Rset^{+}_{0} \ni t 
  \mapsto 
      x_t \in \Lapneu
\]
(which are locally \nmot-bounded by the uniform boundedness
principle). Notice that $\cV$ includes locally bounded
\nmot-continuous adapted processes, in particular centred additive
cocycles.

The key for this extension is provided by the
 It\^o identity. We introduce on $\cV$ the
family of seminorms
\begin{align*}
      \d_{\xi_{0},I}(x)
  :=  \Bigl[\int_{I}\bnm{\abo{x_{s}}\xi_{0}}^{2}\d s\Bigr]^{1/2},
\end{align*}
where $\xi_{0}\in\Ho$ and $I$ ranges over all compact intervals in
$\Rset^{+}$. We have
\[
      \d_{\xi_{0},I}(x)
  \le \nm{\xi_{0}}\sup_{s\ge0}\nmo{x_{s}}\,\ab{I}.
\]
Notice that for a simple adapted process $x \in \cV$, the estimate
\begin{align*}
      \Bnm{\Babo{\int_{I}\d \cac_{s}^{}x_{s}}\xi_{0}}^2
 & =  \Bskp{\xi_{0}}{\int_{I} 
      \skpo{x_s}{\La(\1)\,x_s}\d s\, \xi_{0} } 
   =  \int_{I} 
      \bskp{\xi_{0}}{\skpo{x_s}{\La(\1)\,x_s} \xi_{0} }\d s          \\
 &\le \nm{\Lambda(\1)} \d_{\xi_{0},I}(x)^2.
\end{align*}
is valid. Here, we used the It\^o identity
and $\La(\1) \le \nm{\La(\1)}\1$.

Next, we show that processes in $\cV$ can be approximated by simple
processes in $\cV$. By a simple reduction argument, we may assume that
the process $x\in \cV$ is stop-continuous on $I$. We define the simple
adapted processes $x^{Z} := \sum_{i\in\Zset^{+}}
x_{s_{i}}\,\chi_{[s_{i},s_{i+1})} \in \cV$, where
$Z:=\set{s_{i}}{s_{i}<s_{i+1},\,i\in\Nset_{0}}$ denotes a partition of
$\Rset^{+}$. The continuity of $x$ and
\begin{align*}
      \d_{\xi_{0},I}(x-x^{Z})^{2}
 &\le \sum_{s_{i}\in I} 
      \int_{s_{i}}^{s_{i+1}}
      \bnm{\abo{x_{r}^{}-x_{s_{i}}^{}}\xi_{0}}^{2}\d r
     +\int_{s_{i-1}}^{s_{i}}
      \bnm{\abo{x_{r}^{}-x_{s_{i-1}}^{}}\xi_{0}}^{2}\d r\,,
\end{align*} 
imply $\d_{\xi_{0},I}(x-x^{Z}) \rightarrow 0$, whenever the grid of the
partition $Z$ tends to $0$. From this and from 
\begin{align}
      \Bnmo{\int_{I}\!\d \cac_{s}^{}x_{s}^{Z}}
 & =  \sup_{\nm{\xi_{0}}\le1}
      \Bnm{\Babo{\int_{I}\!\d \cac_{s}^{}x_{s}^{Z}}\xi_{0}}
  \le \nm{\La(\1)}^{1/2}\sup_{s\ge0}\nmo{x_{s}^{Z}}\ab{I}     \notag\\
 &\le \nm{\La(\1)}^{1/2}
      \sup_{s\ge0}\nmo{x_{s}^{}}\ab{I},               \label{mart-cont}
\end{align}
and furthermore 
\begin{align}      
      \Bnm{\Babo{\int_{I}\!\d \cac_{s}^{}x_{s}^{Z} 
     -\int_{I}\!\d \cac_{s}^{}x_{s}^{Z'}}\xi_{0}}
 &\le \nm{\Lambda(\1)}^{1/2} 
      \big(\d_{\xi_{0},I}(x^{Z}-x) 
     +\d_{\xi_{0},I}(x-x^{Z'})\big)\,,                         \notag
\end{align}
we conclude that $\big(\int_{I}\!\d \cac_{s}^{}x_{s}^{Z} \big)_{Z}$ is
a bounded Cauchy net in the stop topology on $\Lapneu$. Its limit in
${\Lapneu}$ is denoted by $\int_{0}^{t}\!\d \cac_{s}x_{s}$. Finally,
one verifies by routine arguments that the definition of the integral
is independent from the chosen net of partitions. Moreover, it is
elementary to see from the definition that $\int_{0}^{t}\!\d
\cac_{s}x_{s}$ carries all the usual properties, as $\int_{t_0}^{t}\!\d
\cac_{s}x_{s} = \int_{t_0}^{t_1}\!\d \cac_{s}x_{s} + \int_{t_1}^{t}\!\d
\cac_{s}x_{s} $ ($ t_0 \le t_1 \le t$), and as linearity with respect
to vector space structure of $\cV$. We summarize the above discussion
as follows:

\begin{Prop}\label{proposition:nc-ito-integral}
  The non-commutative It\^o integral extends from simple adapted
  processes in $\Lapneu$ to the vector space $\cV$ of piecewise
  stop-continuous adapted processes. Moreover, the
   It\^o identity
  \begin{align}
        \Bskpo{\int_{t_0}^t \!\d \cac_s x_s}{\int_{t_0}^t 
        \!\d \cac_s x_s} 
     =  \int_{t_0}^t
        \bskpo{x_s}{\Lambda(\1)\, x_s}\d s         \label{ito-identity}
  \end{align}
  is valid for any $x \in \cV$.
\end{Prop}

We remark that the integral on the right-hand side of the
 It\^o identity is a weak* integral in
$\Ao$.

Notice also that $X_t :=\int_{t_0}^{t}\! \d \cac_{r} x_{r}$ is a
\nmot-continuous non-commutative martingale with respect to the
filtration $(\cA_{(-\infty, s]})_{s\ge 0}$, i.e., $E_{(-\infty, s]}X_t
= X_s$ for any $t_0 \le s \le t$. This is seen easily for a simple
process $x$:
\begin{align*}
      E_{(-\infty,s]}\int_{s}^{t} \!\d \cac_{r} x_{r}
 & =  \int_{s}^{t} E_{(-\infty,s]}\circ E_{(-\infty,r]}
      \d \cac_{r} x_{r}                                              \\
 & =  \int_{s}^{t} E_{(-\infty,s]}E_{(-\infty,r]}(\!\d \cac_{r}) x_{r}
   =  0\,.
\end{align*}
This equation extends from simple processes to processes $x \in \cV$ by
approximation and thus proves that $(X_t)_{t_0 \le t}$ is a martingale.
The \nmot-continuity  follows directly from inequality
\eqref{mart-cont}, which, by \eqref{ito-identity}, extends to the
piecewise stop-continuous adapted processes. Moreover, $X_t$ is locally
adapted if $x$ is locally adapted.

Since a (non)-centred additive cocycle $\ncac$  decomposes uniquely
into its centred part $\cac:= \ncac-\Eo \ncac$ and its drift
$\Eo \ncac_{t}=t\Eo \ncac_{1}$, we let for any $x \in \cV$
\begin{align*} 
      \int_{t_0}^t \!\d \ncac_s\, x_s 
  :=  \int_{t_0}^t \!\d \cac_s  x_s 
     +\int_{t_0}^t (\Eo \ncac_1)\, x_s\d s\,. 
\end{align*} 
By routine calculations it is shown that
\begin{align*}  
      \BnmO{\int_{t_0}^t \!\d \ncac_s\, x_s}^2 
  \le 2(1+ t-t_0)\nmO{\ncac_1}^2 \int_{t_0}^t \nmO{x_s}^2\d s\,. 
\end{align*}

\begin{Rem}\normalfont 
  The following argument ensures that the vector space $\cV$ contains,
  roughly speaking, many processes. Let $P:=(E_{(-\infty,t]})_{t\ge0}$
  be the projection from the vector space of processes onto the vector
  space of adapted processes in $\Lapneu$. For a stop-continuous
  process $x$, the map $s\mapsto E_{(-\infty,s]}x_{s}$ is also
  stop-continuous, since $s\mapsto E_{(-\infty,s]}$ is pointwise
  stop-continuous. Consequently, a stop-continuous process is mapped
  by the projection $P$ to a stop-continuous adapted process. Notice
  that these arguments do not apply to locally adapted processes and
  the family $(E_{[0,s]})_{s\ge0}$ if $s \mapsto E_{[0,s]}$ fails to
  be pointwise stop-continuous.
  
  A further extension of the non-commutative It\^o integral to the
  class of $L^{2}$-integrable adapted processes is possible, but this
  requires tremendously more technical efforts. These processes are
  imposed to the condition of integrability of
  $t\mapsto\nm{x_{t}\xi}^{2}$ on any compact interval, aside of
  measurability conditions for any $\xi\in\Ho$. The relevant tools for
  the development of such a theory are well-known and can be found,
  e.g., in \cite{Take03a}, Chapter IV. We will not elaborate further
  this direction, because the integration class of stop-continuous
  processes will be sufficient for the purpose of this paper.
\end{Rem}

\subsection{Non-commutative It\^o differential equations}
\label{subsection:IDE}
We will close our digression on non-commutative It\^o integration with
an existence and uniqueness theorem for solutions of 
(non-commutative) It\^o differential equations (IDEs).

Let $\cac$ be a centred additive cocycle in $\Ca(\cE,+)$. We say
that the process $x \in \cV$ has the differential $\d x_{t} =
\al_{t}\d t + \d \cac_{t} \be_{t}$, if it has the form $x_{t} = x_{0}
+ \int_{0}^{t} \al_{s}\d s +\int_{0}^{t} \!\d \cac_{s} \be_{s}$ for
any $t\ge t_0$, where \al and \be are processes in $\cV$. We call a
function $\al:\,\Rset^+\times{\Lapneu}\to{\Lapneu}$ \emph{adapted}, if
$\al(t,y)\in \LapXneu{(-\infty,t]}$ for any $t\ge0$ and $y\in
\LapXneu{(-\infty,t]}$. Moreover, we say that the function $\al$ is
\emph{locally adapted}, if the previous statement is satisfied with
respect to $\LapXneu{[0,t]}$, instead of
$\LapXneu{(-\infty,t]}$. Finally, we say that such a function is
stop-continuous, if $t\mapsto \al(t,y)$ is stop-continuous for any
$y\in{\Lapneu}$ and $y\mapsto \al(t,y)$ is stop-stop-continuous for
any $t\ge0$.
 
If $\al$ and $\be$ are adapted, a process $x \in \cV$ is called a
solution of the non-commutative It\^o differential equation (IDE)
\begin{align}\label{sdgl}
      \d x_{t}
 & =  \al(t,x_{t})\d t + \d \cac_{t}\be(t,x_{t})\,,
      \qquad
      x_{t_{0}}\in \LapXneu{(-\infty,t_{0}]}\,,
\end{align}
if $x$ solves the integral equation
\begin{align}\label{sdgl-int}
      x_{t}
 & =  x_{0} 
     +\int_{t_{0}}^{t}\al(t,x_{s})\d s 
     +\int_{t_{0}}^{t}\!\d \cac_{s}\be(s,x_{s})\,.
\end{align}
In particular, the following result insures that the IDE, as stated
in  Theorem \ref{main-theorem}, has a unique solution.

\begin{Thm}\label{ex-eind}
  Let the $\Ao$-expected continuous GNS Bernoulli shift $\CBSE$ be
  given. Let $\cac \in \Ca(\cE,+)$ be a centred additive cocycle
  and $\al$, $\be$ adapted stop-continuous functions from
  $\Rset^+\times{\Lapneu}$ to ${\Lapneu}$. Furthermore, assume for
  $\alpha$ and $\beta$ that for any compact interval $[t_0,t_1]$ with
  $t_0\ge 0$ and $\xi \in \cH_0$, there exists a constant $C_{\xi}\ge
  0$ such that for any $x,y\in{\Lapneu}$ and $t,s\in [t_0,t_1]$
  \begin{align}
        \label{ex-eind-a}
        \d_{\xi}(\al(t,x)-\al(s,x))
   &\le C_{\xi}\,(\nm{\xi}^2+\d_{\xi}(x))\,
        \ab{f_{\xi}(t)-f_{\xi}(s)}\,,                                \\
        \label{ex-eind-b}
        \d_{\xi}(\be(t,x)-\be(s,x))
   &\le C_{\xi}\,(\nm{\xi}^2+\d_{\xi}(x))\,
        \ab{g_{\xi}(t)-g_{\xi}(s)}\,,                                \\
        \label{ex-eind-a1}
        \d_{\xi}(\al(t,x)-\al(t,y))
   &\le C_{\xi} \d_{\xi}(x-y)\,,                                     \\
        \label{ex-eind-b1}
        \d_{\xi}(\be(t,x)-\be(t,y))
   &\le C_{\xi} \d_{\xi}(x-y)\,. 
  \end{align}
  Here, $f_{\xi}$ and $g_{\xi}$ are continuous, real-valued functions
  on $\Rset^+$.
  
  If all these assumptions are satisfied, then there is a unique
  process $x$ in $\cV$ that solves the IDE
  \begin{align}\label{Sdgl}
        \d x_t 
   & =  \al(t,x_t)\d t + \d \cac_t \,\be(t,x_t) \,,
        \qquad x_{t_0}\in \LapXneu{(-\infty,t_{0}]}\,. 
  \end{align}
  Moreover, this solution $x$ is \nmot-continuous. In addition, if
  $x_{t_{0}}\in \LapXneu{[0,t_{0}]}$ and if \al, \be are
  locally adapted functions, then this solution $x$ is locally
  adapted.
\end{Thm}

\begin{Rem}\normalfont
  A similar theorem is valid for IDEs which contain both left and
  right non-commutative It\^o integrals.
\end{Rem}

\begin{proof}
  By the Picard iteration, we will construct on the interval
  $[t_0,t_1]$ a unique stop-continuous solution $x$ of the IDE
  \eqref{Sdgl}. In particular, we will produce from the initial data
  $(t_0, x_{t_0})$ the new data $(t_1, x_{t_1})$. The latter ones will
  serve as initial data for the Picard iteration to produce a solution
  on the interval $[t_1,t_2]$ (with $t_2 \le t_1+1$). By this iterative
  method, we will cover the interval $[t_0,t]$ by finitely many
  intervals $[t_i, t_{i+1}]$, since the iteration procedure will show
  that we can choose all intervals $[t_i, t_{i+1}]$ to be of the same
  length $0< \Delta t\le 1$. Thus, we will produce successively a
  unique solution for each interval $[t_i, t_{i+1}]$, and consequently
  a solution for any time $t\ge 0$.
  
  We already know from \eqref{ex-eind-a} and \eqref{ex-eind-a1} that,
  for a stop-continuous process $x$, the function $t\mapsto \al(t,x_t)$
  is stop-continuous and consequently \nmot-bounded on any compact
  interval. The same is true for the function $t\mapsto \be(t,x_t)$.
  
  We choose some fixed $\Delta t$ with $0 < \Delta t \le 1$ and start
  the iteration on the interval $[t_0,t_1]$ with $t_1:= t_0 + \Delta
  t$. Let $x_t^0:= x_{t_{0}}$ for any $t\in [t_0,t_1]$. The functions
  $t\mapsto\al(t,x_t^0)$ and $t\mapsto\be(t,x_t^0)$ are stop-continuous
  and (locally) adapted if \al and \be are (locally) adapted and if
  $x_{t_{0}}$ is an element in $L^{2}(\cA_{(-\infty,t_{0}]},\Eo )$
  (resp.~$L^{2}(\cA_{[0,t_{0}]},\Eo )$). Thus the $n$th iteration step
  \begin{align*}
        x_t^{n}  
   &:=  x_{t_{0}}
       +\int_{t_0}^{t} \!\d s  \,\al(s,x_s^{n-1}) 
       +\int_{t_0}^{t} \!\d \cac_s \be(s,x_s^{n-1})
  \end{align*}
  is well-defined by induction. Notice that $x^{n}$ is stop-continuous
  and (locally) adapted. We let
  $M:=\sup_{s\in[t_{0},t_{1}]}\nmo{x_s^1-x_s^0}$. Moreover, we define
  $q_{\xi} := C_{\xi}((\De t)^{\1/2} + \nm{\La(\1)}^{\!1/2})$.
  From \eqref{ex-eind-a1} and \eqref{ex-eind-b1} we find for any
  $t\in[t_{0},t_{1}]$ and normalized $\xi \in \cH_0$ the estimate
  \begin{align*}
        \d_{\xi}(x_t^{n+1}-x_t^n) 
   &\le \int_{t_0}^{t} 
        \d_{\xi}(\al(s,x_s^n)-\al(s,x_s^{n-1}))\d s                  \\
   &   +\nm{\La(\1)}^{1/2}
        \left[\int_{t_0}^{t} 
        \d_{\xi}(\be(s,x_s^n)-\be(s,x_s^{n-1}))^2\d s
        \right]^{\1/2}                                               \\
   &\le q_{\xi}^{} 
        \Big[
        \int_{t_{0}}^{t}
        \d_{\xi}(x_s^n-x_s^{n-1})^{2}\d s
        \Big]^{\1/2}
    \le M\,q_{\xi}^{n}\frac{(\De t)^{n/2}}{\sqrt{n!}}. 
  \end{align*}
  For any $n>q_{\xi}^{2}$ follows
  \begin{align*}
        \max_{t\in[t_{0},t_{1}]}
        \d_{\xi}(x_{t}^{n+m}-x_{t}^{n})
   &\le \sum_{k=0}^{m-1}
        \max_{t\in [t_0,t_1]}\d_{\xi}(x_t^{n+k+1}-x_t^{n+k})         \\
   &\le  M\,q_{\xi}^{n}\frac{(\De t)^{n/2}}{\sqrt{n!}}
        \sum_{k=0}^{m-1}
        \sqrt{\frac{n!}{(n+k)!}}(\De t)^{k/2}q_{\xi}^{k}             \\
   &\le  M\,q_{\xi}^{n}\frac{(\De t)^{n/2}}{\sqrt{n!}}
        \sum_{k=0}^{\infty}
        \Big(\frac{q_{\xi}}{\sqrt{n}}\Big)^{k}(\De t)^{k/2}
        \xrightarrow{n\to\infty}\,0\,.
  \end{align*}
  We conclude that the limit $x_{t}:=\lim_{n\to\infty}x_{t}^{n}$ exists
  uniformly on $[t_{0},t_{1}]$ in the stop-topology and defines a
  stop-continuous, (strongly) adapted process on this interval. Next,
  we obtain from \eqref{ex-eind-a1} and \eqref{ex-eind-b1} that
  $\int_{t_{0}}^{t}\al(s,x_{s})\d s =
  \lim_{n\to\infty}\int_{t_{0}}^{t}\al(s,x_{s}^{n})\d s$ and
  $\int_{t_{0}}^{t}\!\d \cac_{s} \be(s,x_{s}) =
  \lim_{n\to\infty}\int_{t_{0}}^{t}\!\d \cac_{s}\be(s,x_{s}^{n})$ in
  the strong operator topology. Now one concludes with routine
  arguments that $x$ solves the IDE on $[t_{0},t_{1}]$. Since
  $q_{\xi}$ depends only on $\De t$, this procedure applies iteratively
  to all right next neighbor intervals with the same length $\De t$.
  The \nmot-continuity of the solution $x$ follows from the fact that
  the integrals, as they appear in the integral equation for $x$,
  contain stop-continuous integrands and thus define \nmot-continuous
  functions. Finally, the uniqueness of the solution follows with the
  help of the estimates \eqref{ex-eind-a} - \eqref{ex-eind-b1} in the
  usual way.\qedopt
\end{proof}

\begin{Rem}\normalfont 
  A more complete treatment of non-commutative It\^o integration,
  compatible with the setting of continuous GNS Bernoulli shifts is
  contained in \cite{Hell01a}. If the shift is scalar-expected, then
  our approach reduces to non-commutative It\^o integration in the GNS
  Hilbert space of the shift, as it is already contained in
  \cite{Prin89a}. If the state of the $\Ao$-expected shift is tracial,
  then the continuous GNS Bernoulli shifts are realized as subspaces of
  non-commutative $L^2$-spaces in \cite{Koes00a}. This approach
  produces similar results on non-commutative It\^o integration. Notice
  also that one-sided integrands in \cite{BiSp98a} fall into our
  setting.
\end{Rem}

\section{Non-commutative exponentials and logarithms} 
\label{section:nc-ln-exp}

This section is devoted to the development of non-commutative
exponentials  (Subsection \ref{subsection:nc-exp}) and non-commutative
logarithms  (Subsection \ref{subsection:nc-ln}). The presented
constructions generalize beyond the frame of non-commutative
exponentials $\Exp$ of additive cocycles and non-commutative
logarithms $\Ln$ of unital cocycles. Here we will refrain to
include all these, in parts immediate generalizations and will just
focus on the development of sufficient tools to complete the proof of
Theorem \ref{theorem:ln-exp} and Theorem \ref{main-theorem}. This
will be done in  Subsection \ref{subsection:proof-correspondence},
where we will show that the mappings $\Exp$ and $\Ln$ are injective
and each other's inverse.

We assume throughout this section that a fixed $\Ao$-expected continuous
Bernoulli shift $\CBS$ and its GNS representation $\CBSE$ are given.

\subsection{Non-commutative exponentials of additive cocycles}
\label{subsection:nc-exp}

Non-commutative exponentials will be obtained as solutions 
of non-commutative It\^o differential equations (IDEs).

\begin{Thm}\label{sdgl-u}
  Let the centred additive cocycle $\cac \in \Ca(\cE,+)$ and the
  operator $K\in\Ao$ be given. Then the (unique) solution $u$ of
  the IDE
  \begin{align}
         u_{t} 
   & =  \1 + \int_{0}^{t}\!\d \cac_{s}  u_{s} 
        +\int_{0}^{t}\!\d t\,Ku_{s}\,.          \label{theorem-sdgl-u}
  \end{align}
  is a locally adapted \nmot-continuous process in $\Lapneu$ that
  satisfies the cocycle identity $u_{s+t}=(S_tu_s)u_t$ $(s,t \ge
  0)$. Moreover, it enjoys the following additional properties:
  \begin{enumerate}
  \item\label{sdgl-u-i}  The compression  
                         $A:=(\Eo u_{t})_{t\ge0}$ 
                         defines a
                         \nmt-continuous semigroup with 
                         $A_t = \e^{tK}$.
  \item\label{sdgl-u-ii} The compression 
                         $R := (\skpo{u_{t}}{\cdot\,u_{t}})_{t\ge0}$
                         defines a \nmt-continuous semigroup of 
                         completely positive mappings
                         on $\Ao$. It has the 
                         Christensen-Evans generator 
                         \[
                               \cL(a)
                            =  \La(a) + K^*a + aK,
                         \]
                         where $\La=\skpo{\cac_1}{\cdot\,\cac_1}$ is
                         the covariance of $b$.
  \end{enumerate}  
  In particular, the following are equivalent:
  { \renewcommand{\theenumi}{\alph{enumi}} 
  \begin{enumerate}
  \item\label{sdgl-u-a}  $u$ is a unital cocycle;
  \item\label{sdgl-u-b}  $R_{t}(\1) =\1$ for any $t\ge 0$, or
                         equivalently $\cL(\1) = 0$;
  \item\label{sdgl-u-c}  $\abO{\cac_{t}}^{2} +(K+K^{*})t =0$.
  \end{enumerate} }
  If the conditions (\ref{sdgl-u-a}) to (\ref{sdgl-u-c}) are satisfied,
  then the semigroup $R$ is contractive.
\end{Thm}

Let $\ncac$ denote the additive cocycle defined by $\ncac_t=
\cac_t +Kt$, $t\ge 0$.

\begin{Defn}\label{definition:nc-exp}
  The solution $u$ of the IDE \eqref{theorem-sdgl-u} is called the
  \emph{exponential of the additive cocycle $\ncac$} and
  is denoted by $\Exp(\ncac)$.
\end{Defn}

We remark that in stochastic analysis the additive cocycle $\ncac$ 
falls into the class of semi-martingales and that a solution of
\eqref{theorem-sdgl-u} is called an exponential semi-martingale.

\begin{proof}[Proof of Theorem \ref{sdgl-u}.]   
  Putting $\alpha(t,u_t) = K u_t$ and $\beta(t,u_t)= u_t$, the
  functions $\alpha$ and $\beta$ are locally adapted and evidently
  satisfy the Lipschitz conditions in Theorem \ref{ex-eind}. Moreover,
  the initial conditions is locally adapted, i.e., $u_{0} =\1\in
  \Ao$. Thus the IDE~\eqref{theorem-sdgl-u} has a unique
  \nmot-continuous locally adapted solution $u \subset \Lapneu$.
  
  In the following we will verify the cocycle property of the solution
  $u$. Let $(Z_n[a,b])_{n\in\Nset}$ be a sequence of partitions of the
  interval $[a,b]$ with grid tending to zero. In a first step we
  determine for an arbitrary $w\in \LapXneu{[0,t]}$
  \begin{align*}
        \Big(S_{t}\int_{0}^{s} \!\d \cac_r u_r\Big)w 
   & =  \lim_{n\to\infty}
        \sum_{r_{i}\in Z_n[0,s]}
        (S_t(\cac_{r_{i+1}}-\cac_{r_{i}})) (S_t u_{r_i})w          \\
   & =  \lim_{n\to\infty}
        \sum_{r_{i}\in Z_n[0,s]}  
        (\cac_{t+r_{i+1}}-\cac_{t+r_{i}}) (S_tu_{t+r_i-t})w        \\
   & =  \lim_{n\to\infty}
        \sum_{t_{i}\in Z_n[t,s+t]}
        (\cac_{t_{i+1}}-\cac_{t_{i}})\,(S_tu_{t_i-t})w             \\
   & =  \int_t^{t+s}\!\d \cac_r (S_tu_{r-t})w\,,
  \end{align*}
  where the limits are taken in the stop-topology on ${\Lapneu}$.
  Throughout these calculations we have used the continuity properties
  of the product of $\Ao$-independent elements in {\Lapneu}
  (Proposition~\ref{new-product}). Moreover, we used triple products
  (cf.~\ref{subsection:CBS-GNS}) of the three $\Ao$-independent,
  increasingly ordered factors $w$, $S_t u_{r_i}$ and
  $\cac_{t+r_{i+1}}-\cac_{t+r_{i}}$. (Here, `increasingly ordered'
  means $[0,t] \le [t,t+r_{i}] \le [t+r_{i},t+r_{i+1}]$.)  We define
  for a fixed $t>0$
  \begin{align*}
        v_r
    :=  \begin{cases}
         u_r             \,, & 0 \le r \le t, \\
         S_t(u_{r-t})u_t \,, & t  <  r      .
        \end{cases}
  \end{align*} 
  The continuity of the product of $\Ao$-independent elements implies
  the continuity of $r\mapsto v_r$ in the stop-topology and hence the
  integrability. Clearly, $v_r$ solves the IDE \eqref{theorem-sdgl-u},
  whenever $0\le r\le t$. Hence we conclude
  \begin{align*}
        v_{t+s} 
   & =  S_t\Big(
        \1 
       +\int_{0}^{s} \!\d \cac_r u_r 
       +\int_{0}^{s} \!\d r   Ku_r\Big)u_t                           \\
   & =  u_t 
       +\int_t^{t+s} \!\d \cac_r (S_t u_{r-t})u_t
       +\int_{0}^{s} \!\d r   K (S_t u_r)u_t                         \\
   & =  \1
       +\int_{0}^{t} \!\d \cac_r u_r 
       +\int_t^{t+s} \!\d \cac_r (S_t u_{r-t})u_t
       +\int_{0}^{t} \!\d r   Ku_r 
       +\int_t^{t+s} \!\d r   K (S_tu_{r-t})u_t                      \\
   & =  \1 
       +\int_0^{t+s} \!\d \cac_r v_r 
       +\int_0^{t+s} \!\d r   Kv_r\,.
  \end{align*}
  This calculation shows that $v$ solves \eqref{theorem-sdgl-u}. Now,
  the uniqueness of the solution implies $u_{r}=v_{r}$ for any $r\ge 0$
  and consequently the cocycle identity
  $u_{t+s}=v_{t+s}=(S_{t}u_{s})u_{t}$ of the solution $u$.
  
  (\ref{sdgl-u-i}) The compression $\Eo  u$ of the solution $u$
  defines a semigroup. Indeed, we apply $\Eo $ on both sides of
  \eqref{theorem-sdgl-u} and obtain
  \begin{align*}
        \Eo u_{t}
     =  \1 + \int_{0}^{t}\!\d s\,K\Eo u_{s}\,.
  \end{align*} 
  This integral equation has the unique solution $t\mapsto \Eo u_{t} =
  \e^{Kt}.$
  
  (\ref{sdgl-u-ii}) $R_{t}$ is completely positive by construction. By
  the cocycle property, it is shown immediately that $R_{s+t}(a)=
  R_s(R_t (a))$ for any $s,t \ge 0$ and $a\in \Ao$. Since $u$ is
  \nmot-continuous and $u_{0}=\1$, the uniform continuity of $R$
  follows from \eqref{cont-u-R}. Thus $R$ has a bounded generator
  $\cL$ (such that $R_t = \e^{t \cL}$). Next, we will identify the
  form of $\cL$. For this purpose we rewrite $R_t$ with the help of
  the IDE \eqref{theorem-sdgl-u}. Recall that
  $\La=\skpo{\cac_{1}}{\cdot\,\cac_{1}}$. An elementary calculation
  shows that for any $a\in\Ao$
  \begin{alignat*}{3}
          R_{t}(a)
   &={}&& \lefteqn{\skpo{u_t}{a u_t}}&&                              \\
   &={}&& a{}+{}                     && 
          \int_0^t \skpo{\1}{a K u_s}\d s
         +\int_0^t \skpo{Ku_s}{a}\d s                               
         +\int_0^t \skpo{u_s}{\La(a)u_s}\d s                         \\
   &   && \phantom{a}{}+{}           &&
          \int_0^t \Bskpo{K u_s}{a\int_0^s \!\d \cac_r u_r}\d s
         +\int_0^t \Bskpo{\int_0^s \!\d \cac_r u_r}{a K u_s}\d s     \\
   &   && \phantom{a}{}+{}           &&
          \int_0^t\int_0^s \skpo{Ku_s}{a Ku_r}\d r\d s
         +\int_0^t\int_0^r \skpo{Ku_s}{a Ku_r}\d s\d r\,. 
  \end{alignat*}
  Notice that, due to the \nmot-continuity of $u$, all integrals of the
  form $\int \skpo{\cdot}{\cdot}\d s$ are Bochner integrals on $\Ao$.
  Consequently, we are allowed to differentiate separately each term 
  and obtain for any $a\in \Ao$
  \begin{alignat*}{2}  
             \frac{\d}{\d t}R_t(a) 
   &{}={}& & \skpo{\1}{a K u_t} + \skpo{Ku_t}{a} 
            +\skpo{u_t}{\La(a)u_t}                                   \\
   &     & &+\Bskpo{Ku_t}{a\int_0^t \!\d \cac_r u_r}
            +\Bskpo{\int_0^t \!\d \cac_r u_r}{a K u_t}               \\
   &     & &+\Bskpo{Ku_t}{a\int_0^t \!\d r\,K u_r}
            +\Bskpo{\int_0^t \!\d r\,K u_r}{a K u_t}                 \\
   \intertext{-- we use again the IDE \eqref{theorem-sdgl-u} --}
   &{}={}& & \skpo{\1}{a K u_t} + \skpo{Ku_t}{a} +R_t(\La(a))        \\
   &     & &+\skpo{Ku_t}{a(u_t-\1)} + \skpo{u_t-\1}{a K u_t}         \\
   &{}={}& & R_t(\La(a)) + R_t(K^*a + a K)\,.
  \end{alignat*}
  Thus, the generator $\cL$ is identified as $\cL(a) = \La(a) + (K^*a +
  a K)$, where $a \in \Ao$.
  
  Finally, the equivalence of (\ref{sdgl-u-a}) to (\ref{sdgl-u-c}) and
  the contractivity of $R$ is evident from the form of the generator
  $\cL$ and the definition of a unital cocycle.\qedopt
\end{proof} 

\subsection{Non-commutative logarithms of unital cocycles}
\label{subsection:nc-ln}
We will construct an additive cocycle in $\Lapneu$ from a
\nmot-continuous unital cocycle. Motivated by Corollary
\ref{stone} and probability theory, we shall call the constructed
additive cocycle the non-commutative  logarithm of the unital
cocycle. Our main result of this subsection is stated in Theorem
\ref{unit-add-main}.

\begin{Notation}\normalfont
  $\cZ(t)$ will denote a net of partitions $Z:=\set{\tii\ge0}{0 = t_0 <
    t_1 < \ldots < t_{n_{Z}} = t}$ of the interval $[0,t]$. The set of
  partitions in $\cZ(t)$ is partially ordered by inclusion such that
  their grid $\ab{Z}:=\max\set{|\tiemti|}{i=0,\ldots,n_{Z}-1}$ tends to
  zero.
\end{Notation}

\begin{Thm}\label{unit-add-main}
  A \nmot-continuous unital cocycle $u\subset{\Lapneu}$ with the
  associated contractive semigroup $A_{t}:=\Eo u_{t}=\e^{Kt}$ defines
  via
  \begin{align}       
        \ncac_{t} 
   &:=  \nmolim{Z\in\cZ(t)}\sum_{i\,=\,0}^{n_{Z}-1} 
        S_{\tii}(u_{\tiemti}-\1),\\
        \cac_{t}   
   &:=  \nmolim{Z\in\cZ(t)}\sum_{i\,=\,0}^{n_{Z}-1}\label{std-constr-c}
        S_{\tii}u_{\tiemti}- A_{\tiemti}
  \end{align}
  two additive cocycles $\ncac,\cac$ in $\Ca(\cE,+)$. They are
 related  uniquely by $\ncac_t = \cac_t + K t$ and $K=\Eo \ncac_{1}$.
  Moreover, the additive cocycle $\ncac$ is in $\Cao(\cE,+)$, in
  other words, it satisfies the structure equation
  \begin{align*}
        \abO{\ncac_{t}-\Eo \ncac_{t}}^{2}
       +t(\Eo  \ncac_{1})^* + t\Eo  \ncac_{1}  
     =  0
  \end{align*}
  and the pair $(\cac,K)$, containing the centred additive cocycle
  $\cac$ and the drift $K$, satisfies the structure equation
  \begin{align*}
        \abO{\cac_{t}}^{2} +t(K^{*}+K)
     =  0\,.
  \end{align*}
\end{Thm}

\begin{Defn}\label{definition:nc-ln}
  The additive cocycle $\ncac$ (resp.\ $\cac$), constructed in
  Theorem \ref{unit-add-main}, is called the \emph{(centred)
  non-commutative logarithm of the unital cocycle} $u$ and
  denoted by $\Ln(u)$ (resp. by $\Lno(u)$).
\end{Defn}

For brevity, we will also say that $\ncac$ is the logarithm and that
$\cac$ is the centred logarithm of $u$. Notice that for a trivial
unital cocycle $u$, i.e.~$u_t=\exp{(itH)}$ with $H=H^*\in \Ao$,
one finds $\Ln(u_t)= itH$ and $\Lno(u_t)=0$. Let us further motivate
this definition:

\begin{Cor} \label{corollar-Ln}
  Let $u,v \subset \cA$ be two \nmot-continuous unital cocycles. If
  $u$ and $v$ are $\Ao$-independent and satisfy the commutation
  relation $u_t(S_tv_s) = (S_tv_s)u_t$ for any $s,t\ge 0$, then
  $uv:=(u_{t}v_{t})_{t\ge0}$ is again a unital \nmot-continuous
  cocycle and
  \begin{align*}
        \Ln(uv)
   & =  \Ln(u)+\Ln(v) 
     =  \Ln(vu)\,,                                                   \\
        \Lno(uv)
   & =  \Lno(u)+\Lno(v)
     =  \Lno(vu)\,. 
  \end{align*}
\end{Cor}

\begin{proof}
  The $\Ao$-independence of $u$ and $v$ ensures that the product
  $u_{t}v_{t}$ is well-defined in {\Lapneu} for any $t\ge0$. The
  commutation property for $u$ and $v$ guarantees that $uv$ is a
  cocycle. The \nmot-continuity of $t\mapsto u_tv_t$ is also
  elementary to check. Due to
  \begin{align*}
        S_{\tii}(u_{\tiemti}v_{\tiemti}-u_0v_0)
   & =  S_{\tii}(u_{\tiemti}-u_0)v_0 + u_0S_{\tii}(v_{\tiemti}-v_0)  \\
   & \quad  +S_{\tii}((u_{\tiemti}-u_0)(v_{\tiemti}-v_0)), 
  \end{align*}
  it is sufficient to establish 
  \begin{align*}
        \nmolim{Z\in\cZ(t)}\ \sum_{i\,=\,0}^{n_{Z}-1} 
        S_{\tii}((u_{\tiemti}-u_0)(v_{\tiemti}-v_0))
     =  0. 
  \end{align*}
  Indeed, this follows from estimates, similar as they appear in the
  proof of Theorem \ref{unit-add-main}. We leave these details to the
  reader.\qedopt
\end{proof}

We prepare the proof of Theorem \ref{unit-add-main} with a technical
result.

\begin{Lem}\label{rst} 
  For $r,s,t\ge 0$ and $a\in\Ao$ hold the following identities:
 \begin{enumerate}
 \item\label{rst-i} $\skpo{S_r u_s}{au_t} = R_s(aA_{t-r-s})\,A_r\,,$
   \quad $0\le s\le t-r,$
 \item\label{rst-ii} $\skpo{S_r u_s}{a u_t} = R_{t-r}(A_{s-t+r}^*a)
   \,A_r\,,$ \quad $0\le t-r\le s.$
 \end{enumerate}
\end{Lem}

\begin{proof}
  For $t-r\ge 0$ one calculates
  \begin{align*}
        \skpo{S_r u_s}{a u_t}
     &=  \skpo{S_r u_s}{(S_r a u_{t-r})\,u_r}
     =  \skpo{\1}{\skpo{u_s}{a u_{t-r}}\,u_r}\\
     &=  \skpo{u_s}{a u_{t-r}}\,A_r^{}\,.
  \end{align*}
  In case (\ref{rst-i}) with $0\le s\le t-r$ we conclude further 
  \begin{align*}  
        \skpo{S_r u_s}{a u_t} 
   & =  \skpo{u_s}{a(S_s u_{t-r-s})\,u_s}\,A_r^{}                    \\
   & =  \skpo{u_s}{a\skpo{\1}{u_{t-r-s}}\,u_s}\,A_r^{} 
     =  R_s(a A_{t-r-s})\,A_r\,.
  \end{align*}
  In the case (\ref{rst-ii}) with $0\le t-r\le s\,$ one finds
  \begin{align*}  
        \skpo{S_r u_s}{a u_t}
   & =  \skpo{(S_{t-r}u_{s-(t-r)})\,u_{t-r}}{a u_{t-r}}\,A_r         \\
   & =  \skpo{u_{t-r}}{A_{s-t+r}^*a u_{t-r}^{}} \,A_r^{}
     =  R_{t-r}(A_{s-t+r}^*a)\,A_r^{}\,. \qedopteqn
\end{align*}
\end{proof}

\begin{proof}[Proof of Theorem~\ref{unit-add-main}.] 
  We let
  \begin{align}\label{unit-add-ansatz}
        \cac_{Z}(t) 
    :=  \sum_{i\,=\,0}^{n_{Z}-1} S_{\tii}u_{\tiemti}-A_{\tiemti}
        \in \LapXneu{[0,t]}\,
  \end{align}
  and will prove that $(\cac_{Z}(t))_{Z\in\cZ(t)}$ is a \nmot-Cauchy
  net in {\Lapneu}. This convergence implies immediately that
  \begin{align}\label{unit-add-ansatz-nc}
        \ncac_{Z}(t) 
    :=  \sum_{i\,=\,0}^{n_{Z}-1} S_{\tii}u_{\tiemti}-\1
        \in \LapXneu{[0,t]}\,
  \end{align}
  converges to $\ncac_t$, since the difference
  $\ncac_{Z}(t)-\cac_{Z}(t)$  converges evidently to $Kt$. Moreover the
  equivalence of the structure equations of $\ncac$ and $\cac$ is
  ensured, since the decomposition of an additive cocycle in its
  centred part and its drift is unique. Thus, we will concentrate in
  the following on the convergence of $\cac_{Z}(t)$ to $\cac_t$, the
  cocycle property of $\cac_t$ and its structure equation.
  
  Let us denote the joint refinement of $Z,W\in\cZ(t)$ by $ZW$. Since
  \begin{align}\label{unit-add-unique}
        \nmo{\cac_{Z}(t)-\cac_{W}(t)} 
    \le \nmo{\cac_{Z}(t)-\cac_{ZW}(t)} 
       +\nmo{\cac_{W}(t)-\cac_{ZW}(t)},
  \end{align}
  it is sufficient to investigate
  $\nmO{\cac_{Z}(t)-\cac_{ZW}(t)}^2$. This expression is controlled by
  $\skpo{\cac_{Z}(t)}{\cac_{ZW}(t)}$. For the refinement $ZW$ of $Z$
  we let $t_i + \sij$, $j=1,\ldots,n^i-1$, be the additional points in
  the interval $[t_i,t_{i+1}]$. Moreover, we put $s_{i,0}:=0$ and
  $s_{i,n^{i}}:=t_{i+1}-t_{i}$. Setting $v_{s}$ $:=$ $u_{s} - A_{s}$
  for $s\in [0,t]$, one calculates
  \begin{align*}
        \lefteqn{\skpo{\cac_{Z}(t)}{\cac_{ZW}(t)}
     =  \sum_{i=0}^{n_{Z}-1}
        \sum_{k=0}^{n_{Z}-1}
        \sum_{j=0}^{n^k-1}\,
        \skpo{S_{\tii}v_{\tiemti}}{
        S_{t_k+s_{k,j}}v_{s_{k,j+1}-s_{k,j}}}}\hspace{1.5cm} 
   &                                                                 \\
   & =  \sum_{i=0}^{n_{Z}-1}\sum_{j=0}^{n^i-1}\,
        \skpo{u_{\tiemti}-A_{\tiemti}}{
        S_{\sij}u_{\sijemij}-A_{\sijemij}}                           \\
   & =  \sum_{i=0}^{n_{Z}-1}\sum_{j=0}^{n^i-1}\,
        \skpo{u_{\tiemti}}{S_{\sij}u_{\sijemij}}
       -A_{\tiemti}^*A_{\sijemij}^{}                                 \\
   & =  \sum_{i=0}^{n_{Z}-1}\sum_{j=0}^{n^i-1}\,
        A_{\sij}^*R_{\sijemij}^{}(A_{\tiemti-\sije}^*)
       -A_{\tiemti}^*A_{\sijemij}^{} \,                              \\
   \intertext{-- $\sijemij\le\tiemti-\sij$ admits the 
        application of Lemma \ref{rst}\spc(\ref{rst-i}) --} 
   & =  \sum_{i=0}^{n_{Z}-1}\sum_{j=0}^{n^i-1}\,
        A_{\sij}^*(R_{\sijemij}^{}-\id)(A_{\tiemti-\sije}^*-\1)      \\
   &   +\sum_{i=0}^{n_{Z}-1}\sum_{j=0}^{n^i-1}\,
        A_{(\tiemti)-(\sijemij)}^*(\1-A_{\sijemij}^*A_{\sijemij}^{})\,.
  \end{align*}
  We will investigate the two double sums separately. 
  
  From the norm continuity of the semigroups $R$ and $A$ we conclude
  that there exists an upper bound $M > 0$ such that $\nm{R_{t}-\id}
  \le M t$, $\nm{A_{t}-\1} \le M t$ and
  $\nm{A_{t}^*A_{t}^{}-\1} \le M t$. With this bound we produce
  the following estimate of the first summand:
  \begin{align*} 
        \lefteqn{\Bnm{\sum_{i=0}^{n_{Z}-1}\sum_{j=0}^{n^i-1}\,
        A_{\sij}^*(R_{\sijemij}^{}-\id)(A_{\tiemti-\sije}^*-\1)}} 
        \hspace{2.5cm} 
   &                                                                 \\
   &\le M^2
        \sum_{i=0}^{n_{Z}-1}\sum_{j=0}^{n^i-1}\,
        (\sijemij)(\tiemti-\sije)                                    \\
   &\le M^2
        \sum_{i=0}^{n_{Z}-1}\sum_{j=0}^{n^i-1}\,
        (\sijemij)(\tiemti)                                          \\
   & =  M^2
        \sum_{i=0}^{n_{Z}-1}\,(\tiemti)^2                              
    \le M^2\max\,\{\ab{Z},\ab{W}\} t\,.
  \end{align*}
  The second summand converges to $\la t$ with $\la:= -(K+K^*)$:
  \begin{alignat*}{2}
           \lefteqn{
           \Bnm{\sum_{i=0}^{n_{Z}-1} \sum_{j=0}^{n^i-1}\,
           A_{(\tiemti)-(\sijemij)}^*
           (\1 - A_{\sijemij}^*A_{\sijemij}^{}) - \la t }} 
           \hspace{1.5cm} 
   &   & &                                                           \\
   &\le& & \sum_{i=0}^{n_{Z}-1}\sum_{j=0}^{n^i-1}\,
           \bnm{A_{(\tiemti)-(\sijemij)}^*-\1}
           \bnm{\1 - A_{\sijemij}^*A_{\sijemij}^{}}                  \\
   &   &+& \sum_{i=0}^{n_{Z}-1}\sum_{j=0}^{n^i-1}\,
           \Bnm{\frac{\1-A_{\sijemij}^*A_{\sijemij}^{}}{\sijemij}-\la}
           (\sijemij)                                                \\
   &\le& & \sum_{i=0}^{n_{Z}-1}\sum_{j=0}^{n^i-1}\,
           M^2[(\tiemti)-(\sijemij)](\sijemij)                       \\
   &   &+& \sum_{i=0}^{n_{Z}-1}\sum_{j=0}^{n^i-1}\,
           \ve (\sijemij)                                            \\
   &\le& & \sum_{i=0}^{n_{Z}-1}\sum_{j=0}^{n^i-1}\,
            M^2(\tiemti)(\sijemij) + \ve t                           \\
   & = & & \sum_{i=0}^{n_{Z}-1}\,M^2(\tiemti)^2 + \ve t              \\
   &\le& & (M^2\ab{Z} + \ve) t
    \le    (M^2\max\,\{\ab{Z},\ab{W}\} + \ve) t \,.
  \end{alignat*}
  During this calculations we used the fact that $\bnm{\frac{\1 -
  A_{\de}^*A_{\de}^{} }{\de}-\la }\le\ve$ for any $0<\de\le\de_{\ve}$,
  with some appropriate $\de_{\ve}>0$ such that
  $\max\,\{\ab{Z},\ab{W}\} <\de_{\ve}$. Consequently, we obtain for
  any partitions $Z,W\in\cZ(t)$ with
  $\max\,\{\ab{Z},\ab{W}\}<\min\{\de_{\ve},\ve\}$ the estimates
  \begin{align*}
         \bnm{\skpo{\cac_{Z}(t)}{\cac_{ZW}(t)}-\la t}
    \le  \ve (M^2+1) t
  \end{align*} 
  and 
  \begin{align*}
         \nmO{\cac_{Z}(t)-\cac_{ZW}(t)}^2 
    \le  4\ve(M^2+1) t\,.
  \end{align*}
  A similar estimate is produced along the same line of arguments for
  $\nmO{\cac_{W}(t)-\cac_{ZW}(t)}^2$. According to the inequality
  \eqref{unit-add-unique}, we conclude that
  $(\cac_{Z}(t))_{Z\in\cZ(t)}$ is a Cauchy net with limit $\cac_t: =
  \nmolim{Z\in\cZ(t)} \cac_{Z}(t)\in \LapXneu{[0,t]}$. The independence
  of the limit $\cac_{t}$ from the used net $\cZ(t)$ is concluded
  immediately with the help of this inequality.
 
  We are left to prove the cocycle property of
  $\cac:=(\cac_{t})_{t\ge0}$. Given the nets $\cZ(s)$ and $\cZ(t)$ of
  partitions of the intervals $[0,s]$ resp.\ $[0,t]$, we define the
  net $\cZ(s+t)$ of partitions associated to the interval $[0,s+t]$ in
  the following manner. For $Z_{s} :=
  \set{s_i\ge0}{0=s_0<s_1<\ldots<s_{n_s}=s}\in\cZ(s)$ let $t+Z_{s}$
  denote the partition $\set{t+s_i}{0=s_0<s_1<\ldots<s_{n_s}=s}$ of
  the interval $[t,t+s]$. Now we define for any pair $(Z_{t},Z_{s})$
  of partitions $Z_{t}\in\cZ(t)$ and $Z_{s}\in\cZ(s)$ an element $Z$
  of the net $\cZ(t+s)$ by $Z := Z_t \cup (t+Z_s)$. Obviously, the
  grid of this net tends to zero. We observe
  \begin{align*}
        \cac_{Z}(t+s) 
     &=  \sum_{i=0}^{n_t-1}\,S_{\tii}(v_{\tiemti})
       +\sum_{i=0}^{n_s-1}\,S_{t+s_i}(v_{s_{i+1}-s_i})\\ 
     &=  \cac_{Z_t}(t) + S_t \cac_{Z_s}(s)\,.
  \end{align*}
  From the net convergence of the left-hand side of this equality and
  the convergence of each summand of the right-hand side to
  $\cac_{t+s}$, $\cac_t$ resp.\ $S_{t}\cac_s$, we obtain the cocycle
  identity $\cac_{t+s} = \cac_t + S_{t}\cac_s$.
 
  We have already identified in Corollary \ref{CE-generator} the
  Christensen-Evans generator $\cL$ of $R$ as
  \begin{align*}
        \cL(a)
     =  \La(a) + K^{*}a + a K\,.
  \end{align*}
  Since $\cL(\1)=0$ and $\abO{\cac_{t}}^{2} = t\abO{\cac_{1}}^{2}$ by
  \eqref{bt-abs}, the structure equation follows immediately from
  \begin{align*}
        \La(\1)
     =  \abO{\cac_{1}}^{2} 
     = -(K + K^*)\,.\qedopteqn
  \end{align*}
\end{proof}   

We close this section with a technical result which will be needed in
Section \ref{subsection:proof-correspondence} to finish the proof of
Theorem \ref{main-theorem}.

\begin{Lem}\label{b-u}
  Let $u$ be a \nmot-continuous unital cocycle and $\cac:= \Lno(u)$
  the centred logarithm of $u$. Setting $A_{t}:=\Eo u_{t}$ and
  $\La:=\skpo{\cac_{1}}{\cdot\,\cac_{1}}$, it follows
  \begin{align}\label{b-a-u}
        \skpo{\cac_t}{au_t}
   & =  \int_{0}^{t}\La(aA_{t-s})A_s\d s\,,
        \qquad a\in\Ao\,.
  \end{align}
\end{Lem}

\begin{proof}
  We recall Lemma~\ref{rst}\spc(\ref{rst-i}) to see that $\skpo{S_r
    u_s}{a u_t} = R_s(aA_{t-r-s})A_r$, whenever $0\le s\le t-r$ and
  $a\in\Ao$. Next, we choose for $[0,t]$ the equidistant partition
  $Z_{n}:=\set{i\de_{n}}{\de_{n}:=2^{-n}t,i=0,\ldots ,2^{n}}\in\cZ(t)$.
  Since the approximants $\cac_{Z_n}$ in \eqref{std-constr-c} converge
  to $\cac_t$ in the \nmot-topology, we obtain in the uniform topology
  on $\Ao$ that
  \begin{alignat*}{2}
           \skpo{\cac_t}{a u_t}
   & = & & \lim_{n\to\infty}
           \sum_{i=0}^{2^n-1}
           \skpo{S_{i\de_n}u_{\de_n}-A_{\de_n}}{au_t}                \\
   & = & & \lim_{n\to\infty}
           \sum_{i=0}^{2^n-1} 
           \big[
           R_{\de_n}(aA_{t-(i+1)\de_n})A_{i\de_n}
          -A_{\de_n}^*aA_{\de_n}^{}
           \big]                                                     \\
   & = & & \lim_{n\to\infty}
           \sum_{i=0}^{2^n-1} 
           \Big[
           \frac{R_{\de_n}-id}{\de_n}-\cL
           \Big]
           (aA_{t-(i+1)\de_n}) A_{i\de_n}\de_n                       \\
   &   & &+\lim_{n\to\infty}
           \sum_{i=0}^{2^n-1} 
           \cL(aA_{t-(i+1)\de_n})A_{i\de_n}\de_n                     \\
   &   & &+\lim_{n\to\infty}
           t
           \left[
           a \frac{A_{t-\de_n}-A_t}{\de_n} 
          +\frac{\1-A_{\de_n}^*}{\de_n} a A_t
           \right]                                                   \\
   & = & & \int_{0}^{t}\cL(aA_{t-s})A_s\d s 
          -t (aK+K^*a)A_t                                            \\
   \intertext{-- we use the form of the generator 
                $\cL$ from Corollary \ref{CE-generator} --} 
   & = & & \int_0^t \La(aA_{t-s})A_s\d s + t K^*aA_t                 \\
   &   & &  +\int_0^t aA_{t-s}KA_s\d s - t(aK+K^*a)A_t               \\
   & = & &        \int_0^t \La(aA_{t-s})A_s\d s\,.\qedopteqn
  \end{alignat*}
\end{proof}

\subsection{Proof of the correspondence}
\label{subsection:proof-correspondence}

In this section we will complete the proof of Theorem
\ref{theorem:ln-exp} and Theorem \ref{main-theorem}. We will show
that the mappings $\Ln$ and $\Exp$, as introduced in Section
\ref{subsection:nc-ln} resp.~Section \ref{subsection:nc-exp}, are
injective. This will show that $\Ln$ and $\Exp$ are each-other inverse.
Thus we will have completed the proof of Theorem \ref{theorem:ln-exp}.
Notice that the injectivity of both mappings will also provide the
bijectivity  of the correspondence in Theorem \ref{main-theorem}. We
divide its proof in  several intermediate results.

Let the $\Ao$-expected continuous GNS Bernoulli shift \CBSE be given. 

\begin{Prop} \label{Exp-injectivity}
  For any additive cocycle $\ncac$ in $\Cao(\Lapneu,+)$ holds 
  \begin{align*}
        \Ln(\Exp(\ncac))
     =  \ncac\,.
  \end{align*} 
  Consequently, the mapping $\Exp\colon \Cao(\Lapneu,+) \to
  \Cmo(\Lapneu,\boldsymbol{\cdot}\,)$ is injective.
 \end{Prop}

We will need the following Lemma for the proof of Proposition
\ref{Exp-injectivity}.

\begin{Lem}\label{u-formeln}
  Let $\ncac$ be an additive cocycle with centred part $\cac$ and
  drift $K$, satisfying the structure equation $\skpo{\cac_t}{\cac_t}+
  t(K^* + K)=0$. Let $u_t := \Exp(\cac_t +Kt)$ and $A_t :=
  \e^{Kt}$. For any $a\in\Ao$ holds
  \begin{align*} 
        \Bskpo{\cac_t}{a\int_{0}^{t}\!\d \cac_s u_s} 
     =  \int_{0}^{t} \La(a)A_s \d s\,.
  \end{align*}
\end{Lem} 

\begin{proof}
  Since $u$ is the solution of the IDE \eqref{theorem-sdgl-u}, the
  map $t\mapsto u_t$ is \nmot-continuous and thus the It\^o integral
  $\int_{0}^{t}\!\d \cac_s u_s$ is well-defined. Moreover,
  $(a^{*}\cac_{t})_{t\ge0}$ is a centred additive cocycle. We
  calculate with the It\^o identity \eqref{ito-identity}
  \begin{align*}
        \Bskpo{\cac_t}{a \int_0^t \!\d \cac_s  u_s}
   & =  \Bskpo{\int_0^t\! \d(a^{*}\cac_s)}{\int_0^t \!\d \cac_s  u_s}
     =  \int_0^t \skpo{\1}{\La(a)u_{s}}\d s                          \\
   & =  \int_{0}^{t}\La(a)A_s \d s\,.                        \qedopteqn
  \end{align*}
\end{proof}

\begin{proof}[Proof of Proposition \ref{Exp-injectivity}.] 
  By Theorem \ref{sdgl-u} $u:= \Exp(c)$ is a unital cocycle and by
  Theorem~\ref{unit-add-main} we know that $\Ln(\Exp(c))$ is again
  an
  additive cocycle satisfying the structure equation. We are left
  with the task to identify this cocycle as $c$. For this we
  consider, as stated in the previous Lemma, the unique decomposition
  $\ncac_t = \cac_t + Kt$ of the additive cocycle $\ncac$ and $A_t
  = \Eo  \Exp(\ncac_t)$. We will prove
  \begin{align*}
        \Ln(\Exp(\ncac_t))-\ncac_t
     =  \lim_{n\to \infty} \sum_{i=0}^{2^n-1}S_{i\de_n}
        (\Exp(\ncac_{\de_n})-\1)-\ncac_{t} 
     =  0\,, 
  \end{align*} 
  where  $\de_n:=t2^{-n}$, by showing 
  \begin{align}\label{exp-1}
         \lim_{n\to \infty} 
         \sum_{i=0}^{2^n-1}S_{i\de_n} (A_{\de_n}-\1) - K t 
     =   0
  \end{align}
  and 
  \begin{align}\label{exp-2}
        \lim_{n\to \infty} \sum_{i=0}^{2^n-1}S_{i\de_n}
        (\Exp(\ncac_{\de_n})-A_{\de_n}) - \cac_{t} 
     =  0
  \end{align}
  in the \nmot-topology. The equation \eqref{exp-1} is obvious, since
  $\Ao$ is the fixed point algebra of $S$ and $A_t$ is a uniformly
  continuous semigroup with generator $K$. Next we focus on the limit
  stated in equation \eqref{exp-2} and let $u_t:= \Exp(\ncac_t)$ for
  shortness:
  \begin{align*}
   \lefteqn{\sum_{i=0}^{2^n-1}
        S_{i\de_n}(u_{\de_n}-A_{\de_n}) - \cac_{t}
     =  \sum_{i=0}^{2^n-1}
        S_{i\de_n}(u_{\de_n}-\1-\cac_{\de_n})
       -\sum_{i=0}^{2^n-1} (A_{\de_n}-\1)}\hspace{0.5cm}
   &                                                                 \\
   & =  \sum_{i=0}^{2^n-1}
        S_{i\de_n} \int_{0}^{\de_n}\d \cac_s (u_s-\1)
       -\sum_{i=0}^{2^n-1}
        S_{i\de_n}\int_{0}^{\de_n} \d s\,Ku_s
       -\sum_{i=0}^{2^n-1}
        \frac{A_{\de_n}-\1}{\de_n}\de_n\,. 
  \end{align*}
  Obviously, the last summand tends to $-Kt$. The second summand
  tends to $Kt$, since 
  \begin{align*}
         \Bnmo{\int_0^{\de_n}K(u_{\de_n}-\1)\d s} 
    \le  \sqrt{2}\nm{K}\int_0^{\de_n}\nm{A_s-\1}^{1/2}\d s  
  \end{align*}
  and thus $\int_0^{\de_n}K(u_{\de_n}-\1)\d s$ tends to zero with
  order $\de_n^{3/2}$. Finally, the first summand converges to zero: 
  \begin{alignat*}{2}
    \lefteqn{\BabO{\sum_{i\,=\,0}^{2^n-1}
           S_{i\de_n}\int_{0}^{\de_n}\d \cac_s\,(u_s - \1)}^2
     =     \sum_{i\,=\,0}^{2^n-1}
           \BabO{\int_{0}^{\de_n}\d \cac_s\,u_s - \cac_{\de_{n}})}^2}
           \hspace{1.5cm}
   &   & &                                                           \\
   & = & & \sum_{i=0}^{2^n-1}
           \Big(
           \int_{0}^{\de_n}R_s(\La(\1))\d s + \de_n\La(\1)
           \Big)
          -\sum_{i=0}^{2^n-1}2\RE
           \Bskpo{\cac_{\de_n}}{\int_{0}^{\de_n}\!\d \cac_s\,u_s}    \\
   \intertext{-- we use Lemma \ref{u-formeln} --}
   & = & & \sum_{i=0}^{2^n-1}
           \int_{0}^{\de_n}R_s(\La(\1))\d s + \La(\1) t 
          -\sum_{i=0}^{2^n-1}2\RE
           \int_{0}^{\de_n}\La(\1) A_s\d s
  \end{alignat*}
  which tends to zero for $n \to \infty$, since both sums converge to
  $2\La(\1)t$ in the \nmt-topology on $\Ao$.\qedopt
\end{proof}   

We prove next the injectivity of the mapping $\Ln$.

\begin{Prop} \label{Ln-injectivity}\label{u-sdgl}
  For any unital cocycle $u$ in $\Cmo(\Lapneu,\boldsymbol{\cdot}\,)$ 
  holds
  \begin{align*}
        \Exp(\Ln(u))
     =  u\,.
  \end{align*} 
  Consequently, the mapping $\Ln\colon
  \Cmo(\Lapneu,\boldsymbol{\cdot}\,) \to \Cao(\Lapneu,+) $ is
  injective.
\end{Prop}

We will show that each \nmot-continuous unital cocycle $u$ is a
solution of the IDE~\eqref{main-t-dgl}, where $\cac = \Ln_0(u)$, the
centred non-commutative logarithm of $u$.

\begin{proof}
  We prepare the proof by some results which we will need in the
  sequel. By Theorem~\ref{unit-add-main} we know already: the semigroup
  $R:=(\skpo{u_{t}}{\cdot\,u_{t}})_{t\ge0}$ has the
  generator $\Ao\ni a\mapsto \cL(a):=\La(a)+K^{*}a + aK$, where
  $\La:=\skpo{\cac_{1}}{\cdot\,\cac_{1}}$ and $A_{t}:=\Eo
  u_{t}=\e^{Kt}$, as well as $\La(\1)=-(K^{*}+K)$. The non-commutative
  logarithm of $u$ we obtain as $\cac_t =\nmolim{n\to\infty}
  \sum_{i=0}^{2^n-1} S_{i\de_n}(u_{\de_n}-A_{\de_n})$, with
  $\de_n:=t2^{-n}$. To prove the Lemma we will show
  \begin{align}\label{star}
        \BabO{u_t-\1-\int_0^t \!\d \cac_s u_s 
       -\int_0^t \!\d s\,Ku_s}^2 
     =  0\,.
  \end{align}
  The inner product leads to the following ten expressions:

  \begin{xalignat*}{5}
     (i) &\quad   \abO{u_t}^2 = \1      \quad &
    (ii) &\quad   \abO{\1}^2  = \1            &            & & &     \\
   (iii) &\quad   \BabO{\int_0^t\! \d \cac_s u_s}^2 
                = \int_0^t R_s(\La(\1))\d s                          \\
    (iv) &\quad   \lefteqn{\BabO{\int_0^t\! \d s\,Ku_s}^2 
                = \int_0^t R_s(\La(\1))\d s}
                                                                     \\
         &\quad   \qquad  \qquad  \qquad 
          \quad   \lefteqn{-2\RE\Big[\int_0^t R_s^{}(\La(A_{t-s}^{}))
                  \d s - (\1-A_t^{})\Big]}  
                                                                     \\
     (v) &\quad   \lefteqn{-2\RE \skpo{u_t}{\1} = -2\RE A_t^{}}
                                                                     \\
    (vi) &\quad   2\RE\Bskpo{\1}{\int_0^t \!\d \cac_s u_s} = 0
                                                                     \\
   (vii) &\quad   \lefteqn{-2\RE\Bskpo{u_t}{\int_0^t \!\d \cac_s u_s} 
               = -2\RE\int_0^t R_s^{}(\La(A_{t-s}^{}))\d s}
                                                                     \\
  (viii) &\quad   \lefteqn{-2\RE\Bskpo{u_t}{\int_0^t \!\d s\,Ku_s}
               =  2\RE\Big[\int_0^t R_s^{}(\La(A_{t-s}^{}))\d s
                 -(\1-A_t^{})\Big]}
                                                                     \\
    (ix) &\quad   \lefteqn{2\RE\Bskpo{\1}{\int_0^t \!\d s\,Ku_s} 
               = -2\RE(\1-A_t^{})}
                                                                     \\
     (x) &\quad   \lefteqn{2\RE\Bskpo{\int_0^t \!\d s\,Ku_s}{\int_0^t
                  \!\d \cac_su_s} 
               =  2\RE\Big[\int_0^t R_s^{}(\La(A_{t-s}^{}))\d s}     \\
         &\quad   & & \qquad \quad - \int_0^t R_s(\La(\1))\d s\Big]
  \end{xalignat*}
  Equation \eqref{star} is obtained by adding (i) to (x). The equations
  (i), (ii), (iii), (v) and (vi),  as well as  (ix), are evident.    \\
  Ad (iv):\hspace{2mm} We use Lemma \ref{rst}\spc(\ref{rst-i}) and
  (\ref{rst-ii}) to obtain
  \begin{align*}
        \BabO{\int_0^t \!\d s\,Ku_s}^2
   & =  \int_0^t\!\int_0^{r}
        \skpo{u_s}{K^{*}Ku_r}\d s\d r
       +\int_0^t\!\int_{r}^{t}
        \skpo{u_s}{K^{*}Ku_r}\d s\d r                                \\
   & =  \int_0^t\!\int_0^r 
        R_s^{}(K^*KA_{r-s}^{})\d s\d r
       +\int_0^t R_r^{}\Big(\int_{r}^{t}
        A_{s-r}^*K^* K\d s\Big)\d r\,.
  \end{align*}
  The second integrals evaluates to $\int_0^t
  R_r^{}((A_{t-r}^*-\1)K)\d r$. In the first integrals we exchange the
  order of integration and obtain $\int_0^t R_s^{}\big(\int_{s}^{t}
  K^*KA_{r-s}\d r\big)\d s = \int_0^t R_s^{}(K^*(A_{t-s}-\1))\d s $.
  We collect all expressions and obtain
  \begin{align*}    
        \BabO{\int_0^t \!\d s\,Ku_s}^2 
     =  \int_0^t R_s^{}(\La(\1))\d s
       +\int_0^t R_s^{}(A_{t-s}^*K + K^*A_{t-s}^{})\d s\,. 
  \end{align*}
  The next calculation will make use of
  \begin{align*}
       -\int_0^t R_s^{}(K^*A_{t-s}^{*})\d s
     =  \int_0^t R_s^{}\big(\frac{\d}{\d s}A_{t-s}^*\big)\d s 
     =  \1-A_t^* -\int_0^t R_s^{}(\cL(A_{t-s}^*))\d s\,.
  \end{align*}
  and of a similar equation for the adjoint expression  
  $-\int_{0}^{t} R_s^{}( A_{t-s}^{}K)\d s$. 
  We take into account these two equations and the form of the
  generator $\cL$ to find
  \begin{alignat*}{2}
             \BabO{\int_0^t \!\d s\,Ku_s}^2 
   &  =  & & \int_0^t R_s^{}(\La(\1))\d s                            \\
   &     &+& \int_0^t R_s^{}(A_{t-s}^*K+K^*A_{t-s}^*)\d s
            +\int_0^t R_s^{}(A_{t-s}^{}K+K^*A_{t-s}^{})\d s          \\
   &     &-& \int_0^t R_s^{}(K^*A_{t-s}^*)\d s
            -\int_{0}^{t} R_s^{}( A_{t-s}^{}K)\d s                   \\
   &  =  & & \int_0^t R_s^{}(\La(\1))\d s                            \\
   &     &+& \int_0^t R_s^{}((A_{t-s}^{}+A_{t-s}^*)K
            +K^*(A_{t-s}^{}+A_{t-s}^*))\d s                          \\
   &     &-& \int_0^t R_s^{}(\cL(A_{t-s}^{}+A_{t-s}^*))\d s
            +\1-A_t^{} +\1-A_t^*                                     \\
   &{}={}& & \int_0^t R_s^{}(\La(\1))\d s
            -2\RE\int_0^t R_s^{}(\La(A_{t-s}^{}))\d s
            +2\RE(\1-A_t^{})\,. 
  \end{alignat*}
  Ad (vii):\hspace{2mm} Next of all we calculate for $0\le s\le t-\de$:
  \begin{alignat*}{2}
             \skpo{u_t^{}}{(S_s^{}\cac_{\de}^{})\,u_s^{}} 
   &{}={}& & \skpo{(S_s^{}u_{t-s}^{})u_s^{}}{(
             S_s^{}\cac_{\de}^{})u_s^{}}
      =      R_s^{}(\skpo{u_{t-s}^{}}{\cac_{\de}^{}})                \\
   &{}={}& & R_s^{}(\skpo{A_{t-s-\de}^{}u_{\de}^{}}{\cac_{\de}^{}})
   \stackrel{\makebox{\scriptsize\eqref{b-a-u}}}{=}
             R_s^{} \Big( \int_{0}^{\de} 
             A_r^*\La(A_{\de-r}^*A_{t-s-\de}^*)\d r \Big)            \\
   &{}={}& & R_s^{} \Big( 
             \int_0^{\de} A_r^*\La(A_{t-s-r}^*)\d r \Big)\,. 
  \end{alignat*}
  Thus we obtain:
  \begin{align*}
        \Bskpo{u_t^{}}{\int_0^t \!\d \cac_s^{}\,u_s^{}}
   & =  \lim_{n\to\infty}\sum_{i=0}^{2^n-1} 
        \skpo{u_t^{}}{(S_{i\de_n}^{}
        \cac_{\de_n}^{})u_{i\de_n}^{}}                               \\
   & =  \lim_{n\to\infty}\sum_{i=0}^{2^n-1} 
        \underbrace{
        R_{i\de_n}^{}
        \Big(\frac{1}{\de_n}\int_0^{\de_n}
        A_r^*\La(A_{t-i\de_n-r}^*)\d r - \La(A_{t-i\de_n}^*)
        \Big)}_{o(\de_n)}\de_n                                       \\
   &   +\lim_{n\to\infty}\sum_{i=0}^{2^n-1} 
        R_{i\de_n}^{}(\La(A_{t-i\de_n}^*))\de_n 
     =  \int_0^t R_s^{}(\La(A_{t-s}^*))\d s\,.
  \end{align*}
  Ad (viii):
  \begin{align*}
        \lefteqn{\Bskpo{u_t^{}}{\int_0^t \!\d s\,Ku_s^{}}
     =  \int_0^t \skpo{(S_s^{}u_{t-s}^{})u_s^{}}{Ku_s^{}}\d s 
     =  \int_0^t R_s^{}(A_{t-s}^*K)\d s}\hspace{2.0cm}
   &                                                                 \\
   & =  \int_0^t R_s^{}(A_{t-s}^*K + K^*A_{t-s}^*)\d s
       -\int_0^t R_s^{}(K^*A_{t-s}^*)\d s                            \\
   & =  \int_0^t R_s^{}\cL(A_{t-s}^*)\d s
       -\int_0^t R_s^{}\La(A_{t-s}^*)\d s
       -\int_0^t R_s^{}(K^*A_{t-s}^*)\d s                            \\
   & =  \1-A_t^* - \int_0^t R_s^{}(\La(A_{t-s}^*))\d s \,.
  \end{align*}
  The last equation is found, similar to (iv), by partial 
  integration.\\
  Ad (x):\hspace{2mm} The same calculation as for (vii) shows
  $\skpo{Ku_s^{}}{\int_0^s \!\d \cac_r^{}\,u_r^{}} = \int_0^s
  R_r^{}(\La(K^*A_{s-r}^*))\d r$. Consequently, we have
  \begin{align*}
        \Bskpo{\int_0^t \!\d s\,Ku_s^{}}{\int_0^t \!\d \cac_s^{}u_s^{}}
   & =  \int_0^t\int_0^s R_r^{}(\La(K^*A_{s-r}^*))\d r\d s           \\
   & =  \int_0^t\int_r^t R_r^{}(\La(K^*A_{s-r}^*))\d s\d r           \\
   & =  \int_0^t R_r^{}(\La(A_{t-r}^*))\d r 
       -\int_0^t R_r^{}(\La(\1))\d r \,.                     \qedopteqn
  \end{align*}
\end{proof}

Now we have provided all results to finish the proof of Theorem
\ref{theorem:ln-exp} and Theorem~\ref{main-theorem}.

\begin{proof}[Proof of Theorems~\ref{theorem:ln-exp} and
  \ref{main-theorem}.]                                      
  According to Theorem~\ref{sdgl-u} we can associate to each additive
  cocycle $\ncac$ in $\Cao(\Lapneu, +)$ a unital cocycle $u =
  \Exp(\ncac)$ in $\Cmo(\Lapneu,\boldsymbol{\cdot}\,)$. By
  Proposition \ref{Exp-injectivity} the mapping $\Exp$ is injective and
  by Proposition \ref{Ln-injectivity} it is also surjective. This
  completes the proof of the two theorems. \qedopt
\end{proof}                                                   
\appendix
\renewcommand{\theequation}{\thesection.\arabic{equation}}
\numberwithin{equation}{section}

\section{Hilbert W*-modules}                    \label{sec-appendix-hm}

For the convenience of the reader we provide some background results on
Hilbert W*-modules, as far as we will need them throughout this paper.

We start with a (pre-) Hilbert C*-module $\cE$ over a von Neumann
algebra $\Ao\subseteq\cB(\Ho)$ and present the construction of a
Hilbert W*-module which we will need within the framework of continuous
Bernoulli shifts. For a detailed approach to Hilbert modules we refer
to \cite{Lance1} and for the notion of Hilbert W*-modules to
\cite{Frank,Paschke,Schweizer}. The \Ao-valued inner product
$\skpo{\cdot}{\cdot}$ induces by $\abo{x}:=\skpO{x}{x}^{1/2}$ an
\Ao-valued `norm' which gives rise to the norm $\nmo{x}:=\nm{\abo{x}}$
on \cE. The completion of $\cE$ in this norm is a Hilbert C*-module
which we will also denote by $\cE$. The Cauchy-Schwarz inequality for
the inner product is valid in the following form:
\begin{align}\label{CS-ineq}
      \ab{\skpo{x}{y}}^2 
 &\le \nmO{x}^2\,\abO{y}^2, 
      \qquad\qquad\qquad 
      x,y\in\cE\,.
\end{align}
$\cB(\cE)$ is the Banach algebra of bounded module maps, i.e., the
continuous \Ao-linear maps $T\,\colon \,\cE\to\cE$ such that $T(xa) =
T(x)a\,$ for any $x\in\cE$ and $a\in\Ao$. A bounded linear map
$T:\cE\to\cE$ is $\Ao$-linear if and only if the inequality
\begin{align}\label{A-linear}
      \skpo{Tx}{Tx}
 &\le M\,\skpo{x}{x}, 
      \qquad\qquad\qquad 
      x\in\cE\,.
\end{align}
is satisfied \cite{Paschke}. $\cL(\cE)$ denotes the C*-algebra of
adjointable $\Ao$-linear maps, i.e., maps $T\in\cB(\cE)$ for which a
linear map $T^{*}$ exists such that $\skpo{x}{Ty}=\skpo{T^{*}x}{y}$ for
any $x,y\in\cE$. In general an \Ao-linear map is not adjointable (but
this property will always be the case for W*-modules). The topological
dual $\cE'$ of $\cE$ is given by the \Ao-linear maps from $\cE$ into
\Ao. The map $\hat{x}:\cE\ni y\mapsto\skp{x}{y}$ defines an isometric
embedding of $\cE$ into $\cE'$. Notice that, in general, the image
$\widehat{\cE}$ is not identical to $\cE'$. A Hilbert module $\cE$ with
the property $\widehat{\cE}=\cE'$ is called selfdual.

$\cK(\cE) := \mathrm{lh}\set{\Theta_{x,y}}{x,y\in\cE}^{-\nmt}
\subseteq\cL(\cE)$, with $\Theta_{x,y}z:=x\skpo{y}{z}$ is the two-sided
*-ideal of `compact operators' on \cE. The ideal $\cK(\cE)$ has an
approximate unit $(e_{\al})_{\al\in I}$ of the form $e_{\al} :=
\sum_{i=1}^{n_{\al}}\Theta_{y_{i}^{\al},y_{i}^{\al}}$ (which can be
constructed by an obvious adaption of the proof in
\cite[Prop. 2.2.18]{BrRo1}).\\
In the following, we present a concrete realization of the Hilbert
W*-module via the minimal Kolmogorov decomposition of the kernel
$\cE\times\cE\to\Ao$; $(x,y)\mapsto\skpo{x}{y}$ \cite{EvansLewis,
  Murp97a}. By this decomposition the Hilbert module $\cE$ is realized
as a subspace of $\cB(\Ho,\cH')$ for some Hilbert space $\cH'$
such that $\overline{\cE\Ho}=\cH'$ and thus, by an
embedding, as a subspace of $\cB(\Ho \oplus\cH')$. Consequently, we
assume that the inner product of $\cE$ takes the form $\skpo{x}{y}
=x^*y$.

\begin{DefnApp}\label{definition:Hilbert-W-module}
  The closure of $\cE$ in the stop topology of $\cB(\Ho\oplus\cH')$ is
  called a Hilbert W*-module.
\end{DefnApp}


The strong operator (stop) topology resp.\ the $\si$-strong operator
($\si$-stop) topology on $\cE$ is induced by the seminorms
$x\mapsto\nm{\abo{x}\xi_{0}}$, $\xi_{0}\in\Ho$, resp.\ 
$x\mapsto|\phi(\skpo{x}{x})|^{1/2}$, $\phi\in {\Ao}_{*}$. The weak*
topology on $\cE_{1}$ is induced by the seminorms
$x\mapsto\ab{\phi(\skpo{y}{x})}$, $\phi\in{\Ao}_{*}$, $y\in\cE$.

\begin{ThmApp}\label{self-dual}
  A Hilbert W*-module $\cE$ over the von Neumann algebra $\Ao$ has the
  predual $\cE_{*}=\mathrm{lh}\set{\phi(\skpo{y}{\cdot})}{y\in\cE,
    \phi\in{\Ao}_{*}}^{\text{--}\nmt}$ and is selfdual in the sense of
  Hilbert modules.
\end{ThmApp}

\begin{proof}
  Since $\cE$ is a weakly* closed subspace of $\cB(\Ho\oplus\cH')$, it
  is the dual of the Banach space $\cT(\Ho\oplus\cH')/\cE^{\circ}$.
  Here, $\cT(\Ho\oplus\cH')$ denotes the trace class operators on
  $\Ho\oplus\cH'$ and $\cE^{\circ}$ the polar of $\cE$ in
  $\cT(\Ho\oplus\cH')$. Obviously, the functionals
  $x\mapsto\skp{\xi}{x\xi_{0}}$, $\xi\in\cH'$, $\xi_{0}\in\Ho$ form a
  total set in $\cE_{*}$. Due to the minimality of the Kolmogorov
  decomposition, they
    are approximated by functionals $x\mapsto\skp{y\eta_{0}}{x\xi_{0}}
    = \skp{\eta_{0}}{\skpo{y}{x}\xi_{0}}$, $\eta_{0},\xi_{0}\in\Ho$,
    $y\in\cE$. The latter can be used to approximate
    $x\mapsto\phi(\skpo{y}{x})$, $\phi\in{\Ao}_{*}$. Thus, we conclude
  $\cE_{*}=\mathrm{lh}\set{\phi(\skpo{y}{\cdot})}{y\in\cE,
    \phi\in\cA_{*}}^{\text{--}\nmt}$.
  
  The self duality of $\cE$ is proved with the help of the approximate
  unit $(e_{\al})_{\al\in I}$ of $\cK(\cE)$. For $\Psi\in\cE'$ we
  check
  \begin{align*}
        \Psi(e_{\al}x) 
   & =  \sum_{i=1}^{n_{\al}}
        \Psi(y_{i}^{\al})\,\skpo{y_{i}^{\al}}{x}  
     =  \bskpo{\sum_{i=1}^{n_{\al}}
        y_{i}^{\al}\,\Psi(y_{i}^{\al})^{*}}{x}
     =: \skpo{z_{\al}}{x}\,.
  \end{align*}
  This expression converges in norm to $\Psi(x)$. Thus
  $\phi(\Psi(e_{\al}x))$ converges to $\phi(\Psi(x))$ for any
  $\phi\in{\Ao}_{*}$. From this we conclude that
  $(z_{\al})_{\al\in I}$ converges to some $z\in\cE$ in the weak*
  topology on $\cE$, since the net is bounded: $\nmo{z_{\al}} =
  \nm{\skpo{z_{\al}}{\cdot}} = \nm{\Psi\circ e_{\al}} \le
  \nm{\Psi}\,\nmo{e_{\al}} \le \nm{\Psi}$. Therefore we get
  $\phi(\Psi(x))=\phi(\skpo{z}{x})$ for any $\phi\in{\Ao}_{*}$ and any
  $x\in\cE$, i.e., $\Psi=\skpo{z}{\cdot}$.\qedopt
\end{proof}

\begin{CorApp}\label{BE-LE}
  For a Hilbert W*-module $\cE$ is $\cB(\cE)=\cL(\cE)$.
\end{CorApp}

\begin{proof}
  For $T\in\cB(\cE)$ and $y\in\cE$ defines $\cE\ni
  x\mapsto\skpo{y}{Tx}$ an element in $\cE'$, i.e., there exists a
  unique element $z_{y}\in\cE$ such that
  $\skpo{z_{y}}{x}=\skpo{y}{Tx}\,$ for any $x\in\cE$. From
  \cite{Lance1} it is known that $T^{*}:y\mapsto z_{y}$ is the adjoint
  of $T$. This proves $\cB(\cE)\subseteq\cL(\cE)$. The inverse
  inclusion is obvious.\qedopt
\end{proof}

\begin{CorApp}\label{yinE}
  Let $\cE\subseteq\cB(\Ho\oplus\cH')$ be a Hilbert W*-module over the
  von Neumann algebra $\Ao\subseteq\cB(\Ho)$. Then
  $y\in\cB(\Ho\oplus\cH')$ and $y^*x\in\Ao\,$ for any $x\in\cE$ implies
  $y\in\cE$.
\end{CorApp}

\begin{proof}
  $\cE\ni x\mapsto y^*x$ defines an element in $\cE'$, i.e., it is
  $y^*x\xi_{0}=\skpo{z}{x}\xi_{0}=z^*x\xi_{0}$ for any $x\in\cE$,
  $\xi_{0}\in\Ho$ and some element $z\in\cE$. Since $\cH'$ is generated
  by $\cE\Ho$ we conclude $y^*=z^*$, and thus $y\in\cE$. \qedopt
\end{proof}

\begin{ThmApp}[Kaplansky]\label{Kaplansky}
  Let $\cE$ be the Hilbert W*-module generated by the pre-Hilbert
  module $\cE_{0}$. Then the unit ball $\cE_{0,1}$ of $\cE_{0}$ is
  $\sigma$-stop dense in the unit ball $\cE_{1}$ of $\cE$.
\end{ThmApp}

\begin{proof}
  The so-called linking algebra $\cL_{\cE}:=
  \begin{bmatrix} \Ao & \cE^{*} \\ \cE & \cL(\cE)\end{bmatrix}$ is
  a W*-algebra. Here we have $\cE^{*}:=\set{x^{*}}{x\in\cE}$. The
  density follows from the contractivity of the embedding of $\cE$ in
  $\cL_{\cE}$ and Kaplansky's density theorem for $\cL_{\cE}$. \qedopt
\end{proof}

\section{The $\psi$-adjoint of morphisms}
\label{sec-appendix-psiadjoint}

We provide results which are needed in the proof of Theorem
\ref{lap-morph}.

\begin{LemApp}\label{morph-cont}
  Let $(\cA,\psi)$ be a probability space and $T\colon \cA\to\cA$ a
  completely positive map such that there exists an $a\in\cA$ with the
  property
  \begin{align}\label{Tofone}
        \psi(T(x))
   & =  \psi(xa)
  \end{align}
  for all $x\in\cA$. Then $T$ is normal (i.e.
  weakly*-weakly*-continuous), $\si$stop-$\si$stop- and
  $\si$stop*-$\si$stop*-continuous.
\end{LemApp}

Note that the operator $a$ is necessary unique, if it exists. For the
proof of the lemma we introduce the so called $\psi$-norm on $\cA$ by
$\nm{x}_{\psi}:=\psi(x^{*}x)^{1/2}$.

\begin{proof}
  Using \cite[Prop. 1.24.1]{Sakai}, we obtain the inequality
  \begin{align}\label{T-Sakai}
        \psi(T(x))
   &\le \nm{a}\psi(x)\,,
        \quad x\in\cA^{+}
  \end{align}
  from the positivity of $\psi\circ T$. The Kadison-Schwarz inequality
  for $T$ leads to $\nm{T(x)}_{\psi} \le \nm{a}^{1/2}\nm{x}_{\psi}$.
  Now the Cauchy-Schwarz inequality shows for each $x\in\cA$ the
  $\si$\,stop-continuity of $\cA\ni y\mapsto\psi_{x}(T(y))$, where
  $\psi_{x}(y):=\psi(xy)$. Since $\set{\psi_{x}}{x\in\cA}$ is uniformly
  dense in $\cA_{*}$ by the bipolar theorem, the usual $\ve/2$-argument
  shows the stop-continuity of the maps $y\mapsto\phi(T(y))$,
  $\phi\in\cA_{*}$, on bounded subsets of $\cA$. From \cite[Thm.
  II.2.6]{Take03a} we conclude $\phi\circ T\in\cA_{*}$ for any
  $\phi\in\cA_{*}$, i.e. $T$ is normal. Now let $(x_{\al})_{\al\in I}$
  be a net converging $\si$-strongly to zero. Then for any positive
  $\phi\in\cA_{*}$ we have by Kadison-Schwarz's inequality
  $\phi(T(x_{\al}^{*})T(x_{\al}^{}))\le \phi\circ
  T(x_{\al}^{*}x_{\al}^{})\xrightarrow{\al}0$, i.e.
  $T(x_{\al}^{})\xrightarrow{\al}0$ $\si$-strongly. The
  $\si$\,stop*-$\si$\,stop*-continuity of $T$ is shown
  analogously. \qedopt 
\end{proof}

For the readers convenience we include a direct proof of a result from
\cite{Kuem84aUP}, which is needed in Theorem \ref{lap-morph} (cf.\
\cite{Hell01a}). In the following $\De$ and $J$ denotes, as usual,
the modular operator and the modular conjugation with respect to $\psi
= \skp{\Om}{\cdot\,\Om}$.

\begin{ThmApp}\label{psi-adj}
  For a completely positive map $T$ on the probability space
  $(\cA,\psi)$ with property \eqref{Tofone} the following are
  equivalent:
  \begin{enumerate}
  \item\label{psi-adj-i} 
    There exists a completely positive map $T^{*}$ on $\cA$ with the 
    property $\psi(xT(y))=\psi(T^{*}(y)x)$ f.a.\ $x,y\in\cA$.
  \item\label{psi-adj-ii} 
    $T$ commutes with the modular automorphism group $\si^{\psi}$ 
    of $\psi$.
  \end{enumerate}
  Consequently, the operator $a$ in equation \eqref{Tofone} belongs to
  the centralizer $\cA^{\psi}$ of $\psi$ and is given by $T^{*}(\1)$.
\end{ThmApp}
\begin{Defn}
$T^{*}$ is called the \emph{$\psi$-adjoint} of $T$.
\end{Defn}
\begin{proof}
  By \eqref{T-Sakai}, $T$ has an extension $\overline{T}$ to a bounded
  operator on the GNS Hilbert space $\cH_{\psi}$, which is defined by
  $\overline{T}x\Om := T(x)\Om$, $x\in\cA$.
  
  (\ref{psi-adj-ii}) $\Rightarrow$ (\ref{psi-adj-i}): Since $T$ and
  $\si^{\psi}$ commute, $T$ maps the $\stp$-dense subspace of entire
  analytic elements $\cA_{\text{a}}$ for $\si^{\psi}$ into itself.
  Hence, for any $y\in\cA_{\text{a}}$, the analytic continuation of
  $\Rset\ni t\mapsto \overline{T}\De^{\i t}y^{*}\Om =
  \De^{\i t}T(y^{*})\Om$, evaluated at $t=-\i/2$, yields
  \[
        \overline{T}Jy\Om
     =  \overline{T}\De^{1/2}y^{*}\Om
     =  \De^{1/2}_{}\overline{T}y^{*}\Om
     =  \De^{1/2}_{}T(y^{*})\Om
     =  JT(y)\Om
     =  J\overline{T}y\Om\,.
  \]
  It follows $[\overline{T},J]=0$, since $\cA_{\text{a}}\Om$ is dense
  in $\cH_{\psi}$. Next we will show $\overline{T}^{*}\cA^{+}\Om
  \subseteq \cA^{+}\Om$. To begin with, by the duality of the cones
  $\overline{\cA^{+}\Om}$ and $\overline{\cA'^{+}\Om}$ and the
  estimation
  \[
        \skp{Jx\Om}{\overline{T}^{*}y\Om}
     =  \skp{JT(x)\Om}{y\Om}
     =  \skp{T(x)^{1/2}\Om}{JyJT(x)^{1/2}\Om}
    \ge 0\,,
  \]
  for all $x,y\in\cA^{+}$, we obtain
  $\overline{T}^{*}y\Om\in\overline{\cA^{+}\Om}$. Next we observe that
  the positive linear functional
  $x\mapsto\skp{Jx\Om}{\overline{T}^{*}y\Om}$ is dominated by $\psi$:
  For $x\in\cA^{+}$ we have
  \[
        \skp{T(x)^{1/2}\Om}{JyJT(x)^{1/2}\Om}
    \le \nm{y}\,\skp{\Om}{T(x)\Om}
    \stackrel{\eqref{T-Sakai}}{\le}
        \nm{y}\,\nm{a}\,\psi(x)\,.
  \]
  Hence, by the Radon-Nikodym theorem, there is a unique $z\in\cA^{+}$
  with the property $\skp{Jx\Om}{\overline{T}^{*}y\Om} =
  \skp{Jz\Om}{x\Om} = \skp{Jx\Om}{z\Om}$ for any $x\in\cA$. That is,
  $\overline{T}^{*}y\Om=
  z\Om\in\cA^{+}\Om$.\\
  Therefore, $T^{*}(y)\Om:=\overline{T}^{*}y\Om$, $y\in\cA$, defines a
  positive linear map $T^{*}$ on $\cA$, which fulfills obviously 
  $\psi(xT(y))=\psi(T^{*}(x)y)$. Thus the uniqueness of $T^{*}$ and
  $a=T^{*}(\1)\ge0$ are evident. $T^{*}(\1)\in\cA^{\psi}$ follows from
  \[
        \psi(T^{*}(\1)x)
     =  \overline{\psi(T(x^{*}))}
     =  \overline{\psi(T^{*}(\1)x^{*})}
     =  \psi(xT^{*}(\1))\,.
  \]
  Up to now we  have used only the positivity of $T$. We are left to
  show that $T^{*}$ is completely positive. To this end
  consider the map $T_{(n)}:=\id_{n}\otimes T$ on the probability space
  $(M_{n}\otimes\cA,\psi_{(n)})$, where $\psi_{(n)} :=
  \tau_{n}\otimes\psi$, with $\tau_{n}$ the normed trace on $M_{n}$.
  The modular operator $\De_{(n)}$ and the modular conjugation
  $J_{(n)}$ of $\psi_{(n)}$ are, respectively, given by $\De_{(n)}: =
  \1_{n}\otimes\De$ and $J_{(n)}: = J_{n}\otimes J$, with $J_{n}$ the
  modular conjugation of $\tau_{n}$. Obviously, $T_{(n)}$ commutes with
  the modular automorphism group of $\psi_{(n)}$, given by
  $\id_{n}\otimes\si^{\psi}$. An elementary calculation shows
  $\psi_{(n)}(xT_{(n)}(y)) = \psi_{(n)}(T^{*}_{(n)}(x)y)$ for all $x,
  y\in M_{n}\otimes\cA$, with $T^{*}_{(n)} := \id_{n}\otimes T^{*}$.
  Choosing $x=\1_{n}\otimes\1$, we arrive at \eqref{Tofone} for
  $\psi_{(n)}$ and $T_{(n)}$, with $a$ replaced by $\1_{n}\otimes a$.
  Since $T$ is completely positive, $T_{(n)}$ is positive. Hence, our
  considerations above show the existence of the adjoint map
  $(T_{(n)})^{*}$, which by uniqueness is given by $T^{*}_{(n)}$. Since
  $(T_{(n)})^{*}$ is positive, $T^{*}$ is completely positive.
  
  (\ref{psi-adj-i}) $\Rightarrow$ (\ref{psi-adj-ii}): Let $y\in
  \operatorname{Dom}(\De)$ and $x\in\cA$. Then we have, using
  $\overline{T}^{*} = \overline{T^{*}}$:
  \begin{align*}
       \skp{x\Om}{\overline{T}\De y\Om}
   & = \skp{\De^{1/2}T^{*}(x)\Om}{\De^{1/2}y\Om}
     = \skp{JT^{*}(x^{*})\Om}{Jy^{*}\Om}                             \\
   & = \skp{y^{*}\Om}{T^{*}(x^{*})\Om}
     = \skp{T(y^{*})\Om}{x^{*}\Om}                                   \\
   & = \skp{J\De^{1/2}T(y)\Om}{J\De^{1/2}x\Om}
     = \skp{\De^{1/2}x\Om}{\De^{1/2}\overline{T}y\Om}\,.
  \end{align*}
  Since $\cA\Om$ is a core for $\De^{1/2}$, $\De^{1/2}\overline{T}y\Om$
  is in the domain of $\De^{1/2}$ and we have
  $\De\overline{T}y\Om=\overline{T}\De y\Om$. Hence $\overline{T}$ and
  $\De$ commute strongly. From this we obtain
  $[\overline{T},\De^{\i t}]=0$ and finally
  $T\circ\si_{t}^{\psi}=\si_{t}^{\psi}\circ T$.\qedopt
 \end{proof}



\bibliographystyle{alpha}                 
\label{section:bibliography}
\bibliography{ncbs-lit}                   


\end{document}